\documentclass[11pt]{article}
\usepackage{amsmath,amsbsy,amsfonts}
\usepackage[frenchb]{babel}
\usepackage{makeidx}

 \newtheorem{thm}{Th\'eor\`eme}[section]
 
 \newtheorem{lem}[thm]{Lemme}
 \newtheorem{prop}[thm]{Proposition}
 \newtheorem{rem}[thm]{Remarque}

 \title{Image des op\'erateurs d'entrelacements normalis\'es et p\^oles des s\'eries d'Eisenstein}
 \author{C. M{\oe}glin\\
 Institut de Math\'ematiques de Jussieu\\
 CNRS
 }
 \date{}
 
\begin{document}

 \maketitle
 \begin{abstract}
 On a montr\'e en \cite{holomorphie} comment normaliser les op\'erateurs d'entrelacement standard de fa\c{c}on \`a assurer l'holomorphie dans le domaine positif. L'objet de ce travail est  de d\'ecrire l'image des ces op\'erateurs et on montre que cette image est soit nulle soit irr\'eductible. Dans les cas ayant des applications globales on d\'ecrit l'image en termes de param\`etres dans les paquets d'Arthur. On termine l'article en donnant explicitemment les r\'esidus de s\'eries d'Eisenstein dans le cas d'un parabolique maximal et d'une induite de repr\'esentation automorphes de carr\'e int\'egrable modulo le centre. Ce travail admet les r\'esultats annonc\'es par Arthur sur la classification en paquet des repr\'esentations automorphes de carr\'e int\'egrable, mais pas les formules de multiplicit\'es fines.
 \end{abstract}

\section{Introduction}
Ici $G$ est un groupe classique, disons $Sp(2n)$ ou $SO(2n+1)$ et nos m\'ethodes s'appliquent aussi \`a $O(2n)$ (mais il faut faire quelques v\'erifications  li\'ees \`a la non connexit\'e) et aux groupes unitaires avec des changements de notations pour se mettre dans le cadre de \cite{unitaire}. On note $G^*$ le groupe classique alg\'ebrique dual de $G$; $G^*$ est vu comme un sous-groupe complexe d'un certain $GL(m^*_{G},{\mathbb C})$ (ce qui d\'efinit $m^*_{G}$) via sa repr\'esentation naturelle. On fixe aussi $F$ un corps p-adique et $k$ un corps de nombres. Dans la partie locale, on suppose que $G$ est d\'efini sur $F$ et on confondra $G$ et $G(F)$.

 On fixe $\psi$ une repr\'esentation unitaire irr\'eductible de $W_{F}\times SL(2,{\mathbb C})\times SL(2,{\mathbb C})$ et on suppose que $\psi$ est de dimension $m^*_{G}$ et se factorise par $G^*$. Gr\^ace \`a la correspondance de Langlands local pour les groupes $GL$, on sait associer \`a $\psi$ une repr\'esentation unitaire de $GL(m^*_{G},F)$, not\'ee $\pi^{GL}(\psi)$.

 On fixe aussi une repr\'esentation $\rho$ cuspidale unitaire irr\'eductible d'un groupe lin\'eaire $GL(d_{\rho},F)$ (cela d\'efinit $d_{\rho}$) et un entier $a_{0}$; cela permet de former la repr\'esentation de Steinberg $St(\rho,a_{0})$ du groupe $GL(d_{\rho}a,F)$. Au groupe classique $G$ et \`a  tout entier $d$,  on associe une repr\'esentation $r_{G,d}$ (en g\'en\'eral on enl\`eve le $d$) de $GL(d,{\mathbb C})$ qui est $\wedge^2 {\mathbb C}^{d}$ si $G$ est un groupe de type $SO(2n+1)$ et $Sym^2 {\mathbb C}^{d}$ si $G$ est de type $Sp(2n)$ ou $O(2n)$.  
 
 On  d\'efinit la fonction m\'eromorphe: $$r(s,\rho,a_0,\psi):=
 \frac{L(St(\rho,a_{0})\times \pi^{GL}(\psi),s)}{L(St(\rho,a_{0})\times \pi^{GL}(\psi),s+1)}\times \frac{L(St(\rho,a_{0}),r_{G,ad_{\rho}},2s)}{L(St(\rho,a_{0}),r_{G,ad_{\rho}},2s+1)}.
 $$
 Cette fonction se calcule tr\`es facilement; on \'ecrit $\pi^{GL}(\psi)$ sous-forme d'induite:
 $$
 \pi^{GL}(\psi)=\times_{(\rho',a',b')\in Jord(\psi)}Speh(St(\rho',a'),b'),
 $$o\`u $\rho'$ est une repr\'esentation cuspidale unitaire d'un groupe $GL$ et $a',b'$ sont des entiers
 (ce qui d\'efinit aussi $Jord(\psi)$). Alors (cf. \cite{manuscripta} 2.1.1) $r(s,\rho,a_0,\psi):=$
 $$
 \prod_{(\rho',a',b')\in Jord(\psi)}\frac{L(St(\rho,a_{0})\times St(\rho',a'),s-(b'-1)/2)}{L(St(\rho,a_{0})\times St(\rho',a'),s+(b'+1)/2)} \times  \frac{L(St(\rho,a_{0}),r_{G,ad_{\rho}},2s)}{L(St(\rho,a_{0}),r_{G,ad_{\rho}},2s+1)}.
 $$
Il est clair que les d\'enominateurs de la fonction $r(s,\rho,a_0,\psi)$ n'ont ni z\'ero ni p\^ole pour $s\in {\mathbb R}_{>-1/2}$. Les num\'erateurs ont \'eventuellement des p\^oles que l'on rappellera ult\'erieurement.

On fixe $\pi$ dans le paquet d'Arthur associ\'e \`a $\psi$; on rappelle dans le texte ce que cela veut dire et ceci suppose donc que  l'on connaisse la classification des repr\'esentations cuspidales dans les termes de \cite{elementaire}; on a montr\'e en \cite{unitaire} les hypoth\`eses minimales qu'il fallait mais \'evidemment la classification de Langlands des s\'eries discr\`etes annonc\'ees par Arthur est plus que suffisante comme montr\'e en \cite{pourshahidi}. 

Et on consid\`ere l'op\'erateur d'entrelacement standard, bien d\'efini apr\`es un choix d'\'el\'ement du groupe de Weyl, diff\'erents choix \'etant reli\'es par la multiplication par une fonction holomorphe inversible:
 $$M(s,\rho,a_{0},\pi):\qquad
 St(\rho,a_{0})\vert\,\vert^{s}\times \pi \rightarrow St(\rho^*,a_{0})\vert\,\vert^{-s}\times \pi$$

On pose $N_{\psi}(s,\rho,a_{0},\pi):= r(s,\rho,a_0,\psi)^{-1}M(s,\rho,a_{0},\pi)$. Le r\'esultat principal de \cite{holomorphie}, dans le cas o\`u $\rho\simeq \rho^*$ (hypoth\`ese lev\'ee ici) est de montrer que $N_{\psi}(s,\rho,a_{0},\pi)$ est holomorphe pour tout $s\in {\mathbb R}_{\geq 0}$ et le premier but de cet article est de d\'ecrire l'image de cet op\'erateur.

On montre que soit cet op\'erateur, calcul\'e en un point $s=:s_{0}\in {\mathbb R}_{>0}$ est  identiquement 0 soit son image est une repr\'esentation irr\'eductible. De plus  on d\'ecrit alors cette image dans les cas o\`u il y a  des applications globales, c'est \`a dire le cas o\`u $s_{0}$ est un demi-entier donc de la forme $(b_{0}-1)/2$ avec $b_{0}\geq 2$ et $(\rho,a_{0},b_{0}-2)\in Jord(\psi)$ (si $b_{0}=2$, il n'y a pas de condition). On a dans ce cas un r\'esultat tr\`es pr\'ecis: on note $\psi^+$ le morphisme analogue \`a $\psi$ mais tel que $Jord(\psi^+)$ se d\'eduit de $Jord(\psi)$ en rempla\c{c}ant $(\rho,a_{0},b_{0}-2)$ par $(\rho,a_{0},b_{0})$ (si $b_{0}=2$ on ajoute $(\rho,a_{0},b_{0})$); \`a l'aide des param\`etres de $\pi$ comme repr\'esentation dans le paquet associ\'e \`a $\psi$, on d\'efinit des param\`etres pour le paquet associ\'e \`a $\psi^+$ et on note $\pi^+$ la repr\'esentation du paquet associ\'e \`a $\psi^+$; elle peut \^etre nulle. Et on montre l'\'egalit\'e $$Im\, N_{\psi}(s,\rho,a_{0},\pi)=\pi^+$$ ce qui veut dire que l'image de l'op\'erateur est identiquement 0 exactement quand $\pi^+$ est nulle et sinon cette image est la repr\'esentation irr\'eductible $\pi^+$. Bien s\^ur la construction des repr\'esentations \`a partir des param\`etres est compliqu\'ee et le r\'esultat n'est donc pas simple mais si on trouve ult\'erieurement une meillleure param\'etrisation on pourra facilement retraduire ce r\'esultat; c'est donc \`a mon avis le meilleur possible. Quand les articles d'Arthur seront compl\`etement disponibles, on pourra calculer explicitement le caract\`ere associ\'ee \`a cette repr\'esentation et intervenant dans les formules de multiplicit\'e.

Dans le texte on g\'en\'eralise en rempla\c{c}ant $St(\rho,a_{0})$ par une composante en la place p-adique fix\'ee d'une repr\'esentation cuspidale autoduale d'un groupe  $GL$ et on consid\`ere aussi la g\'en\'eralisation des paquets d'Arthur de fa\c{c}on \`a couvrir le cas o\`u la seule hypoth\`ese sur $\pi$ est d'\^etre une composante locale d'une forme automorphe de carr\'e int\'egrable; le point ici est d'\'eviter l'utilisation de la conjecture de Ramanujan et donc de consid\'erer certaines repr\'esentations non unitaires de $W_{F}\times SL(2,{\mathbb C}) \times SL(2,{\mathbb C})$.

Pour avoir des applications globales, 
on traite aussi le cas des places archim\'ediennes dans des cas tr\`es particuliers:  on se limite \`a $F={\mathbb R}$ et surtout on met l'hypoth\`ese forte que $\pi$ a de la cohomologie pour un bon syst\`eme de coefficient et que le caract\`ere infinit\'esimal de l'induite $\rho\vert\,\vert^{s_{0}}\times \pi$  est entier r\'egulier, c'est surtout r\'egulier qui compte. On d\'emontre alors le m\^eme r\'esultat que ci-dessus, \`a savoir que l'op\'erateur normalis\'e est holomorphe en $s=s_{0}$ et que son image est une repr\'esentation irr\'eductible explicite, en termes de param\`etre de Langlands. Ce sont les hypoth\`eses qui simplifient le r\'esultat.

On termine l'article en d\'ecrivant les points de non holomorphie des s\'eries d'Eisenstein dans le cadre suivant. Le corps de base est un corps de nombres $k$ (totalement r\'eel) et on suppose que $G$ est d\'efini sur $k$. On consid\`ere $H$ un groupe de m\^eme type que $G$ d\'efini sur  $k$. On fixe $P$ un sous-groupe parabolique maximal de $H$, d\'efini sur $k$, et on suppose que les sous-groupes de Levi de $P$ sont isomorphes \`a  $GL(d)\times G$ (\'eventuellement, ici, $G$ peut \^etre le groupe trivial). On fixe $\pi_{0}$ une repr\'esentation de carr\'e int\'egrable irr\'eductible de $G$ et $\tau$ une repr\'esentation cuspidale irr\'eductible et unitaire de $GL(d)$. Et on consid\`ere les s\'eries d'Eisenstein de la forme $E_{P}^{H}(\tau\vert\,\vert^{s}\times \pi_{0},f)$ o\`u $s\in {\mathbb C}$. On sait qu'une telle s\'erie d'Eisenstein est holomorphe pour $Re\, s=0$ et on s'int\'eresse au cas o\`u $Re\, s>0$; quitte \`a tordre $\tau$ par un caract\`ere unitaire, on suppose en fait que $s\in {\mathbb R}$. On suppose que les r\'esultats d'Arthur annonc\'es dans \cite{arthurnouveau} sont disponibles pour $\pi_{0}$; ainsi \`a $\pi_{0}$ est associ\'e un ensemble fini de couples $(\rho,b)$ o\`u $\rho$ est une repr\'esentation cuspidale d'un groupe de type $GL$ et $b$ est un entier; on note $Jord(\pi_{0})$ cet ensemble. La propri\'et\'e qui d\'etermine uniquement cet ensemble est qu'en notant $Speh(\rho,b)$ la repr\'esentation du GL convenable dans l'espace des r\'esidus:
$$Speh(\rho,b)=\biggl((\prod_{i\in [1,b[}(s_{i}-s_{i+1}))E(\rho\vert\,\vert^{(b-1)/2+s_{1}}\times \cdots \times \rho\vert\,\vert^{-(b-1)/2+s_{b}}, s_{1}, \cdots, s_{b})\biggr)_{s_{1}=0, \cdots, s_{b}=0}$$
la repr\'esentation $\pi^{GL}_{0}:==\times_{(\rho,b)\in Jord(\pi_{0})}Speh(\rho,b)$ et la repr\'esentations $\pi$ ont leurs composantes locales non ramifi\'ees qui se correspondent presque partout via la correspondance non ramifi\'ee de Langlands.
On d\'efinit l'analogue global de $r(s,\rho,a_{0},\psi)$ en posant
$$
r(s,\tau,\pi^{GL}_{0}):=\frac{L(\tau\times \pi^{GL}_{0},s)}{L(\tau\times \pi^{GL}_{0},s+1)}\frac{L(\tau,r_{G},2s)}{L(\tau,r_{G},2s+1)};
$$
ici tous les facteurs comptent.
On fixe un r\'eel positif $s_{0}$ et on suppose que le caract\`ere infinit\'esimal de $\tau\vert\,\vert^{s_{0}}\times \pi_{0}$ est entier et r\'egulier et que $\pi_{0}$ a de la cohomologie \`a l'infini. Gr\^ace \`a \cite{manuscripta} 1.2.2 et  4.4.4 (i) (avec le r\'esultat d'holomorphie de \cite{holomorphie} compl\'et\'e ici ), on sait que les s\'eries d'Eisenstein $E_{P}^H(\tau\vert\,\vert^{s}\times \pi)$ sont holomorphe en $s=s_{0}\in {\mathbb R}_{>0}$ sauf \'eventuellement si  $r(s,\tau,\pi^{GL}_{0})$ a un p\^ole n\'ecessairement d'ordre 1 en $s=s_{0}$; ceci s'explicite en:

\

{\sl soit $s_{0}=1/2$ et $L(\tau,r_{G,d},2s)$ a un p\^ole en $s=1/2$, soit $s_{0}$ est un demi-entier sup\'erieur ou \'egal \`a $1$ avec $(\tau,2s_{0}-2)\in Jord(\pi_{0})$ et, dans les 2 cas, pour tout $(\rho,b)\in Jord(\pi_{0})$ avec $b=2s_{0}-1$, $L(\tau\times \rho,1/2)\neq 0$. }

\

On donne dans cet article des conditions n\'ecessaires et suffisantes pour que $s_{0}$ soit un p\^ole
On  exprime ces conditions ainsi: on a globalement la repr\'esentation d'un groupe lin\'eaire convenable $\pi^{GL}_{0}:=\times_{(\rho,b)\in Jord(\pi_{0})}Speh (\rho,b)$; en toute place $v$, avec la param\'etrisation de Langlands \'etendu en une param\'etrisation d'Arthur, on a un morphisme $\psi_{0,v}$ qui correspond \`a la composante locale en la place $v$ de la repr\'esentation $\pi^{GL}_{0}$. On consid\`ere la composante locale en toute place $v$ de $\pi_{0}$, cela joue le r\^ole de la repr\'esentation $\pi$ dans ce qui pr\'ec\`ede; l'analogue de $St(\rho,a_{0})$ est (avec la g\'en\'eralisation d\'ej\`a annonc\'ee) la composante locale de $\tau_{0}$; comme expliqu\'e ci-dessus, pour $s_{0}$ v\'erifiant les conditions n\'ecessaires pour avoir non holomorphie des s\'eries d'Eisenstein, on a d\'efini en toute place $v$ une repr\'esentation $\pi_{0,v}^+$ qui \'eventuellement peut \^etre nulle (sauf aux places archim\'ediennes o\`u la non nullit\'e est assur\'ee avec  l'hypoth\`ese mise).  On pose formellement $\pi^+_{0}=\otimes_{v}\pi^+_{0,v}$.   On suppose aussi que l'on sait a priori que  la repr\'esentation $\pi^+_{0}$ si elle est non nulle intervient avec multiplicit\'e au plus 1 dans le spectre discret. On montre alors le th\'eor\`eme suivant:

\

{\sl
Les s\'eries d'Eisenstein $E(\tau\vert\,\vert^{s}\times \pi_{0},f)$ ont un p\^ole  en $s_{0}$  si et seulement si les conditions n\'ecessaires ci-dessus sont satisfaites et $\pi^+_{0}$ est non nulle et alors
$$ \biggl((s-s_{0})E(\tau\vert\,\vert^{s}\times \pi_{0},f)\biggr)_{s=s_{0}}=\pi^+_{0};$$
}

Ceci revient \`a dire que l'existence d'un p\^ole est \'equivalent \`a ce que l'op\'erateur d'entrelacement global:
$$
\tau\vert\,\vert^{s}\times \pi_{0}\rightarrow \tau^*\vert\,\vert^{-s}\times \pi_{0}
$$ait un p\^ole. Cela ne veut pas dire que les autres op\'erateurs d'entrelacement n'ont pas de p\^ole mais cela dit que s'ils ont des p\^oles celui \'ecrit en a aussi. Evidemment ce sont les hypoth\`eses de r\'egularit\'e du caract\`ere infinit\'esimal prises qui permettent un r\'esultat si simple.

L'hypoth\`ese de multiplicit\'e 1 est une hypoth\`ese raisonable car elle r\'esulte de la formule de multiplicit\'e globale d'Arthur et des formules de multiplicit\'e 1 locale; celles-ci sont d\'emontr\'es aux places p-adiques (\cite{pourshahidi}) d\`es qu'on les a pour les paquets de s\'eries discr\`etes (r\'esultat annonc\'e par Arthur) et aux places archim\'ediennes elles r\'esultent de \cite{aj} si on sait que $\pi^+_{0}$ a de la cohomologie pour un bon syst\`eme de coefficient et si on sait que les paquets d'Adams-Johnson de \cite{aj} sont bien ceux d'Arthur.

On montre aussi que si $\pi'$ est une repr\'esentation automorphe irr\'eductible de carr\'e int\'egrable qui n'est pas cuspidale, et si $\tau$ est une repr\'esentation cuspidal unitaire, $s_{0}\in 1/2{\mathbb N}$ sont tel que les termes constants de $\pi'$ pour un parabolique maximal de la forme $GL(d_{\tau})\times G'$ ($G'$ un groupe de m\^eme type que $G$ convenable) ont une projection non nulle sur l'espace de $\tau\vert\,\vert^{-s_{0}}$ vu comme repr\'esentation automorphe cuspidale du facteur $GL(d_{\tau})$, alors $\pi'\simeq \pi^+_{0}$ pour un bon choix de $\pi_{0}$. On discute en \ref{commentaires} pourquoi on s'est approch\'e sans l'atteindre de l'objectif de \cite{manuscripta}, \`a savoir donner des conditions n\'ecessaires et suffisantes pour d\'ecrire les repr\'esentations automorphes de carr\'e int\'egrable non cuspidales; en fait on \'ecrit en \ref{commentaires} de telles conditions n\'ecessaires et suffisantes mais l'une de ces conditions me semble redondante et devrait dispara\^{\i}tre quand on aura explicit\'e les formules de multiplicit\'e d'Arthur.

\

Les applications que l'on peut esp\'erer de ces r\'esultats concernent  la cohomologie des repr\'esentations automorphes; pour cela, il faut utiliser la description faite par Franke en \cite{franke98} des formes automorphes comme d\'eriv\'ees de s\'erie d'Eisenstein \`a partir de repr\'esentations de carr\'e int\'egrable et essayer de g\'en\'eraliser  la premi\`ere partie de \cite{franke}. Sans travail suppl\'ementaire, on ne peut esp\'erer de description explicite comme ce qui a \'et\'e fait dans certaines situations en particulier   r\'ecemment en  \cite{g}, \cite{hg}, \cite{gg} et \cite{gs}.

\

Ce travail a \'et\'e expos\'e lors de la p\'eriode sp\'eciale se d\'eroulant \`a l'Institut Erwin Schr\"odinger de Vienne d\'ebut 2009 et je remercie l'ESI pour son hospitalit\'e, les organisateurs de cette p\'eriode, G. Henniart, G. Muic et J. Schwermer ainsi que les auditeurs et en particulier H. Grobner pour ses remarques sur le cas archim\'edien.
\tableofcontents

\section{Notations et description des repr\'esentations $\pi^+$, cas local p-adique}
\subsection{Notations g\'en\'erales}

On se place sur un corps p-adique, $F$. On param\'etrise les repr\'esentations irr\'eductibles de $W_{F}\times SL(2,{\mathbb C})\times SL(2,{\mathbb C})$ par un triplet $(\rho,a,b)$ form\'ee d'une repr\'esentation irr\'eductible de $W_{F}$ et de 2 entiers qui param\'etrisent chacun une repr\'esentation irr\'eductible de $SL(2,{\mathbb C})$. On note $m^*_{G}$ la dimension de la repr\'esentation naturelle du $L$-groupe de $G$; cette notation servira peu.

\subsubsection{Bonne Parit\'e\label{bonneparite}}

On fixe $G$ un groupe classique d\'efini sur $F$; ici on se limite aux groupes orthgonaux ou symplectique. On dit qu'une repr\'esentation irr\'eductible de $W_{F}\times SL(2,{\mathbb C})\times SL(2,{\mathbb C})$ se factorise par un groupe de m\^eme type que le groupe dual de $G$, si  elle est autoduale et symplectique si $G$ est un groupe orthogonal impair et orthogonal sinon. On dit alors que cette repr\'esentation a bonne parit\'e. Soit $\psi$ une repr\'esentation semi-simple de $W_{F}\times SL(2,{\mathbb C})\times SL(2,{\mathbb C})$ de dimension finie dont toutes les sous-repr\'esentations irr\'eductibles sont de bonne parit\'e, on dit alors que $\psi$ est de bonne parit\'e.

\subsubsection{Param\`etre g\'eom\'etrique des paquets d'Arthur et $Jord(\psi)$\label{jord}}

Les param\`etres g\'eom\'etriques des paquets d'Arthur sont des  repr\'esentations irr\'eductibles de $W_{F}\times SL(2,{\mathbb C})\times SL(2,{\mathbb C})$ que l'on suppose semi-simple, se factorisant par le groupe dual de $G$. On suppose que ces repr\'esentations sont continues quand on les restreint \`a $W_{F}$ et alg\'ebriques sur les 2 copies de $SL(2,{\mathbb C})$. En g\'en\'eral, c'est le cas le plus difficile, on suppose ces repr\'esentations unitaires quand on les restreint \`a $W_{F}$; on dit simplement unitaire. On fera les g\'en\'eralisations n\'ecessaires pour \'eviter la conjecture de Ramanujan en \ref{descriptiongenerale}. En \ref{descriptionpaquets} on expliquera comment on associe \`a un tel morphisme un ensemble fini de repr\'esentations lisses irr\'eductibles de $G$. On note $Jord(\psi)$ l'ensemble des sous-repr\'esentations irr\'eductibles incluses dans $\psi$ en tenant compte des multiplicit\'es. Avec les notations d\'ej\`a introduites, on a donc:
$$
\sum_{(\rho,a,b)\in Jord(\psi)}abdim\,\rho=m^*_{G}
$$

\subsubsection{Correspondance de Langlands pour les groupes $GL$\label{correspondance}}
On utilisera librement la correspondance de Langlands pour $GL$ d\'emontr\'ee par Harris-Taylor en \cite{harris} et par Henniart en \cite{henniart}. Cela permet d'identifier toute repr\'esentation continue irr\'eductible de dimension finie $d$ de $W_{F}$, not\'ee $\rho$, \`a une repr\'esentation cuspidale irr\'eductible de $GL(d,F)$, not\'ee encore $\rho$. Soit $\psi$ comme ci-dessus, on note $\pi^{GL}(\psi)$ la repr\'esentation induite $$\pi^{GL}(\psi)=\times_{(\rho,a,b)\in Jord(\psi)}Speh(St(\rho,a),b),$$ o\`u $St(\rho,a)$ est la repr\'esentation de Steinberg, unique sous-repr\'esentation irr\'eductible de l'induite $\rho\vert\,\vert^{(a-1)/2}\times \cdots \times \rho\vert\,\vert^{-(a-1)/2}$ et o\`u $Speh(St(\rho,a),b)$ est l'unique quotient irr\'eductible de l'induite $St(\rho,a)\vert\,\vert^{(b-1)/2}\times \cdots \times St(\rho,a)\vert\,\vert^{-(b-1)/2}$. Si $\psi$ se factorise par le groupe dual de $G$, cette repr\'esentation $\pi^{GL}(\psi)$ est, \`a fortiori, autoduale.

Soit $\rho$ une repr\'esentation cuspidale d'un groupe $GL(d,F)$ (ce qui d\'efinit $d$) et soit $[x,y]$ un segment, c'est-\`a-dire $x,y\in {\mathbb R}$ avec $x-y \in {\mathbb Z}$. Dans le texte $x,y$ seront toujours des demi-entiers mais cela ne sert pas ici. On note $<x,\cdots,y>_{\rho}$ l'unique sous-repr\'esentation irr\'eductible de $GL(d\vert x-y\vert +1,F)$ incluse dans l'induite 
$\rho\vert\,\vert^{x}\times \cdots \times \rho\vert\,\vert^{y}$. Cette repr\'esentation est une s\'erie discr\`ete tordue si $x-y\geq 0$ et est une repr\'esentation de Speh si $x-y\leq 0$. On g\'en\'eralise parfois cette notation dans la situation suivante: soit une matrice dont les lignes et les colonnes sont des segments mais de croissance oppos\'ee; en d'autres termes soit $x,y,z$ tel que $[x,y]$ soit un segment et $[x,z]$ soit aussi un segment mais avec $(x-y)(x-z)\leq 0$ et on consid\`ere la matrice
$$A:=
\begin{matrix} x&\cdots &y\\
\vdots &\vdots &\vdots\\
z&\cdots &z+y-x
\end{matrix}$$
et la repr\'esentation $<A>_{\rho}$ est alors l'unique sous-repr\'esentation irr\'eductible incluse dans l'induite $<x,y>_{\rho}\times \cdots \times <z,z+y-x>_{\rho}$. La possibilit\'e de d\'efinir ainsi une unique repr\'esentation irr\'eductible r\'esulte des r\'esultats de Bernstein-Zelevinsky et Zelevinsky que nous utilisons librement dans l'article.
\subsubsection{Op\'erateurs d'entrelacement\label{entrelacement}}
Soit $\tau$ une repr\'esentation, en g\'en\'eral irr\'eductible mais ce n'est pas indispensable pour cette d\'efinition, d'un groupe $GL(d_{\tau})$ et soit $\pi$ une repr\'esentation irr\'eductible de $G$. Pour $s\in {\mathbb C}$, on consid\`ere la repr\'esentation induite $\tau\vert\,\vert^{s}\times \pi$ et l'op\'erateur d'entrelacement standard:
$$
M(s,\tau,\pi): \tau\vert\,\vert^{s}\times \pi\rightarrow \tau^*\vert\,\vert^{-s}\times \pi.
$$
Supposons que $\pi$ est dans un paquet d'Arthur associ\'e \`a une repr\'esentation $\psi$ de $W_{F}\times SL(2,{\mathbb C})\times SL(2,{\mathbb C})$ (cf. \ref{descriptionpaquets} et \ref{descriptiongenerale}) et on reprend la notation $r(s,\tau,\psi)$ de l'introduction:
$$
r(s,\tau,\psi)=\frac{L(\tau\times \pi^{GL}(\psi),s)}{L(\tau\times \pi^{GL}(\psi),s+1)}\frac{L(\tau,r_{G,d_{\tau}},2s)}{L(\tau,r_{G,d_{\tau}},2s+1)}
$$
$r_{G,d_{\tau}}$ a \'et\'e d\'efini dans l'introduction et ne compte pas quand on est dans le cas local, le facteur correspondant n'a ni z\'ero ni p\^ole pour $\tau$ temp\'er\'ee et $Re\, s>0$. On pose
$$
N_{\psi}(s,\tau,\pi):=r(s,\tau,\pi)^{-1}M(s,\tau,\pi).
$$
Quand $\tau$ est une repr\'esentation de Steinberg, c'est-\`a-dire de la forme $St(\rho,a)$, on pose
$$
N_{\psi}(s,St(\rho,a),\pi)=:N_{\psi}(s,\rho,a,\pi).
$$
\subsection{Rappel sur les paquets d'Arthur\label{descriptionpaquets}}
On fixe $\psi$ une repr\'esentation unitaire de $W_{F}\times SL(2,{\mathbb C})\times SL(2,{\mathbb C})$ \`a valeurs dans $GL(m^*_{G},{\mathbb C})$ comme dans \ref{jord}, d'o\`u aussi $\pi^{GL}(\psi)$ comme dans \ref{correspondance}.

Supposons momentan\'ement que $G$ est quasid\'eploy\'e. 
Arthur a annonc\'e dans le dernier chapitre de \cite{arthurnouveau}  qu'il existe un ensemble fini de repr\'esentation irr\'eductible de $G$, not\'e $\Pi(\psi)$ uniquement d\'etermin\'e par le fait que pour tout $\pi\in \Pi(\psi)$ il existe une nombre complexe non nul $a_{\pi}$ tel que la distribution $\sum_{\pi\in \Pi(\psi)}a_{\pi}tr\, \pi$ se transf\`ere, via l'endoscopie  la $\theta$-trace  de $\pi^{GL}(\psi)$.  Ne supposons plus $G$ quasid\'eploy\'e, on a alors une d\'efinition de $\Pi(\psi)$ en utilisant le transfert \`a la forme quasid\'eploy\'ee du groupe. 

En admettant ce r\'esultat d'Arthur pour les morphismes $\psi$ triviaux sur la 2e copie de $SL(2,{\mathbb C})$ et tel que $Jord(\psi)$ soit sans multiplicit\'e, ce sont ceux qui param\'etrisent les paquets de s\'eries discr\`etes, on a  retrouv\'e le r\'esultat g\'en\'eral d'Arthur en \cite{discret}; ce qui nous importe est que dans ces r\'ef\'erences, on a donn\'e une construction combinatoire des \'el\'ements de $\Pi(\psi)$ \`a partir des s\'eries discr\`etes. C'est cette description que l'on va rappeler.

Cette description est simple dans le cas o\`u la restriction de $\psi$ \`a $W_{F}$ fois la diagonale de $SL(2,{\mathbb C})\times SL(2,{\mathbb C})$ est sans multiplicit\'e. Dans ce cas, on a montr\'e en \cite{discret} que $\Pi(\psi)$ est en bijection avec les couples de fonctions $(\underline{t},\underline{\eta})$ de $Jord(\psi)$ \`a valeurs dans ${\mathbb N}\times \{\pm 1\}$ satisfaisant \`a
$$
\forall (\rho,a,b)\in Jord(\psi), \underline{t}(\rho,a,b)\in [0, [inf(a,b)/2]]$$
$$\hbox{  si $\underline{t}(\rho,a,b)=inf(a,b)/2$, alors $\underline{\eta}(\rho,a,b)=+$};\eqno(1)
$$
$$
\times_{(\rho,a,b)\in Jord(\psi)}\underline{\eta}(\rho,a,b)^{inf(a,b)}(-1)^{[inf(a,b)/2]+\underline{t}(\rho,a,b)}=
\epsilon_{G},\eqno(2)
$$
o\`u $\epsilon_{G}$ vaut $+$ quand $G$ est quasid\'eploy\'e. En g\'en\'eral, il vaut mieux d\'efinir $\epsilon_{G}$  comme valant l'invariant de Hasse de la forme bilin\'eaire servant \`a d\'efinir $G$; comme le m\^eme groupe peut correspondre \`a 2 formes bilin\'eaires avec des invariants de Hasse diff\'erents, dans ces cas, la param\'etrisation d\'epend de la forme bilin\'eaire et non du groupe et la notation est incorrecte.

Dans le cas particulier, o\`u en plus des hypoth\`eses d\'ej\`a faites, pour tout $(\rho,a,b)\in Jord(\psi)$, on a $a\geq b$, la param\'etrisation que nous avons d\'ecrites est particuli\`erement agr\'eable car elle donne directement les param\`etres de Langlands des repr\'esentations cherch\'ees. Si \`a l'inverse pour tout $(\rho,a,b)\in Jord(\psi)$, on a $b\geq a$, la param\'etrisation donn\'ee est une g\'en\'eralisation de la param\'etrisation de Zelevinsky. En g\'en\'eral, notre param\'etrisation est une interpolation des param\'etrisations de Langlands et de Zelevinsky adapt\'ee \`a $\psi$ et il n'y a pas de formule simple pour retrouver la param\'etrisation de Langlands. On n'aura pas besoin du d\'etail pr\'ecis de la construction ici, ils sont r\'esum\'es en \cite{holomorphie} 2.2.

\

Dans cet article, on a besoin d'utiliser pr\'ecis\'ement le passage du cas g\'en\'eral \`a ce cas particulier. On note $\psi_{bp}$ la somme des sous-repr\'esentations irr\'eductibles de $\psi$ ayant bonne parit\'e (cf. \ref{bonneparite}) et $\psi_{mp}$ la somme des autres. Comme $\psi$ se factorise par le groupe dual de $G$, on peut d\'ecouper (de fa\c{c}on non unique) $\psi_{mp}$ en la somme directe de 2 sous-repr\'esentations $\psi_{1/2,mp}\oplus \psi_{-1/2,mp}$ o\`u $\psi_{-1/2,mp}\simeq \psi_{1/2,mp}^*$. A $\psi_{1/2,mp}$ on associe comme pr\'ec\'edemment une repr\'esentation de $GL(d_{1/2,mp},F)$ (o\`u $d_{1/2,mp} $ est la dimension de la repr\'esentation $\psi_{1/2,mp}$); on note $\pi^{GL}(\psi_{1/2,mp})$ cette repr\'esentation. On va construire ci-dessous $\Pi(\psi_{bp})$ l'ensemble des repr\'esentations associ\'ees au morphisme $\psi_{bp}$ et on a montr\'e en \cite{pourshahidi} 3.2 que pour tout $\pi'\in \Pi(\psi_{bp})$ la repr\'esentation induite $\pi^{GL}(\psi_{1/2,mp})\times \pi'$ de $G$ est irr\'eductible et que $\Pi(\psi)$ est pr\'ecis\'ement l'ensemble de toutes ces repr\'esentations, l'induction d\'efinit donc une bijection de $\Pi(\psi_{bp})$ sur $\Pi(\psi)$. Il reste \`a rappeler la description de $\Pi(\psi_{bp})$. Pour cela on suit \cite{holomorphie} 2.8.

Pour simplifier les notations on suppose que $\psi=\psi_{bp}$. La param\'etrisation de $\Pi(\psi)$ se fait ici encore \`a l'aide des couples $\underline{t},\underline{\eta}$ v\'erifiant les m\^emes hypoth\`eses que ci-dessus mais d\'efinis sur l'ensemble $Jord(\psi)$ vu comme ensemble avec r\'ep\'etition (correspondant \`a la multiplicit\'e) mais il faut en plus un ordre total sur $Jord(\psi)$. Evidemment $\Pi(\psi)$ comme ensemble ne d\'epend pas du choix de l'ordre mais la param\'etrisation que nous ne donnons en d\'epend en g\'en\'eral.

On fixe donc un ordre total sur $Jord(\psi)$ et on fait en plus un choix de signe pour tout $(\rho,a,b)\in Jord(\psi)$ tel que $a=b$; pour avoir des notations plus simples on inclut ce choix de signe dans les notations en rempla\c{c}ant les triplets $(\rho,a,b)\in Jord(\psi)$ par des quadruplets $(\rho,A,B,\zeta)$ o\`u le passage se fait par les \'egalit\'es:
$$
A=(a+b)/2-1; B=\vert (a-b)\vert/2; \zeta(a-b)\geq 0.
$$
L'ordre total sur $Jord(\psi)$ doit v\'erifier la propri\'et\'e suivante:

\

${\mathcal P}:$ $
\forall (\rho,A,B,\zeta),(\rho',A',B',\zeta') \in Jord(\psi)$, les propri\'et\'es $\rho\simeq \rho'$, $\zeta=\zeta'$, $A>A'$ et $B>B'$ entra\^{\i}nent $(\rho,A,B,\zeta)>(\rho',A',B',\zeta')$.

\

En particulier l'ordre s\'epare les \'el\'ements de $Jord(\psi)$ qui sont \'egaux. Une fois l'ordre fix\'e, pour $G'$ un groupe de m\^eme type que $G$ mais de rang plus grand et pour $\psi'$ un morphisme analogue \`a $\psi$ mais relativement \`a $G'$, on dit que $\psi'$ domine $\psi$ s'il existe une fonction $\underline{T}: Jord(\psi)\rightarrow {\mathbb Z}_{\geq 0}$ respectant l'ordre sur $Jord(\psi)$ et tel que 
$$Jord(\psi')=\{(\rho,A+\underline{T}(\rho,A,B,\zeta),B+\underline{T}(\rho,A,B,\zeta),\zeta); (\rho,A,B,\zeta)\in Jord(\psi)\}.\eqno(3)$$
Evidemment quand $\psi'$ est connu on trouve la fonction $\underline{T}$, en posant $\underline{T}(\rho,A,B,\zeta)=A'-A$ pour $A'$ l'\'el\'em\'ent du quadruplet $(\rho',A',B',\zeta')$ situ\'e dans $Jord(\psi')$ \`a la m\^eme place que $(\rho,A,B,\zeta)$ dans $Jord(\psi)$ et n\'ecessairement pour ce quadruplet $\rho'=\rho,\zeta'=\zeta$ et $B'=B+A'-A$.

Pour avoir une notation plus expressive, on note plut\^ot $\psi_{>}$ un morphisme dominant $\psi$.  Pour construire $\Pi(\psi)$, on fixe $\psi_{>>}$ un morphisme, de bonne parit\'e,  dominant $\psi$ et tel que la restriction de $\psi_{>>}$ \`a $W_{F}$ fois  la diagonale de $SL(2,{\mathbb C})\times SL(2,{\mathbb C})$ soit sans multiplicit\'e. Ceci est possible car on a suppos\'e $\psi$ de bonne parit\'e. On sait donc construire $\Pi(\psi_{>>})$. On rappelle une notation bien commode; on fixe $\rho$ une repr\'esentation cuspidale autoduale d'un $GL(d)$ et $\tau$ une repr\'esentation d'un groupe classique, $H$, ainsi que $x$ un demi-entier (on n'utilisera que ce cas). Si l'indice de Witt de la forme d\'efinissant $H$ est inf\'erieur \`a $d$ strictement on pose $Jac_{\rho\vert\,\vert^{x}}\tau=0$ et sinon, on fixe un sous-groupe parabolique maximal de $H$ ayant ses sous-groupes de Levi isomorphe \`a $GL(d)\times H'$ pour $H'$ un groupe classique de m\^eme type que $H$ et on
 note $Jac_{\rho\vert\,\vert^{x}}\tau$ l'\'el\'ement du groupe de Grothendieck des repr\'esentations lisses irr\'eductibles de $H'$ tel que dans le groupe de Grothendieck des repr\'esentations lisses irr\'eductibles de $GL(d)\times H'$, le module de Jacquet de $\tau$ le long du radical unipotent du parabolique fix\'e soit de la forme $\rho\vert\,\vert^{x}\otimes Jac_{\rho\vert\,\vert^{x}}\tau$ plus des \'el\'ements de la forme $\rho'\otimes \sigma'$ avec $\rho'\not\simeq \rho\vert\,\vert^{x}$. Quand $\rho$ est fix\'e, on peut l'oublier de la notation. Fixons $\underline{t},\underline{\eta}$ v\'erifiant (1) et (2) ci-dessus; on les transporte en des fonctions sur $Jord(\psi_{>>})$ via l'isomorphisme ordonn\'e de $Jord(\psi_{>>})$ sur $Jord(\psi)$. On a donc d\'efini la repr\'esentation irr\'eductible $\pi(\psi_{>>},\underline{t},\underline{\eta})$. On a d\'emontr\'e en \cite{holomorphie} 2.8 que la repr\'esentation
 $$
 \circ_{(\rho,A,B,\zeta)\in Jord(\psi)}\circ_{\ell\in [\underline{T}(\rho,A,B,\zeta),1]}Jac_{\zeta (A+\ell)}\circ \cdots \circ Jac_{\zeta(B+\ell)} \pi(\psi_{>>},\underline{t},\underline{\eta})=:\pi(\psi,\underline{t},\underline{\eta})
 $$
 o\`u les $(\rho,A,B,\zeta)$ sont pris dans l'ordre d\'ecroissant est un \'el\'ement du groupe de Grothendieck de $G$ qui est soit 0 soit une repr\'esentation irr\'eductible. De plus pour $(\underline{t}',\underline{\eta}')\neq(\underline{t},\underline{\eta})$ soit $\pi(\psi,\underline{t},\underline{\eta})=0=\pi(\psi,\underline{t'},\underline{\eta}')$ soit les 2 repr\'esentations sont in\'equivalentes. Et on a montr\'e en \cite{holomorphie} 2.8 que $\Pi(\psi)$ est exactement l'ensemble des repr\'esentations $\pi(\psi,\underline{t},\underline{\eta}) $ obtenues ainsi.

 Le seul d\'efaut de cette construction est qu'il n'existe pas de crit\`ere simple pour distinguer la nullit\'e ou la non nullit\'e de $\pi(\psi,\underline{t},\underline{\eta})$. Mais par contre on contr\^ole totalement les coefficients n\'ecessaires (des signes) pour que la combinaison lin\'eaire avec coefficients de ces repr\'esentations soit stable et se transf\`ere en la $\theta$-trace de $\pi^{GL}(\psi)$ pour une action totalement explicite de l'automorphisme $\theta$.

Dans le texte on aura besoin de propri\'et\'es suppl\'ementaires que l'on rappelera avec une r\'ef\'erence quand n\'ecessaire et les 2 plus importantes seront red\'emontr\'ees dans la g\'en\'eralit\'e qu'il nous faut en appendice.
\subsection{Notations pour les modules de Jacquet et propri\'et\'e\label{proprietejac}}
Dans le paragraphe pr\'ec\'edent, on a introduit la notation $Jac_{\rho\vert\,\vert^{x}}$ abr\'eg\'ee parfois en $Jac_{x}$. On la g\'en\'eralise ainsi: pour $\rho$ fix\'e et $x_{1}, \cdots, x_{v}$ un ensemble ordonn\'e de nombres r\'eels, on pose
$$
Jac_{x_{1}, \cdots, x_{v}}=\circ_{i\in [v,1]}Jac_{x_{i}}.
$$
On remarque que pour une repr\'esentation irr\'eductible $\pi$, le fait que $Jac_{x_{1}, \cdots, x_{v}}\pi\neq 0$ est \'equivalent \`a ce qu'il existe une repr\'esentation $\sigma$ d'un groupe de m\^eme type que $G$ mais de rang plus petit avec une inclusion:
$$
\pi\hookrightarrow \rho\vert\,\vert^{x_{1}}\times \cdots \times \rho\vert\,\vert^{x_{v}}\times \sigma, \eqno(*)
$$
et on peut \'evidemment supposer $\sigma$ irr\'eductible. On le d\'emontre \`a partir du cas $v=1$ par r\'eciprocit\'e de Frobenius. On v\'erifie aussi que si $x,y$ sont tels que $\vert x-y\vert>1$, on  a l'\'egalit\'e $Jac_{x,y}\pi=Jac_{y,x}\pi$: pour cela on calcule le module de Jacquet pour un parabolique de Levi $GL(2d_{\rho})\times G'$ o\`u $G'$ est un groupe de m\^eme type que $G$ et on d\'ecompose le r\'esultat dans le groupe de Grothendieck en une somme avec coefficients de repr\'esentations $\sigma\otimes \tau$. Ainsi $Jac_{x,y}$ est la somme avec les bons coefficients de $\tau$ qui interviennet avec $\sigma$ qui ont dans leur module de Jacquet le terme $\rho\vert\,\vert^{x}\otimes \rho\vert\,\vert^{y}$. Il n'y a qu'un seul $\sigma$ possible avec cette propri\'et\'e et son module de Jacquet cuspidal est exactement la somme de ce terme avec le terme $\rho\vert\,\vert^{y}\otimes \rho\vert\,\vert^{x}$. D'o\`u le r\'esultat.

\subsection{Le facteur de normalisation et ses p\^oles\label{poles}}
On fixe $\psi$ comme ci-dessus, $\rho$ une repr\'esentation cuspidale unitaire autoduale  d'un $GL$ et $a_{0}$ un entier strictement positif  et on a pos\'e:
$
r(s,\rho,a_{0},\psi):=$
 $$\times_{(\rho,a,b)\in Jord(\psi)}\frac{L(St(\rho,a_{0})\times St(\rho,a), s-(b-1)/2)}{L(St(\rho,a_{0})\times St(\rho,a),s+(b+1)/2)}\times \frac{L(St(\rho,a_{0}),r_{G,ad_{\rho}},2s)}{L(St(\rho,a_{0}),r_{G,ad_{\rho}},2s+1)}.
 $$
 Les facteurs $L$ locaux qui interviennent ci-dessus sont bien connus; par exemple, on sait depuis \cite{jpss} page 445, (6) que pour $\sigma$, $\tau$ des repr\'esentations de carr\'e int\'egrable d'un $GL$, le facteur $L(\sigma\times \tau,s)$ n'a pas de p\^ole pour $s\in {\mathbb R}_{>0}$. On en d\'eduit, a fortiori, que les d\'enominateurs ci-dessus n'ont pas de p\^oles pour $s=s_{0}>0$ et qu'il en est de m\^eme de $L(St(\rho,a_{0}),r_{G,ad_{\rho}},2s)$ . Ainsi l'ordre en $s=s_{0}$ de la fonction m\'eromorphe $r(s,\rho,a_{0},\psi)$ est la somme des ordres en $s=s_{0}$ des fonctions $L(St(\rho,a_{0})\times St(\rho,a), s-(b-1)/2)$ quand $(\rho,a,b)$ parcourt $Jord(\psi)$. Pour d\'eterminer ces ordres, on utilise le calcul explicite de ces facteurs donn\'e en \cite{jpss} th. 8.2:
 $
 L(St(\rho,a_{0})\times St(\rho,a), s')=\times_{\ell\in [(a-1)/2,-(a-1)/2]}L(\rho\times \rho, (a_{0}-1)/2+\ell+s')$ si $a\geq a_{0}$ et une formule sym\'etrique dans le cas inverse. Cela se r\'ecrit
 $$
 L(St(\rho,a_{0})\times St(\rho,a), s-(b-1)/2)=\times_{\ell \in [\vert (a-a_{0})/2\vert, (a+a_{0})/2-1]}L(\rho\times \rho, \ell+s-(b-1)/2).
 $$Une telle fonction a donc au plus un p\^ole simple et elle en a un exactement quand $(b-1)/2- s_{0}\in[\vert (a-a_{0})/2\vert, (a+a_{0})/2-1]$. 
 On pose $b_{0}:=2s_{0}+1$ et on suppose comme pr\'ec\'edemment que $b_{0}$ est un entier sup\'erieur ou \'egal \`a $2$ et on reprend les notations: $A=(a+b)/2-1, B=\vert (a-b)\vert/2$, $\zeta$ est un signe tel que $\zeta(a-b)\geq 0$ et de m\^eme $A_{0}=(a_{0}+b_{0})/2-1,B_{0}=\vert (a_{0}-b_{0})$ et un signe tel que $\zeta_{0}(a_{0}-b_{0})\geq 0$ avec $\zeta_{0}=+$ si $a_{0}=b_{0}$. On a donc un p\^ole si les in\'egalit\'es suivantes sont satisfaites:
 $$
 (b-b_{0})/2\geq (a-a_{0})/2; (b-b_{0})/2 \geq (a_{0}-a)/2; (b-b_{0})/2 \leq (a+a_{0})/2-1.
 $$
 Ou encore avec les notations introduites:
 $$
 \zeta_{0}B_{0}-\zeta B\geq 0; A\geq A_{0}; \zeta B+A_{0}\geq 0.
 $$
On a alors le tableau suivant qui dit quand $(\rho,a,b)$ participe aux p\^oles de la fonction $r(s,\rho,a_0,\psi)$ en $s=s_{0}$ (on rappelle que l'on a suppos\'e que $\zeta_{0}=+$ si $B_{0}=0$):
\begin {center}
\begin{tabular}{| c | c | c |}
\hline
$\zeta\backslash \zeta_{0}$ & + & - 
\\ \hline
 + & $B\leq B_{0} \leq A_{0}\leq A$  &  \hbox {pas de p\^ole}
\\\hline
- & $B\leq A_{0}\leq A$  & $B_{0}\leq B\leq A_{0}\leq A$
\\\hline
\end{tabular}
\end {center}

\subsection{Param\`etres de l'image des op\'erateurs d'entrelacement dans le cas de bonne parit\'e}
On fixe $\rho,a_{0},b_{0}$ un triplet qui param\'etrise donc une repr\'esentation irr\'eductible de $W_{F}\times SL(2,{\mathbb C})\times SL(2,{\mathbb C})$; on suppose que cette repr\'esentation est autoduale \`a valeurs dans un groupe de m\^eme type (orthogonal ou symplectique) que le type du groupe dual de $G$; c'est-\`a-dire qu'elle est de bonne parit\'e. On pose $A_{0}=(a_{0}+b_{0})/2-1$ et $B_{0}=\vert a_{0}-b_{0}\vert/2$; on note $\zeta_{0}$ le signe de $a_{0}-b_{0}$ sauf si ce nombre est 0 o\`u $\zeta_{0}$ vaut alors $+$ par d\'efinition. On consid\`ere aussi le triplet $(\rho,a_{0},b_{0}-2)$ qui n'est d\'efini que si $b_{0}>2$ et on lui associe aussi $A'_{0}:=(a_{0}+b_{0}-2)/2-1=A_{0}-1$, $B'_{0}=\vert (a_{0}-b_{0}+2)/2\vert$ qui vaut $B_{0}+\zeta_{0}$ sauf si $a_{0}=b_{0}-1$ o\`u l'on a $B_{0}=B'_{0}=1/2$. On note $\zeta'_{0}$ le signe de $a_{0}-b_{0}+2$ et si $a_{0}=b_{0}-2$, on pose $\zeta'_{0}=-$; on a donc $\zeta'_{0}=\zeta_{0}$ sauf dans le cas particulier o\`u $B_{0}=B'_{0}=1/2$ avec $\zeta_{0}=-$ o\`u $\zeta'_{0}=+$. Si $b_{0}=2$, on dit que $(\rho,A'_{0},B'_{0},\zeta'_{0})$ n'existe pas.

On fixe un morphisme $\psi$ comme dans les paragraphes pr\'ec\'edents.  On suppose que $\psi$ est de bonne parit\'e.
On suppose de plus que si $b_{0}>2$ alors $\psi$ contient la sous-repr\'esentation irr\'eductible associ\'ee au triplet $(\rho,a_{0},b_{0}-2)$. Dans tous les cas, on note $\psi^+$ le morphisme qui se d\'eduit de $\psi$ en rempla\c{c}ant une copie de $(\rho,a_{0},b_{0}-2)$ (si il y en a une) par une copie de $(\rho,a_{0},b_{0})$ ou en ajoutant simplement $(\rho,a_{0},2)$ si $b_{0}=2$.

En termes de quadruplet, l'hypoth\`ese sur $\psi$ est que  $Jord(\psi)$ contient $(\rho,A'_{0},B'_{0},\zeta'_{0})$ si ce quadruplet est d\'efini.

On veut d\'efinir a priori un \'el\'ement  $\pi^+$ dans le paquet associ\'e \`a $\psi^+$ et montrer que l'image de $N_{\psi}(s,\rho,a_{0})_{s=s_{0}}$ est exactement $\pi^+$ c'est \`a dire est nulle si $\pi^+$ est nulle et est r\'eduite \`a $\pi^+$ sinon. Le guide pour faire cela est que l'on veut une inclusion
$$
\pi^+\hookrightarrow St(\rho,a_{0})\vert\,\vert^{-s_{0}}\times \pi. \eqno(1)
$$
Dans certains cas, cela d\'etermine uniquement les param\`etres de $\pi^+$. C'est ce qui se produit dans le cas o\`u la restriction de $\psi^+$ et la restriction de $\psi$  \`a $W_{F}$ fois la diagonale de $SL(2,{\mathbb C})\times SL(2,{\mathbb C})$ sont sans multiplicit\'e. On retrouvera cela dans la preuve de \ref{identification}. Un cas est compl\`etement \'evident et on va donc le donner d\`es maintenant; il s'agit du cas o\`u $\zeta_{0}=+$ (c'est-\`a-dire $a_{0}\geq b_{0}$) et o\`u la restriction de $\psi^+$ \`a $W_{F}$ fois la diagonale de $SL(2,{\mathbb C})\times SL(2,{\mathbb C})$ est sans multiplicit\'e. On v\'erifie alors que $(\rho,A'_{0},B'_{0};\zeta'_{0})=(\rho,A_{0}-1,B_{0}+1,\zeta_{0})$ et la restriction de $\psi$ \`a $W_{F}$ fois la diagonale de $SL(2,{\mathbb C})\times SL(2,{\mathbb C})$ est aussi sans multiplicit\'e. Notons $\underline{t},\underline{\eta}$ les param\`etres de $\pi$ dans son paquet $\Pi(\psi)$ et $\underline{t}^+,\underline{\eta}^+$ ceux de $\pi^+$ dans $\Pi(\psi^+)$. D'apr\`es (1) et \cite{discret} 4.2, $\underline{t},\underline{\eta}$ et $\underline{t}^+,\underline{\eta}^+$ co\"{\i}ncident sur $Jord(\psi)\cap Jord(\psi^+)$ et l'on a
$$
\underline{t}^+(\rho,A_{0},B_{0},\zeta_{0})=\underline{t}(\rho,A'_{0},B'_{0},\zeta'_{0})+1\hbox{ (resp. $=1$)};$$
$$
\underline{\eta}^+(\rho,A_{0},B_{0},\zeta_{0})=\underline{\eta}(\rho,A'_{0},B'_{0},\zeta'_{0})\hbox{ (resp. $=+$ )}
$$
si $(\rho,A'_{0},B'_{0},\zeta'_{0})$ existe (resp. n'existe pas).

Le cas o\`u $\zeta_{0}=-$ m\^eme avec l'hypoth\`ese que la restriction de $\psi^+$ \`a $W_{F}$ fois la diagonale de $SL(2,{\mathbb C})\times SL(2,{\mathbb C})$ est sans multiplicit\'e est un peu plus compliqu\'e; en effet, l'hypoth\`ese n'entra\^{\i}ne pas que la restriction de $\psi$ \`a $W_{F}$ fois la diagonale de $SL(2,{\mathbb C})\times SL(2,{\mathbb C})$ soit aussi sans multiplicit\'e. Mais si on ajoute aussi cette hypoth\`ese, l'inclusion (1) et \cite{discret} 4.4  forcent les formules
$$
\underline{t}^+(\rho,A_{0},B_{0},\zeta_{0})=\underline{t}(\rho,A'_{0},B'_{0},\zeta'_{0})\hbox{ (resp. $=0$)};$$
$$
\underline{\eta}^+(\rho,A_{0},B_{0},\zeta_{0})=\underline{\eta}(\rho,A'_{0},B'_{0},\zeta'_{0})\hbox{ (resp. $=+$ )}
$$
si $(\rho,A'_{0},B'_{0},\zeta'_{0})$ existe (resp. n'existe pas).

On va maintenant donner a priori la formule pour les param\`etres de $\pi^+$ en g\'en\'eralisant les cas ci-dessus. Auparavant il faut r\'egler le probl\`eme du choix d'un ordre sur $Jord(\psi^+)$ qui n'apparaissait pas ci-dessus \`a cause des hypoth\`eses simplificatrices. Comme on veut des formules analogues \`a celles donn\'ees ci-dessus, l'id\'ee la plus simple, qui est celle que nous prendrons, est de mettre sur $Jord(\psi^+)$ l'ordre tel que $(\rho,A_{0},B_{0},\zeta_{0})$ prenne la place de $(\rho,A'_{0},B'_{0},\zeta'_{0})$ dans l'ordre sur $Jord(\psi)$. Mais  on doit garder la propri\'et\'e donn\'ee en \ref{descriptionpaquets} ce qui n\'ecessite d'imposer \`a l'ordre sur $Jord(\psi)$ certaines propri\'et\'es; si $(\rho,A'_{0},B'_{0},\zeta'_{0})$ n'existe pas, on peut faire le choix que l'on veut \`a condition de respecter cette propri\'et\'e. On a aussi besoin de mettre des conditions sur les ordres pour pouvoir faire les d\'emonstrations et je ne sais pas si ces conditions sont n\'ecessaires, en tout cas, elles compliquent la situation.

A la fin du paragraphe suivant on donne une famille d'ordre qui v\'erifie toutes les conditions voulues mais auparavant on met les conditions minimales pour savoir d\'ecrire l'image de $N_{\psi}(s,\rho,a_{0},\pi)$ en $s=(b_{0}-1)/2$ en termes de param\`etres. Cela est  technique mais laisse plus de flexibilit\'e que de n'accepter que les ordres ayant les propri\'et\'es de la fin du paragraphe.

\subsubsection{Propri\'et\'es de l'ordre sur $Jord(\psi)$ et sur $Jord(\psi^+)$ \label{proprietedelordre}}
On garde les notations pr\'ec\'edentes. On consid\`ere des ordres sur $Jord(\psi)$ et $Jord(\psi^+)$ v\'erifiant ${\mathcal P}$ de \ref{descriptionpaquets} et on demande en plus les conditions ci-dessous:

l'ordre sur $Jord(\psi)$ et celui sur $Jord(\psi^+)$ sont tels que les quadruplets de $Jord(\psi)$ donnant des p\^oles \`a la fonction $r(s,\rho,a_0,\psi)$ en $s=(b_{0}-1)/2$ (cf. \ref{poles}), sont plus grands que $(\rho,A'_{0},B'_{0},\zeta'_{0})$ (si cet \'el\'ement existe) dans $Jord(\psi)$ et que $(\rho,A_{0},B_{0},\zeta_{0})$ dans $Jord(\psi^+)$ et plus grand que les quadruplets $(\rho,A',B',\zeta')$ avec $A'\leq  A'_{0}$ dans $Jord(\psi)$ et dans $Jord(\psi^+)$. Avec des symboles  cela se traduit par:

\

$({\mathcal P})_{p}$: 
pour tout $(\rho,A,B,\zeta)$ contribuant aux p\^oles de $r(s,\rho,a_0,\psi)$ en $s=(b_{0}-1)/2$, on a dans $Jord(\psi)$ (resp. $Jord(\psi^+)$), $(\rho,A,B,\zeta)>(\rho,A'_{0},B'_{0},\zeta'_{0})$ (resp. $(\rho,A_{0},B_{0},\zeta_{0})$).

pour tout $(\rho,A,B,\zeta),(\rho,A',B',\zeta')$ dans $Jord(\psi)\setminus\{(\rho,A'_{0},B'_{0},\zeta'_{0})\}$, tels que $A<A_{0}$ et $(\rho,A',B',\zeta')$ contribue aux p\^oles de $r(s,\rho,a_0,\psi)$ en $s=(b_{0}-1)/2$ (en particulier $A'\geq A_{0}$), on a $(\rho,A,B,\zeta)<(\rho,A',B',\zeta')$.

\

La deuxi\`eme partie de ${\mathcal P}_{p}$  a la  cons\'equence suivante:  notons $(\rho,A_{1},B_{1},\zeta_{1})$ le plus petit \'el\'ement de l'ensemble des quadruplets fournissant des p\^oles \`a $r(s,\rho,a_0,\psi)$ en $s=(b_{0}-1)/2$ (quand cet ensemble est non vide), alors  tout $(\rho,A,B,\zeta)\geq (\rho,A_{1},B_{1},\zeta_{1})$ v\'erifie $A\geq A_{0}$.

\

On impose aussi la condition suivante dans le cas exceptionnel:

supposons que $(\rho,a_{0},b_{0})$ est tel que $b_{0}=a_{0}+1$, c'est-\`a-dire $B_{0}=1/2$ et $\zeta'_{0}=+=-\zeta_{0}$; ici on impose aux ordres d'\^etre tels que $(\rho,A'_{0},B'_{0},\zeta'_{0})$ et $(\rho,A_{0},B_{0},\zeta_{0})$ sont les \'el\'ements minimaux. Si $b_{0}=2$ et $a_{0}=1$, on impose que $(\rho,A_{0},B_{0},\zeta_{0})$ soit le plus petit \'el\'ement de $Jord(\psi^+)$.

\

On a  une condition propre au cas $\zeta_{0}=+$ qui est un peu plus difficile que le cas $\zeta_{0}=-$:

(0) supposons que  $\zeta_{0}=+$; soit $(\rho,A,B,\zeta)\in Jord(\psi)$ diff\'erent de $(\rho,A'_{0},B'_{0},\zeta'_{0})$ si cet \'el\'ement existe. On suppose que $A<A_{0}$ et $\zeta=\zeta_{0}$, alors si $(\rho,A,B,\zeta)>(\rho,A_{0},B_{0},\zeta_{0})$ dans $Jord(\psi^+)$ ou \`a $(\rho,A'_{0},B'_{0},\zeta'_{0})$ dans $Jord(\psi)$, on a $B>B_{0}+1$. (cette propri\'et\'e ne servira que dans \ref{identification})

\

Comme on veut que $(\rho,A,B,\zeta)>(\rho,A_{0},B_{0},\zeta_{0})$ dans $Jord(\psi^+)$ si et seulement si $(\rho,A,B,\zeta)>(\rho,A'_{0},B'_{0},\zeta'_{0})$ dans $Jord(\psi)$, il faut avoir des conditions telles que ${\mathcal P}$ soit respect\'e. Ceci impose des conditions aux limites; on suppose que l'on est pas dans le cas exceptionnel d\'ej\`a r\'egl\'e par la condition forte ci-dessus et on explique ce que l'on entend par conditions aux limites: soit $(\rho,A,B,\zeta)$ tel que $\zeta=\zeta_{0}=\zeta'_{0}$. Par exemple  si $A=A_{0}$, $(\rho,A,B,\zeta)$ peut \^etre plus petit ou plus grand que $(\rho,A_{0},B_{0},\zeta_{0})$ et le signe de $B-B_{0}$ n'intervient pas alors que l'on a $A>A_{0}-1=A'_{0}$ et  si $B>B_{0}+\zeta_{0}=B'_{0}$ il faut n\'ecessairement $(\rho,A,B,\zeta) > (\rho,A'_{0},B'_{0},\zeta_{0})$ dans $Jord(\psi)$ donc on veut aussi l'in\'egalit\'e $(\rho,A,B,\zeta)>(\rho,A_{0},B_{0},\zeta_{0})$ dans $Jord(\psi^+)$. Il y a d'autres conditions du m\^eme ordre, on les \'ecrit toutes ci-dessus, o\`u on a fix\'e $(\rho,A,B,\zeta)$ avec $\zeta=\zeta_{0}=\zeta'_{0}$:

(1) si $A=A_{0}=A'_{0}-1$ et $B>B'_{0}=B_{0}+\zeta_{0}$ alors $(\rho,A,B,\zeta)>(\rho,A_{0},B_{0},\zeta_{0})$ dans $Jord(\psi^+)$

(2) si $A=A'_{0}=A_{0}-1$ et si $B<B_{0}$ alors $(\rho,A,B,\zeta)<(\rho,A'_{0},B'_{0},\zeta_{0})$ dans $Jord(\psi)$

(3) si $B=B_{0}=B'_{0}-\zeta_{0}$: si $\zeta_{0}=+$ et si $A< A'_{0}=A_{0}-1$, alors $(\rho,A,B,\zeta)<(\rho,A_{0},B_{0},\zeta_{0})$ dans $Jord(\psi^+)$; si $\zeta_{0}=-$ et si $A\geq A_{0}=A'_{0}+1$ alors $(\rho,A,B,\zeta)>(\rho,A_{0},B_{0},\zeta_{0})$ dans $Jord(\psi^+)$

(4) si $B=B'_{0}=B_{0}+\zeta_{0}$: si $\zeta_{0}=+$ et si $A>A_{0}$ alors $(\rho,A,B,\zeta)>(\rho,A'_{0},B'_{0},\zeta_{0})$ dans $Jord(\psi)$; si $\zeta_{0}=-$ et si $A<A_{0}$ alors $(\rho,A,B,\zeta)<(\rho,A'_{0},B'_{0},\zeta_{0})$ dans $Jord(\psi)$.

\

Fixons un ordre sur $Jord(\psi)$ en supposant que $b_{0}> 2$; on d\'eduit un ordre sur $Jord(\psi^+)$ en demandant simplement que $(\rho,A_{0},B_{0},\zeta_{0})$ prenne la place de $(\rho,A'_{0},B'_{0},\zeta'_{0})$, c'est-\`a-dire pour tout $(\rho,A,B,\zeta)\in Jord(\psi^+)\setminus \{(\rho,A_{0},B_{0},\zeta_{0})$ on a $(\rho,A,B,\zeta)>(\rho,A_{0},B_{0},\zeta_{0})$ dans $Jord(\psi^+)$ si et seulement si $(\rho,A,B,\zeta)>(\rho,A'_{0},B'_{0},\zeta'_{0})$ dans $Jord(\psi)$. 

R\'eciproquement, sans hypoth\`ese sur $b_{0}$,  ayant un ordre sur $Jord(\psi^+)$ on en d\'eduit un ordre sur $Jord(\psi)$ en demandant que $(\rho,A'_{0},B'_{0},\zeta'_{0})$ prenne la place de $(\rho,A_{0},B_{0},\zeta_{0})$.

Montrons que les op\'erations que l'on a d\'efinies entre les ordres sur $Jord(\psi)$ et $Jord(\psi^+)$ conservent \`a la fois ${\mathcal P}$, ${\mathcal P}_{p}$ et les propri\'et\'es ci-dessus. Pour ${\mathcal P}_{p}$ c'est \'evident ainsi que pour (0) car elles sont sym\'etriques en $Jord(\psi)$ et $Jord(\psi^+)$. Les autres conditions sont exactement faites pour que la propri\'et\'e ${\mathcal P}$ soit respect\'ee dans les cas limites.

\

Remarquons qu'il existe des ordres v\'erifiant les conditions ci-dessus: on prend un ordre sur $Jord(\psi^+)\setminus \{(\rho,A_{0},B_{0},\zeta_{0})\}$ tel que $(\rho,A,B,\zeta)>(\rho,A',B',\zeta')$ d\`es que $A>A'$. On compl\`ete cet ordre en un ordre sur $Jord(\psi)$ en demandant que $(\rho,A,B,\zeta)>(\rho,A'_{0},B'_{0},\zeta'_{0})$ si et seulement si $A>A'_{0}=A_{0}-1$. On le compl\`ete aussi en un ordre sur $Jord(\psi^+)$ en demandant que $(\rho,A,B,\zeta)>(\rho,A_{0},B_{0},\zeta_{0})$ si et seulement si $A\geq A_{0}$. Alors ces ordres se d\'eduisent l'un de l'autre par la proc\'edure expliqu\'ee: $(\rho,A_{0},B_{0},\zeta_{0})$ prend la place de $(\rho,A'_{0},B'_{0},\zeta'_{0})$ ou lui c\`ede sa place. Et ${\mathcal P}$, ${\mathcal P}_{p}$, (1) \`a (3) sont v\'erifi\'ees. En imposant les conditions minimales, on a voulu garder un maximum de flexibilit\'e.

\subsubsection{Description des param\`etres\label{descriptionparametre}}
On fixe un ordre sur $Jord(\psi)$ satisfaisant toutes les propri\'et\'es de \ref{proprietedelordre} et on en d\'eduit un ordre sur $Jord(\psi^+)$ comme expliqu\'e en loc.cite; si $b_{0}=2$, l'ordre sur $Jord(\psi^+)$ n'est pas uniquement d\'etermin\'e par celui sur $Jord(\psi)$.

On a des param\`etres $\underline{t},\underline{\eta}$ pour d\'efinir $\pi$. On d\'efinit $\underline{t}^+, \underline{\eta}^+$ pour d\'efinir un \'el\'ement (\'eventuel\-lement 0) de $\Pi(\psi^+)$: sur tout $(\rho,A,B,\zeta)\in Jord(\psi)-\{ (\rho,A'_{0},B'_{0},\zeta'_{0})\}$, $\underline{t},\underline{\eta}$ co\"{\i}ncide avec $\underline{t}^+, \underline{\eta}^+$ et il faut donc d\'efinir ces applications sur $(\rho,A_{0},B_{0},\zeta_{0})$. On pose $t_{0}:=\underline{t}(\rho,A'_{0},B'_{0},\zeta'_{0})$ et des notations analogues pour $t^+_{0}$, $\eta_{0}$ et $\eta^+_{0}$. On donne d'abord les d\'efinitions dans le cas o\`u $b_{0}>2$ et $\zeta_{0}=\zeta'_{0}$ ainsi que dans le cas o\`u $b_{0}=2$, o\`u on pose $t_{0}=0$ et $\eta_{0}=+$ (ces nombres ne sont pas d\'efinis car $(\rho,A'_{0},B'_{0},\zeta'_{0})$ n'existe pas). Dans tous les cas, on a
$$
\eta_{0}=\eta^+_{0}$$
\begin{center}
\begin{tabular}{|c | c| c|}
\hline$
\zeta_{0}$&+&-
\\
\hline$
t^+_{0}$&$ t_{0}+1 $&$ t_{0}$
\\
\hline
\end{tabular}
\end{center}
Il reste le cas exceptionnel $\zeta_{0}=-$, $\zeta'_{0}=+$, $B'_{0}=B_{0}=1/2$ qui est d\'esagr\'eable dans notre param\'etrisation. On a 
$$
\eta_{0}^+=-\eta_{0}$$
et la valeur de $t^+_{0}$ d\'epend des valeurs de $\eta_{0}$ et $t_{0}$ par la formule:
\begin{center}
\begin{tabular}{|c | c| c|}
\hline$
\eta_{0}$&+&-
\\
\hline$
t^+_{0}$&$ t_{0}+1 $&$ t_{0}$
\\
\hline
\end{tabular}
\end{center}
Il y a 2 v\'erifications \`a faire pour s'assurer que ces d\'efinitions ont les propri\'et\'es de \ref{descriptionpaquets}: on doit v\'erifier que $t_{0}^+\leq [inf(a_{0},b_{0})/2]$ et que 
$$
\eta_{0}^{inf(a_{0},b_{0}-2)}(-1)^{[inf(a_{0},b_{0}-2)/2]+t_{0}}=(\eta_{0}^{+})^{inf(a_{0},b_{0})}(-1)^{[inf(a_{0},b_{0})/2]+t_{0}^+}.
$$
On sait que $t_{0}\leq [inf(a_{0},b_{0}-2)/2]$ et on a $t_{0}=t_{0}^+$ chaque fois que $a_{0}< b_{0}$ et l'on a donc l'in\'egalit\'e cherch\'e sur $t_{0}^+$. Montrons l'\'egalit\'e des signes, en notant $\epsilon_{0}$ le membre de gauche et $\epsilon_{0}^+$ celui de droite. On commence par le cas o\`u soit $b_{0}=2$ soit $\zeta_{0}=\zeta'_{0}$ et on a
\begin{center}
\begin{tabular}{|c | c| c|}
\hline$
\zeta_{0}$&+&-
\\
\hline$
\epsilon_{0}$&$\eta_{0}^{b_{0}-2}(-1)^{[b_{0}/2]-1+t_{0}} $&$\eta_{0}^{a_{0}}(-1)^{[a_{0}/2]+t_{0}}$
\\
\hline
$\epsilon^+_{0} $& $\eta_{0}^{b_{0}}(-1)^{[b_{0}/2]+t_{0}+1}$&$\eta_{0}^{a_{0}}(-1)^{[a_{0}/2]+t_{0}}$
\\ \hline
\end{tabular}
\end{center}
et l'\'egalit\'e cherch\'ee. On fait maintenant le calcul quand $\zeta_{0}=+$ et $\zeta'_{0}=-$, c'est-\`a-dire $a_{0}=b_{0}-1$:
\begin{center}
\begin{tabular}{|c | c| c|}
\hline$
\eta_{0}$&+&-
\\
\hline$
\epsilon_{0}$&$(-1)^{[b_{0}/2]-1+t_{0}} $&$(-1)^{b_{0}-2}(-1)^{[b_{0}/2]-1+t_{0}}$
\\
\hline
$\epsilon^+_{0} $& $(-1)^{b_{0}-1}(-1)^{[(b_{0}-1)/2]+t_{0}+1}$&$(-1)^{[(b_{0}-1)/2]+t_{0}}$
\\ \hline
\end{tabular}
\end{center}
Or on a l'\'egalit\'e pour tout nombre entier $b'$
$$
(-1)^{[b'/2]}=(-1)^{b'-1}(-1)^{[(b'-1)/2]}.
$$
En reportant dans le tableau, avec $b'=b_{0}$ pour la 2e colonne et $b'=b_{0}+1$ pour la 3e colonne, on trouve l'\'egalit\'e cherch\'ee de $\epsilon_{0}$ et $\epsilon^+_{0}$.

\

Avec ces param\`etres, on d\'efinit une repr\'esenation $\pi ^+$ dans $\Pi(\psi^+)$ et
il n'est \'evidemment pas clair que la repr\'esentation $\pi^+$ (et sa nullit\'e \'eventuelle) ne d\'epende pas du choix de l'ordre mis sur $Jord(\psi)$ mais ce sera une cons\'equence de son identification avec l'image de $N_{\psi}(s,\rho,a_{0},\pi)_{s=(b_{0}-1)/2}$.
\section{Description de $\pi^+$ comme sous-module irr\'eductible dans le cas holomorphe }
Le cas o\`u $ordre_{s=s_{0}} r(s,\rho,a_{0},\psi)=0$ est dit le cas holomorphe; on a d\'ecrit les p\^oles et donc leur absence dans \ref{poles}. Cela se traduit par:

 si $\zeta_{0}=+$, pour tout $(\rho,A,B,\zeta)$ soit $A>A_{0}$ soit $\zeta=+$ et $B>B_{0}$ soit $\zeta=-$ et $B>A_{0}$;

 si $\zeta_{0}=-$, pour tout $(\rho,A,B,\zeta)$ soit $\zeta=+$ soit $\zeta=-$ et soit $B>A_{0}$ soit $A<A_{0}$.

\subsection{Etude de l'induite $St(\rho,a_{0})\vert\,\vert^{-s_{0}}\times \pi$ dans le cas holomorphe\label{remarque}}
 Supposons que le facteur de normalisation de \ref{poles} n'a pas de p\^ole en $s=s_{0}=(b_{0}-1)/2$.
 \begin{rem} Avec l'hypoth\`ese faite, la repr\'esentation induite $
St(\rho,a_{0})\vert\,\vert^{-(b_{0}-1)/2}\times \pi
 $ a un unique sous-module irr\'eductible. Ce sous-module irr\'eductible est l'unique sous-quotient, $\pi^0$ de l'induite v\'erifiant $Jac_{(a_{0}-b_{0})/2, \cdots, -(a_{0}+b_{0})/2+1}\pi^0\neq 0$. Et il intervient avec multiplicit\'e 1 en tant que sous-quotient irr\'eductible de cette induite.
 \end{rem}
Pour d\'emontrer la remarque, il suffit de prouver l'\'egalit\'e (dans le groupe de Grothendieck convenable)
$$
Jac_{(a_{0}-b_{0})/2, \cdots, -(a_{0}+b_{0})/2+1}\biggl(St(\rho,a_{0})\vert\,\vert^{-(b_{0}-1)/2}\times \pi
\biggr)=\pi.\eqno(1)
$$Car l'exactitude du foncteur de Jacquet assure alors qu'un unique sous-quotient irr\'eductible, $\sigma$ de l'induite v\'erifie $Jac_{(a_{0}-b_{0})/2, \cdots, -(a_{0}+b_{0})/2+1}\sigma\neq 0$ et la r\'eciprocit\'e de Frobenius  assure que tout sous-module irr\'eductible de l'induite \`a cette propri\'et\'e.

Montrons (1). On reprend les notations commodes, $(a_{0}-b_{0})/2=\zeta_{0} B_{0}$ et $A_{0}=(a_{0}+b_{0})/2-1$.
Le terme de gauche se calcule avec les formules de Bernstein-Zelevinsky: il a une filtration dont les quotients sont index\'es par les d\'ecompositions de l'intervalle $[\zeta_{0}B_{0},-A_{0}]$ en 3 sous-ensembles, ${\mathcal E}_{i}$ pour $i=1,2,3$,   totalement ordonn\'es par l'ordre induit  tels que
$$
Jac_{x\in {\mathcal E}_{1}}<\zeta_{0}B_{0}, \cdots, -A_{0}>_{\rho}\neq 0; \circ_{x\in {\mathcal E}_{2}}Jac_{x}<A_{0}, \cdots, -\zeta_{0}B_{0}>_{\rho}\neq 0; Jac_{x\in {\mathcal E}_{3}}\pi\neq 0.
$$
Ainsi ${\mathcal E}_{i}$ pour $i=1,2$ est un intervalle \'eventuellement vide; si ${\mathcal E}_{1}$ est non vide c'est n\'ecessairement un intervalle commen\c{c}ant par $\zeta_{0}B_{0}$ et si ${\mathcal E}_{2}$ est non vide c'est n\'ecessairement un intervalle commen\c{c}ant par $A_{0}$; comme $A_{0}>B_{0}$ n'est pas dans l'intervalle $[\zeta_{0}B_{0},-A_{0}]$, on sait que ${\mathcal E}_{2}$ est vide. Ainsi 
${\mathcal E}_{3}$ est soit vide soit est un sous-intervalle de la forme $[x,-A_{0}]$ de $[\zeta_{0}B_{0},-A_{0}]$. Ainsi en appliquant \ref{proprietedujac}, on voit que si ${\mathcal E}_{3}$ n'est pas vide, l'une des conditions du tableau ci-dessous est satisfaite:
\begin{center}
\begin{tabular}{|c|c|c|}
\hline$\zeta'\backslash \zeta_{0} $&$+$ &$-$ \\
\hline
$+$ & $B'\leq B_{0}<A_{0}\leq A'$ & {Impossible} \\
\hline
$-$ &$B'\leq A_{0}\leq A' $& $B_{0}\leq B'\leq A_{0}\leq A' $
\\
\hline\end{tabular}
\end{center}
Ce tableau est exactement le m\^eme que celui de \ref{poles} et aucune de ces conditions ne peut donc \^etre satisfaite. On obtient donc toute la remarque.
\subsection{Identification de $\pi^+$ dans le cas holomorphe}
On continue avec les hypoth\`eses de \ref{remarque}, c'est-\`a-dire que $r(s,\rho,a_0,\psi)$ n'a pas de p\^ole en $s=s_{0}=(b_{0}-1)/2$.
\subsubsection{Le cas tr\`es dominant\label{tresdominant}}
On suppose ici en plus que pour tout $(\rho,A',B',\zeta')\in Jord(\psi)\setminus (\rho,A'_{0},B'_{0},\zeta'_{0})$ sup\'erieur  \`a $(\rho,A_{0},B_{0},\zeta_{0})$ dans $Jord(\psi^+)$ on a $B'>>A_{0}$ et si $(\rho,A'',B'',\zeta'')$ v\'erifie les m\^emes propri\'et\'es, on a soit $B''>>A'$ soit $B'>>A''$.
\begin{lem} Avec les hypoth\`eses pr\'ec\'edentes, la repr\'esentation $\pi^0$ de \ref{remarque} est isomorphe \`a $\pi^+$.
\end{lem}
On commence par les cas o\`u soit $\zeta_{0}\neq -$ soit $B_{0}\neq 1/2$. On fixe un morphisme $\psi_{>}^+$ qui s'obtient \`a partir de $\psi^+$ en ne changeant que $(\rho,A_{0},B_{0},\zeta_{0})$ en $(\rho,A_{0}+T,B_{0}+T,\zeta_{0})$ pour $T$ suffisamment grand mais pas trop pour que pour tout $(\rho,A',B',\zeta')<(\rho,A_{0},B_{0},\zeta_{0})$ on ait $B_{0}+T>A'$ et pour tout $(\rho,A',B',\zeta')>(\rho,A_{0},B_{0},\zeta_{0})$ on ait encore $B'>A_{0}+T$. Avec les param\`etres de $\pi^+$, on construit $\pi_{>}^+$ et par d\'efinition, on a:
$$
\pi^+=\circ_{\ell\in [1,T]}Jac_{\zeta_{0}(B_{0}+\ell), \cdots, \zeta_{0}(A_{0}+\ell)} \pi^+_{>}.
$$
On r\`egle le cas o\`u $\zeta_{0}=-$; comme par hypoth\`ese $B_{0}\neq 1/2$, on n'est pas dans le cas $b_{0}=2$ car il faut alors $a_{0}=1$ par $\zeta_{0}=-$ ce qui contredit $B_{0}\neq 1/2$. On a donc aussi $\psi_{>}$ qui s'obtient en rempla\c{c}ant $(\rho,A'_{0},B'_{0},\zeta'_{0})$ par $(\rho,A'_{0}+T+1,B'_{0}+T+1,\zeta'_{0})$, ce qui n'est autre que $\psi^+_{>}$. Par construction $\pi^+_{>}$ a les param\`etres de $\pi$, d'o\`u aussi:
$$
\pi=\circ_{\ell\in [1,T+1]}Jac_{\zeta_{0}(B'_{0}+\ell), \cdots, \zeta_{0}(A'_{0}+\ell)}\pi^+_{>}
$$Mais ici $B'_{0}=B_{0}-1$ et $A'_{0}=A_{0}-1$ d'o\`u
$\pi=Jac_{\zeta_{0}B_{0},\cdots, \zeta_{0}A_{0}}\pi^+$. Cela donne la non nullit\'e de $\pi^+$ et une inclusion $$\pi^+\hookrightarrow \rho\vert\,\vert^{\zeta_{0}B_{0}}\times \cdots \times \rho\vert\,\vert^{\zeta_{0}A_{0}}\times \pi.$$
Comme $Jac_{x, \cdots, \zeta_{0}(A_{0})}\pi^+=0$ pour tout $x\in ]\zeta_{0}B_{0},\zeta_{0}A_{0}]$ (cf \ref{proprietedujac}), cette inclusion se factorise par l'unique sous-module irr\'eductible de l'induite $\rho\vert\,\vert^{\zeta_{0}B_{0}}\times \cdots \times \rho\vert\,\vert^{\zeta_{0}A_{0}}$ qui n'est autre que $St(\rho,a_{0})\vert\,\vert^{-(b_{0}-1)/2}$. Et on obtient l'identification de $\pi^+$ avec $\pi^0$ l'unique sous-module irr\'eductible de cette induite.

On r\`egle maintenant le cas o\`u $\zeta_{0}=+$; 

ici le param\`etre $\underline{t}^+$ d\'efinissant $\pi^+$ v\'erifie $\underline{t}^+(\rho,A_{0},B_{0},\zeta_{0})=\underline{t}(\rho,A'_{0},B'_{0},\zeta_{0})+1$ et vaut $1$ si $b_{0}=2$. On d\'efinit encore $\psi_{>}$ sans changer $\psi$ si $b_{0}=2$ et en rempla\c{c}ant $(\rho,A'_{0},B'_{0},\zeta'_{0})$ par $(\rho,A'_{0}+T, B'_{0}+T,\zeta'_{0})$ qui vaut exactement $(\rho,A_{0}+T-1,B_{0}+T+1,\zeta_{0})$. Et on d\'efinit $\psi^+_{>}$ en rempla\c{c}ant $(\rho,A_{0},B_{0},\zeta_{0})$ par $(\rho,A_{0}+T,B_{0}+T,\zeta_{0})$. On d\'efinit alors $\pi_{>}$ (resp. $\pi^+_{>}$) dans le paquet associ\'e \`a $\psi_{>}$ (resp. $\psi^+_{>}$) en utilisant les param\`etres de $\pi$ (resp. $\pi^+$).
Ainsi si la restriction de $\psi^+_{>}$ \`a $W_{F}$ fois  la diagonale de $SL(2,{\mathbb C})\times SL(2,{\mathbb C})$ est sans multiplicit\'e et on a par d\'efinition une inclusion
$$
\pi^+_{>}\hookrightarrow <B_{0}+T, \cdots, -A_{0}-T>_{\rho}\times \pi_{>}.\eqno(1)
$$
Si $\psi^+_{>}$ n'est pas de restriction discr\`ete \`a la diagonale, on v\'erifie, apr\`es avoir fait ''redescendre'' les quadruplets stritement plus petits que $(\rho,A_{0},B_{0},\zeta_{0})$ que l'on garde la m\^eme inclusion. Par d\'efinition:
$$
\pi^+:=\circ_{\ell \in [1,T]}Jac_{(B_{0}+\ell), \cdots, A_{0}+\ell}\pi^+_{>};
$$
$$
\pi:=\circ_{\ell \in [1,T]}Jac_{(B_{0}+1+\ell), \cdots, A_{0}-1+\ell} \pi_{>}.
$$
On sait donc que l'on a une inclusion
$$
\pi_{>} \hookrightarrow \biggl<\begin{matrix} B_{0}+T+1 &\cdots & A_{0}+T-1\\ \vdots &\vdots &\vdots \\
B_{0}+2 &\cdots & A_{0}\end{matrix}\biggr>_{\rho}\times \pi.
$$
On utilise l'inclusion suivante dans un $GL$ convenable:
$$<B_{0}+T, \cdots, -A_{0}-T>_{\rho}
 \hookrightarrow 
<B_{0}+T, \cdots, B_{0}+1>_{\rho}\times <B_{0}, \cdots, -A_{0}>_{\rho}\times $$
$$<-(A_{0}+1), \cdots, -(A_{0}+T)>_{\rho}
$$
et en composant avec (1) et l'inclusion ci-dessus, on obtient: $\pi^+ \hookrightarrow$
$$
<B_{0}+T, \cdots, B_{0}+1>_{\rho}\times <B_{0}, \cdots, -A_{0}>_{\rho}\times <-(A_{0}+1), \cdots, -(A_{0}+T)>_{\rho} \times $$
$$ \biggl<\begin{matrix} B_{0}+T+1 &\cdots & A_{0}+T-1\\ \vdots &\vdots &\vdots \\
B_{0}+2 &\cdots & A_{0}\end{matrix}\biggr>_{\rho}\times \pi$$
$$
\simeq 
<B_{0}+T, \cdots, B_{0}+1>_{\rho}\times \biggl<\begin{matrix} B_{0}+T+1 &\cdots & A_{0}+T-1\\ \vdots &\vdots &\vdots \\
B_{0}+2 &\cdots & A_{0}\end{matrix}\biggr>_{\rho}\times <B_{0}, \cdots, -A_{0}>_{\rho}$$
$$\times <-(A_{0}+1), \cdots, -(A_{0}+T)>_{\rho} \times \pi
$$
On sait aussi que $\rho\vert\,\vert^{x}\times \pi$ est irr\'eductible si $x\in [A_{0}+1,A_{0}+T]$ gr\^ace \`a \ref{irreductibilite} dont les hypoth\`eses sont satisfaites: puisque $\zeta_{0}=+$,  on a suppos\'e (cf. \ref{proprietedelordre}) que tout $(\rho,A,B,\zeta)\in Jord(\psi)$ v\'erifiant $A\geq A_{0}$ est plus grand que $(\rho,A_{0},B_{0},\zeta_{0})$ dans $Jord(\psi^+)$. Il v\'erifie donc aussi $B>>A_{0}+T$. 
D'o\`u $$ <-(A_{0}+1), \cdots, -(A_{0}+T)>_{\rho} \times \pi \simeq <(A_{0}+T), \cdots, (A_{0}+1)>_{\rho}\times \pi.$$
En remontant ci-dessus, on obtient une inclusion de $\pi^+_{>} \hookrightarrow$
$$
<B_{0}+T, \cdots, B_{0}+1>_{\rho}\times \biggl<\begin{matrix} B_{0}+T+1 &\cdots & A_{0}+T-1\\ \vdots &\vdots &\vdots \\
B_{0}+2 &\cdots & A_{0}\end{matrix}\biggr>_{\rho}\times <B_{0}, \cdots, -A_{0}>_{\rho}\times$$
$$ <(A_{0}+T), \cdots, A_{0}+1>_{\rho}\times \pi \simeq 
$$
$$
<B_{0}+T, \cdots, B_{0}+1>_{\rho}\times \biggl<\begin{matrix} B_{0}+T+1 &\cdots & A_{0}+T-1\\ \vdots &\vdots &\vdots \\
B_{0}+2 &\cdots & A_{0}\end{matrix}\biggr>_{\rho}\times <(A_{0}+T), \cdots, (A_{0}+1)>_{\rho}$$
$$\times
 <B_{0}, \cdots, -A_{0}>_{\rho}\times \pi.
 $$
Par r\'eciprocit\'e de Frobenius
$
\pi^+=\circ_{\ell\in [1,T]}Jac_{B_{0}+\ell, \cdots, A_{0}+\ell}\pi^+_{>}$ est non nul et \`a un morphisme non nul dans l'induite $<B_{0}, \cdots, -A_{0}>_{\rho}\times \pi$; comme $\pi^+$ est  irr\'eductibilit\'e, il  est donc un sous-module de cette induite et co\"{\i}ncide avec $\pi^0$ d\'efini en \ref{remarque} par unicit\'e du sous-module irr\'eductible. D'o\`u aussi l'assertion dans ce cas.

On consid\`ere maintenant le cas o\`u $\zeta_{0}=-$ et $B_{0}=1/2$. Par hypoth\`ese $(\rho,A_{0},B_{0},\zeta_{0})$ est le plus petit \'el\'ement de $Jord(\psi^+)$ et quand $b_{0}\neq 2$, $(\rho,A'_{0},B'_{0},\zeta'_{0})$ est le plus petit \'el\'ement de $Jord(\psi)$. Ainsi on sait que les restrictions de  $\psi^+$ et de $\psi$ \`a $W_{F}$ fois  la diagonale de $SL(2,{\mathbb C})\times SL(2,{\mathbb C})$ sont sans multiplicit\'e. On conna\^{\i}t donc $\pi^+$ et $\pi$ comme certain sous-module irr\'eductible. On reprend ces constructions.

On r\`egle d'abord le cas o\`u $b_{0}=2$ et $a_{0}=1$; comme $\underline{\eta}^+(\rho,1/2,1/2,-)=+$, par d\'efinition $\pi^+$ est l'unique sous-module irr\'eductible de l'induite $\rho\vert\,\vert^{-1/2}\times \pi$ et c'est bien la m\^eme d\'efinition que $\pi^0$. D'o\`u le r\'esultat trivialement dans ce cas.

On suppose donc maintenant que $b_{0}>2$ avec toujours $a_{0}=b_{0}-1$. D'o\`u $\zeta'_{0}=+$, $\zeta_{0}=-$, $B_{0}=B'_{0}=1/2$.

On pose $t:=\underline{t}(\rho,A'_{0},B'_{0},\zeta'_{0})$. On fait une r\'ecurrence sur $t$. Pour initialiser la r\'ecurrence, supposons  que $t=0$. Dans ce cas, on sait d'apr\`es \cite{discret} 4.2 (notre $t$ \'etait $\ell$) que $\pi$ est dans le paquet associ\'e au morphisme $\psi'$ tel que $Jord(\psi')$ se d\'eduit de $Jord(\psi)$ en rempla\c{c}ant $(\rho,A'_{0},B'_{0},+)$ par l'ensemble des \'el\'ements $(\rho,C,C,+)$ pour $C\in [A'_{0},B'_{0}]$ avec comme param\`etre $\underline{t}',\underline{\eta}'$ qui se d\'eduit de $\underline{t},\underline{\eta}$, en posant simplement $\underline{t}'(\rho,C,C,+)=0$ (on n'a pas le choix) et $\underline{\eta}'(\rho,C,C,+)=\underline{\eta}(\rho,A'_{0},B'_{0},\zeta'_{0})(-1)^{C-B'_{0}}$. On distingue suivant les valeurs de $\underline{\eta}(\rho,A'_{0},B'_{0},\zeta'_{0})$.

On suppose d'abord que $\underline{\eta}(\rho,A'_{0},B'_{0},\zeta'_{0})=+$. On note $\psi''$ le morphisme qui se d\'eduit de $\psi'$ en rempla\c{c}ant chaque $(\rho,C,C,+)$ comme ci-dessus par $(\rho,C-1,C-1,+)$ sauf $(\rho,1/2,1/2,+)$ qui dispara\^{\i}t. Et on note $\pi'$ la repr\'esentation dans le paquet associ\'e \`a $\psi''$ ayant les param\`etres de $\pi$, apr\`es la translation par $-1$, c'est-\`a-dire que l'on a un signe sur les  $(\rho,C-1,C-1,+)$ pour $C\in ]1/2,A'_{0}]$  qui  alterne en valant maintenant $-$ sur $(\rho,1/2,1/2,+)$. Alors, d'apr\`es les d\'efinitions  (\cite{elementaire}, 2.3 \bf 1\rm pour le calcul de $Jac_{1/2}\pi$ et  \cite{discret} 4.4 pour le calcul $Jac_{3/2, \cdots A'_{0}}Jac_{1/2}\pi=\pi'$),  $\pi$ est l'unique sous-module irr\'eductible de l'induite
$$
<1/2, \cdots, A'_{0}>_{\rho}\times \pi'. \eqno(1)
$$
On note $\psi'''$ le morphisme qui se d\'eduit de $\psi'$ en rempla\c{c}ant chaque $(\rho,C,C,+)$ pour $C\in [1/2,A'_{0}-1]$ en $(\rho,C,C,-)$. 
Et  $\pi'$ est aussi dans le paquet associ\'e \`a $\psi'''$ avec les m\^emes param\`etres que que pour $\psi''$: une r\'ef\'erence pour cela est que l'on peut se ramener au cas \'el\'ementaire (pour tout $(\rho',A',B',\zeta')\in Jord(\psi''), A'=B'$) et appliquer le th\'eor\`eme 5 de \cite{elementaire} avec $a=2A'_{0}-1$; l'involution qui y est d\'efini est l'identit\'e sur $\pi'$ car tous les modules de Jacquet  \`a cause de l'alternance des signes en commen\c{c}ant par $-$ (premi\`ere propri\'et\'e de \cite{elementaire} 2.2, puiqu'ici d'apr\`es ce que l'on vient de rappeler le $b_{\rho,\psi,\epsilon}$ de loc.cite vaut $2A'_{0}-1$) . On reprend la d\'efinition de $\pi^0$ comme unique sous-module irr\'eductible d'une induite (cf. \ref{remarque}) et on \'ecrit la suite d'inclusion, en se rappelant que $A'_{0}=A_{0}-1$
$$
\pi^0\hookrightarrow <-1/2, \cdots, -A_{0}>_{\rho}\times <1/2, \cdots, A_{0}-1>_{\rho}\times \pi'
\hookrightarrow$$
$$ \rho\vert\,\vert^{-1/2}\times <-3/2, \cdots, -A_{0}>_{\rho}\times <1/2, \cdots, A_{0}-1>_{\rho}\times \pi'$$
$$
\simeq 
\rho\vert\,\vert^{-1/2}\times <1/2, \cdots, A_{0}-1>_{\rho}\times <-3/2, \cdots, -A_{0}>_{\rho}\times \pi'.
\hookrightarrow $$
$$
\rho\vert\,\vert^{-1/2}\times <1/2, \cdots, A_{0}-1>_{\rho}\times <-3/2, \cdots, -A_{0}+1>_{\rho}\times \rho\vert\,\vert^{-A_{0}}\times \pi'.
$$
L'induite $\rho\vert\,\vert^{-A_{0}}\times \pi'$ est irr\'eductible, d'apr\`es \ref{irreductibilite} puisque $\pi'$ est dans le paquet associ\'e \`a $\psi''$ et que tout $(\rho,A',B',\zeta')\in Jord(\psi'')$ v\'erifie soit $A'=B'\leq A'_{0}-1=A_{0}-2$ soit $B'>> A_{0}$. D'o\`u encore
$$
\pi^0\hookrightarrow \rho\vert\,\vert^{-1/2}\times <1/2, \cdots, A_{0}-1>_{\rho}\times <-3/2, \cdots, -A_{0}+1>_{\rho}\times \rho\vert\,\vert^{A_{0}}\times \pi' \simeq$$
$$
\rho\vert\,\vert^{-1/2}\times <1/2, \cdots, A_{0}-1>_{\rho}\times \rho\vert\,\vert^{A_{0}}\times <-3/2, \cdots, -A_{0}+1>_{\rho}\times \pi'.
$$
On veut v\'erifier que cette inclusion se factorise par l'unique sous-module irr\'eductible inclut dans l'induite form\'ee par  les 3 premiers facteurs. Pour cela il suffit de montrer que $Jac_{\rho\vert\,\vert^x}\pi^0=0$ pour $x=1/2,A_{0}$. On revient \`a (1) qui donne
$$
\pi^0\hookrightarrow <-1/2, \cdots, -A_{0}>\times <1/2, \cdots, A_{0}-1>_{\rho}\times \pi'\hookrightarrow$$
$$ \rho\vert\,\vert^{-1/2}\times \rho\vert\,\vert^{1/2}\times <-3/2, \cdots, -A_{0}>_{\rho}\times <3/2, A'_{0}>_{\rho}\times \pi'. \eqno(2)
$$
L'inclusion se factorise soit par $<-1/2,1/2>_{\rho}\times \cdots$ soit par $<1/2,-1/2>_{\rho}\times \cdots$. Dans le premier cas o\`u a s\^urement $Jac_{1/2}\pi^0=0$ car $Jac_{1/2}\pi'=0$ et dans le second on a $Jac_{-1/2}\pi^0=0$ car $Jac_{-1/2}\pi'=0$. Donc la 2e possibilit\'e est absurde et c'est la premi\`ere qui est v\'erifi\'ee; d'o\`u
$$
\pi^0\hookrightarrow <-1/2,1/2>_{\rho}\times <-3/2, \cdots, -A_{0}>_{\rho}\times <3/2, A'_{0}>_{\rho}\times \pi'
$$
et $Jac_{1/2}$ de l'induite de droite vaut 0; d'o\`u a fortiori $Jac_{1/2}\pi^0=0$. Montrons maintenant que $Jac_{A_{0}}\pi^0=0$; ici on se rappelle que $A'_{0}=A_{0}-1$. On v\'erifie que l'induite $\rho\vert\,\vert^{A_{0}}\times \pi$ est r\'eductible, elle a un unique sous-module irr\'eductible, $\pi_{s}$,  et un unique quotient irr\'eductible $\pi_{q}$: en effet avec (1)
$$
\rho\vert\,\vert^{A_{0}}\times \pi\hookrightarrow \rho\vert\,\vert^{A_{0}}\times <1/2, \cdots, A_{0}-1>_{\rho}\times \pi'.
$$
L'inclusion se factorise soit par l'unique sous-module irr\'eductible inclus dans $ \rho\vert\,\vert^{A_{0}}\times <1/2, \cdots, A_{0}-1>_{\rho}$, not\'e $\tau$ soit par $<1/2, \cdots, A_{0}>_{\rho}\times \pi'$. On a 
$$
Jac_{A_{0}}<1/2, \cdots, A_{0}>_{\rho}\times \pi'=0
$$
car $Jac_{A_{0}}\pi'=0$ donc un sous-module irr\'eductible de $\rho\vert\,\vert^{A_{0}}\times \pi$ n'est s\^urement pas dans l'induite $<1/2, \cdots, A_{0}>_{\rho}\times \pi'$. A l'inverse $Jac_{-A_{0}}\tau \times \pi'=0$ et tout quotient irr\'eductible de $\rho\vert\,\vert^{A_{0}}\times \pi$ a la propri\'et\'e oppos\'ee. D'o\`u la r\'eductibilit\'e de $\rho\vert\,\vert^{A_{0}}\times \pi$; l'unicit\'e de $\pi_{s}$ et $\pi_{q}$ utilise alors simplement le fait que $Jac_{\pm A_{0}}\pi=0$. De plus $\pi_{q}$ est le seul sous-quotient irr\'eductible de l'induite v\'erifiant $Jac_{-A_{0}}\neq 0$. D'o\`u
$$
\pi^0 \hookrightarrow <-1/2, \cdots, -A_{0}+1>_{\rho}\times \rho\vert\,\vert^{-A_{0}}\times \pi
$$
se factorise par $<-1/2, \cdots, -A_{0}+1>_{\rho}\times \pi_{q}$ et en particulier $Jac_{A_{0}}\pi^0=0$. On revient \`a l'inclusion (2) et on obtient
$$
\pi^0\hookrightarrow <-1/2, 1/2, \cdots, A_{0}-1, A_{0}>_{\rho}\times <-3/2, \cdots, -A_{0}+1>_{\rho}\times \pi'.\eqno(3)
$$
L'induite $<-3/2, \cdots, -A_{0}+1>_{\rho}\times \pi'_{>>}$ a un unique sous-module irr\'eductible (\cite{discret} 4.4) qui est pr\'ecis\'ement la repr\'esentation du paquet associ\'e au morphisme $\tilde{\psi}$ qui s'obtient en rempla\c{c}ant $(\rho,C,C,-)\in Jord(\psi''_{>>})$ pour tout $C\in [1/2,A'_{0}-2]$ par $(\rho,C+1,C+1,-)$ et en gardant les m\^emes signes. On note $\tilde{\pi}$ cette repr\'esentation et c'est la seule sous repr\'esentation irr\'eductible de l'induite qui v\'erifie $Jac_{-3/2, \cdots, -A_{0}+1}\neq 0$. Donc n\'ecessairement (3) se factorise par 
$$
\pi^0\hookrightarrow <-1/2, 1/2, \cdots, A_{0}-1, A_{0}>_{\rho}\times \tilde{\pi}.
$$
Mais $\pi^+$ v\'erifie la m\^eme inclusion et l'induite de droite a un unique sous-module irr\'eductible d'o\`u l'isomorphisme $\pi^0\simeq \pi^+$. 

Le cas o\`u $\underline{\eta}(\rho,A'_{0},B'_{0},\zeta'_{0})=-$ avec toujours $\underline{t}(\rho,A'_{0},B'_{0},\zeta'_{0})=0$ est plus simple: $\pi$ est dans le paquet obtenu en rempla\c{c}ant $(\rho,A'_{0},B'_{0},+)$ par $\cup_{C\in [1/2,A'_{0}]}(\rho,C,C,+)$ avec un signe qui alterne en commen\c{c}ant par $-$; on peut donc encore remplacer les $(\rho,C,C,+)$ par $(\rho,C,C,-)$ (m\^eme argument que ci-dessus pour passer de $\psi''$ \`a $\psi'''$). Par sa d\'efinition, $\pi^+$ est dans le paquet o\`u on remplace dans $Jord(\psi^+)$ le quadruplet $(\rho,A_{0},B_{0},-)$ par $\cup_{C\in [1/2,A_{0}]}(\rho,C,C,-)$ et avec comme param\`etre un signe qui alterne sur ces \'el\'ements en commen\c{c}ant par $+$. D'o\`u directement une inclusion (cf. ci-dessus):
$$
\pi^+\hookrightarrow <-1/2, \cdots, -A_{0}>_{\rho}\times \pi.
$$
Et l'assertion r\'esulte de l'unicit\'e du sous-module irr\'eductible.

On suppose maintenant que $t>0$ et on admet la remarque pour le morphisme $\psi'$ qui s'obtient  en rempla\c{c}ant $(\rho,A'_{0},B'_{0},+)$ par $(\rho,A'_{0}-2,B'_{0},+)$ pour les param\`etres $\underline{t}'$ valant $t-1$ sur cet \'el\'ement. On v\'erifie d'abord qu'il existe $\pi'$ dans le paquet associ\'e \`a ce morphisme avec les m\^emes param\`etres que $\pi$ sauf $t$ remplac\'e par $t-1$: On v\'erifie que l'on a une inclusion 
$$
\pi\hookrightarrow <1/2, \cdots, -A'_{0}>_{\rho}\times <3/2, \cdots, A'_{0}-1>_{\rho} \times \pi';
$$ l'apparition de $<1/2, \cdots, -A'_{0}>_{\rho}$ vient du fait que $t\geq 1$ et la description rappel\'ee dans l'introduction de \cite{discret} et celle de $<3/2, \cdots, A'_{0}-1>_{\rho}$ est encore l'application de \cite{discret} 4.4.

On sait aussi avec \ref{irreductibilite} que $\rho\vert\,\vert^{A'_{0}}\times \pi'$ est irr\'eductible. D'o\`u encore les inclusions
$$
\pi\hookrightarrow <1/2, \cdots, -A'_{0}+1>_{\rho}\times \rho\vert\,\vert^{-A'_{0}}\times  <3/2, \cdots, A'_{0}-1>_{\rho}\times \pi'\simeq$$
$$
\pi\hookrightarrow <1/2, \cdots, -A'_{0}+1>_{\rho}\times <3/2, \cdots, A'_{0}-1>_{\rho}\times  \rho\vert\,\vert^{-A'_{0}}\times \pi'$$
$$
\simeq <1/2, \cdots, -A'_{0}+1>_{\rho}\times <3/2, \cdots, A'_{0}-1>_{\rho}\times  \rho\vert\,\vert^{A'_{0}}\times \pi'
\hookrightarrow$$
$$
\rho\vert\,\vert^{1/2}\times <-1/2, \cdots, -A'_{0}+1>_{\rho}\times <3/2, \cdots, A'_{0}-1>_{\rho}\times  \rho\vert\,\vert^{A'_{0}}\times \pi'$$
$$
\simeq 
\rho\vert\,\vert^{1/2}\times <3/2, \cdots, A'_{0}-1>_{\rho}\times \rho\vert\,\vert^{A'_{0}}\times <-1/2, \cdots, -A'_{0}+1>_{\rho}\times \pi'.
$$
On sait que $Jac_{x}\pi=0$ pour $x=3/2$ et $x=A'_{0}$ car $\pi$ est dans le paquet associ\'e \`a $\psi$ et que $Jord(\psi)$ ne contient pas d'\'el\'ement de la forme $(\rho,A',B',\zeta')$ avec $\zeta'B'=3/2$ ou $\zeta' B=A'_{0}$ (et on applique \ref{proprietedujac} avec $x=y$) d'o\`u encore
$$
\pi\hookrightarrow <1/2, \cdots, A'_{0}>_{\rho}\times <-1/2, \cdots, -A'_{0}+1>_{\rho}\times \pi'.\eqno(4)
$$
On note $\pi^{'+}$ l'unique sous-module irr\'eductible de l'induite $<-1/2, \cdots, -A'_{0}+1>_{\rho}\times \pi'$; on sait, par l'hypoth\`ese de r\'ecurrence que c'est la repr\'esentation associ\'e au morphisme qui se d\'eduit de $\psi'$ en rempla\c{c}ant $(\rho,A'_{0}-2,1/2,+)$ par $(\rho,A'_{0}-1,1/2,-)$ et avec les param\`etres qui sont ceux de $\pi^+$ sauf le $t^+$ qui est $t^+-1$.

Comme $Jac_{1/2, -1/2, \cdots, -A'_{0}+1}\pi\neq 0$, l'inclusion (4) se factorise n\'ecessairement par
$$
\pi\hookrightarrow <1/2, \cdots, A'_{0}>_{\rho}\times \pi^{'+}.$$
D'o\`u, en se rappelant que $A_{0}=A'_{0}+1$:
$$
\pi^0\hookrightarrow <-1/2, \cdots, -A_{0}>_{\rho}\times <1/2, \cdots, A_{0}-1>_{\rho}\times \pi^{'+} \hookrightarrow$$
$$
\rho\vert\,\vert^{-1/2}\times <-3/2, \cdots, -A_{0}+1>_{\rho}\times \rho\vert\,\vert^{-A_{0}}\times <1/2, \cdots, A_{0}-1>_{\rho}\times \pi^{'+}.
$$
L'induite $\rho\vert\,\vert^{A_{0}}\times \pi^{'+}$ est irr\'eductible (cf. \ref{irreductibilite}) car $Jord(\psi^{'+})$ contient $(\rho, A'_{0}-1,1/2,-)$ avec $A'_{0}=A_{0}-1$ et d'autres termes qui n'ont rien \`a voir avec $A_{0}$; d'o\`u:
$$
\pi^0\hookrightarrow \rho\vert\,\vert^{1/2}\times <-1/2, \cdots, -A_{0}+1>_{\rho}\times \rho\vert\,\vert^{-A_{0}}\times 
<-3/2, \cdots, -A_{0}+1>_{\rho}\times \pi^{'+}.
$$
Comme dans la preuve o\`u $t=0$, on v\'erifie que cette inclusion se factorise par
$$
\pi^0\hookrightarrow <1/2, -1/2, \cdots, -A_{0}+1,-A_{0}>_{\rho}\times <-3/2, \cdots, -A_{0}+1>_{\rho}\times \pi^{'+}.
$$
La repr\'esentation $\pi^+$ a la m\^eme propri\'et\'e et l'induite de droite a un unique sous-module irr\'eductible. D'o\`u $\pi^+\simeq \pi^0$ comme annonc\'e.
\subsubsection{Descente dans le cas holomorphe\label{identification}}
Ici on suppose simplement que $r(s,\rho,a_0,\psi)$ n'a pas de p\^ole.
\begin{lem} La repr\'esentation $\pi^+$ est non nulle et est l'unique sous-module irr\'eductible de l'induite $St(\rho,a_{0})\vert\,\vert^{-(b_{0}-1)/2}\times \pi$.
\end{lem}
On fixe $(\rho,A,B,\zeta)\in Jord(\psi)$ strictement sup\'erieur \`a $(\rho,A'_{0},B'_{0},\zeta_{0})$ et on suppose que pour tout $(\rho,A',B',\zeta')\in Jord(\psi)$ strictement sup\'erieur \`a $(\rho,A,B,\zeta)$ on a $B'>>A$. On note $\psi'$ le morphisme qui se d\'eduit de $\psi$ en rempla\c{c}ant $(\rho,A,B,\zeta)$ par $(\rho,A+1,B+1,\zeta)$ et $\pi'$ la repr\'esentation dans le paquet associ\'e \`a $\psi'$ ayant les m\^emes param\`etres que $\pi$. C'est-\`a-dire
$
\pi=Jac_{\zeta(B+1), \cdots, \zeta (A+1)}\pi'
$ ou encore
$$
\pi'\hookrightarrow <\zeta(B+1), \cdots, \zeta(A+1)>_{\rho}\times \pi .
$$
On admet le lemme pour $\psi',\pi'$ et on le d\'emontre pour $\psi,\pi$. Il est clair que cela permet de redescendre du cas d\'ej\`a d\'emontr\'e en \ref{tresdominant} au cas cherch\'e. On note $\pi^{'+}$ l'analogue de $\pi^+$ et on sait que $\pi^{'+}$ est l'unique sous-module irr\'eductible de l'induite $$\pi^{'+}\hookrightarrow <\zeta_{0}(B_{0}, \cdots, -A_{0})>_{\rho}\times \pi'.$$
En composant avec l'inclusion reliant $\pi'$ et $\pi$, on obtient
$$
\pi^{'+}\hookrightarrow <\zeta_{0}B_{0}, \cdots, -A_{0}>_{\rho}\times <\zeta(B+1), \cdots, \zeta(A+1)>_{\rho}\times \pi.
$$
On va v\'erifier que l'induite $<\zeta_{0}B_{0}, \cdots, -A_{0}>_{\rho}\times <\zeta(B+1), \cdots, \zeta(A+1)>_{\rho}$ dans le $GL$ convenable est irr\'eductible. C'est \'evident si $\zeta_{0}=-$ et $\zeta=+$ car avec $\zeta_{0}=-$ on a suppos\'e que $B_{0}>0$. On consid\`ere donc les autres cas et d'abord celui o\`u $A<A_{0}$; si $\zeta_{0}=-$ on a aussi $\zeta=-$ et par la propri\'et\'e ${\mathcal P}$ de l'ordre sur $Jord(\psi^+)$, $B\geq B_{0}$ puisque $(\rho,A,B,\zeta)>(\rho,A_{0},B_{0},\zeta)$. Ainsi le segment $[-B,-A]$ est inclus dans $[\zeta_{0}B_{0},-A_{0}]$. Si $\zeta_{0}=+$ et $\zeta=-$ on a encore l'inclusion du segment $[-B,-A]$ dans le segment $[\zeta_{0}B_{0},-A_{0}]$.  D'o\`u l'irr\'eductibilit\'e dans ces cas.  Si $\zeta=\zeta_{0}=+$, comme $(\rho,A,B,\zeta)>(\rho,A'_{0},B'_{0},\zeta_{0})$, la condition (0) de \ref{proprietedelordre} assure que $B>B_{0}+1$; d'o\`u  $B+1-B_{0}\geq 2$ ce qui assure encore l'irr\'eductibilit\'e.

Supposons  que $A\geq A_{0}$. Par hypoth\`ese la fonction $r(s,\rho,a_0,\psi)$ n'a pas de p\^ole en $s=(b_{0}-1)/2$.

Si $\zeta_{0}=+$, on a donc, d'apr\`es \ref{poles}, $B>B_{0}$ si $\zeta=+$ et $B>A_{0}$ si $\zeta=-$; d'o\`u si $\zeta=+$, $(B+1)-B_{0}\geq 2$ et si $\zeta=-$, $-A_{0}+(B+1)\geq 2$ et encore l'irr\'eductibilit\'e.

Si $\zeta_{0}=\zeta=-$, l'hypoth\`ese d'holomorphie entra\^{\i}ne que soit $B_{0}>B$ soit $B>A_{0}$; dans le premier cas, le segment $[-B_{0},-A_{0}]$ est inclus dans le segment $[-B,-A]$ et dans le 2e cas, les segments ne sont pas li\'es. On a donc encore l'irr\'eductibilit\'e de l'induite \'ecrite pr\'ec\'edemment. 

\

D'o\`u, par \'echange des facteurs:
$$
\pi^{'+}\hookrightarrow <\zeta(B+1), \cdots, \zeta(A+1)>_{\rho}\times <\zeta_{0}(B_{0}, \cdots, -A_{0})>_{\rho}\times \pi.
$$
Or $\pi^+=Jac_{\zeta(B+1), \cdots, \zeta (A+1)}\pi^{'+}$ (car  par d\'efinition des param\`etres, $\pi^+$ se d\'eduit de $\pi^{'+}$ de la m\^eme fa\c{c}on que $\pi$ se d\'eduit de $\pi'$). Par r\'eciprocit\'e de Frobenius appliqu\'e \`a l'induite ci-dessus, $\pi^+$ s'envoie de fa\c{c}on non nulle dans l'induite $<\zeta_{0}(B_{0}, \cdots, -A_{0})>_{\rho}\times \pi$. Ceci est bien l'assertion cherch\'ee et cela termine la preuve.

\section{Image des op\'erateurs d'entrelacement dans le cas de bonne parit\'e}

\subsection{Le cas holomorphe\label{imagesimple}}
\begin{lem}On suppose que $r(s,\rho,a_0,\psi)$ n'a pas de p\^ole en $s=(b_{0}-1)/2$. L'image de l'op\'erateur $N_{\psi}(s,\rho,a_{0},\pi)$ en $s=(b_{0}-1)/2$ est irr\'eductible et non nulle; c'est la repr\'esenta\-tion $\pi^+$ d\'ecrite en \ref{descriptionparametre}.
\end{lem}
L'hypoth\`ese  assure que  l'op\'erateur d'entrelacement standard $M(s,\rho,a_{0},\pi)$ est  comme $N_{\psi}(s,\rho,a_{0},\pi)$ holomorphe en $s=s_{0}$ et qu'en ce point les 2 op\'erateurs co\"{\i}ncident \`a un scalaire non nul pr\`es. Or un op\'erateur d'entrelacement standard, en un point o\`u il est holomorphe y est certainement non nul. D'o\`u la non nullit\'e annonc\'ee dans l'\'enonc\'e et l'image contient certainement $\pi^+$ comme sous-module irr\'eductible. Il faut donc d\'emontrer l'irr\'eductibilit\'e de l'image.

Avec \cite{mvw} II.1, on sait qu'il existe un automorphisme explicite (venant du groupe des similitudes) tel que toute repr\'esentation irr\'eductible du groupe consid\'er\'e est isomorphe \`a sa repr\'esentation duale transform\'ee par cet automorphisme. Ainsi en dualisant et en appliquant cet automorphisme, on voit que $\pi^+$ est aussi un quotient irr\'eductible de l'induite $<A_{0}, \cdots, -\zeta_{0}B_{0}>_{\rho}\times \pi$. Par multiplicit\'e 1 de $\pi^+$ comme sous-quotient irr\'eductible de l'induite \'ecrite (cf. \ref{remarque} et \ref{identification}), $N_{\psi}(s_{0},\pi,\rho,a_{0})$  n'annulle pas ce quotient et l'envoie dans l'unique sous-module irr\'eductible de l'induite $<\zeta_{0}B_{0}, \cdots, -A_{0}>_{\rho}\times \pi$. Ainsi $\pi^+$ est facteur direct de l'image mais par unicit\'e du sous-module irr\'eductible, sous nos hypoth\`eses, l'image est r\'eduite \`a $\pi^+$. Cela termine la preuve.
\subsection{Descente \`a partir du cas holomorphe\label{imagebonneparite}}
On a fix\'e $\psi, (\rho,a_{0},b_{0})$ et on fixe aussi un morphisme $\psi$ et  $\pi$ dans le paquet de repr\'esentations associ\'e \`a $\psi$; on suppose que $\psi$ est de bonne parit\'e (cf. \ref{descriptionpaquets} pour cette d\'efinition). 
\begin{prop} (i) On suppose que $(\rho,a_{0},b_{0})$ est de bonne parit\'e et que soit $b_{0}=2$, soit $(\rho,a_{0},b_{0}-2)\in Jord(\psi)$ et $\psi^+$ ainsi que $\pi^+$ sont d\'efinis. Alors l'image de l'op\'erateur d'entrelacement $N_{\psi}(s,\rho,a_{0},\pi)$ en $s=(b_{0}-1)/2$ est exactement $\pi^+$ c'est-\`a-dire vaut 0 si $\pi^+=0$ et est irr\'eductible isomorphe \`a $\pi^+$ sinon.

(ii) Sans hypoth\`ese sur $(\rho,a_{0},b_{0})$, l'image de $N_{\psi}(s,\rho,a_{0},\pi)$ en $s=(b_{0}-1)/2$ est soit nulle soit irr\'eductible.
\end{prop}

(i) Ici $\psi^+$ est d\'efini, on pose encore $s_{0}=(b_{0}-1)/2$.
On reprend les notations $(\rho,A_{0},B_{0},\zeta_{0})$, o\`u $\zeta_{0}$ est le signe de $(a_{0}-b_{0})$, $+$ si ce nombre est $0$ et o\`u $A_{0}=(a_{0}+b_{0})/2-1$, $B_{0}=\vert(a_{0}-b_{0})\vert/2$. On fixe $(\rho,A,B,\zeta)\in Jord(\psi)$ et on suppose que $(\rho,A,B,\zeta)>(\rho,A_{0},B_{0},\zeta_{0})$ dans $Jord(\psi^+)$.

 On suppose aussi que pour tout $(\rho,A',B',\zeta')$ dans $Jord(\psi)$ strictement sup\'erieur \`a $(\rho,A,B,\zeta)$ on a $B'>>A$. On d\'efinit $\psi_{>}$ comme le morphisme qui s'obtient \`a partir de $\psi$ en changeant $(\rho,A,B,\zeta)$ en $(\rho,A+1,B+1,\zeta)$. Les hypoth\`eses assurent que $\psi_{>}$ domine $\psi$ et en reprenant la construction des \'el\'ements de $\psi$ on peut partir d'un morphisme qui domine \`a la fois $\psi_{>}$ et $\psi$ et avant de construire les repr\'esentations associ\'ees \`a $\psi$ on construit celles associ\'ees \`a $\psi_{>}$ et on obtient celles qui sont associ\'ees \`a $\psi$ en prenant $Jac_{\zeta (B+1), \cdots, \zeta(A+1)}$ de celles associ\'ees \`a $\psi_{>}$. Ainsi il existe une unique repr\'esentation $\pi_{>}$ dans le paquet associ\'ee \`a $\psi_{>}$ tel que  $\pi=Jac_{\zeta(B+1), \cdots, \zeta(A+1)}\pi$, ou encore, $\pi_{>}$ est l'unique sous-module irr\'eductible de l'induite
$$
\pi_{>}\hookrightarrow <\zeta(B+1), \cdots, \zeta(A+1)>_{\rho}\times \pi. \eqno(1)
$$
Admettons le r\'esultat pour $\psi_{>},(\rho,a_{0},b_{0})$ et $\pi_{>}$ et montrons le pour $\pi$. Avant de faire cela remarquons que de cette fa\c{c}on on d\'emontrera effectivement le th\'eor\`eme sans les hypoth\`eses $>>$ car pour arriver \`a $\pi$, on part du cas holomorphe d\'ej\`a trait\'e et on redescend par des \'etapes comme celles d\'ecrites ci-dessus. De plus, comme le cas holomorphe est d\'ej\`a connu, on peut supposer qu'il existe $(\rho,A',B',\zeta')\leq (\rho,A,B,\zeta)$ (il peut y avoir \'egalit\'e) qui participe aux p\^oles de $r(s,\rho,a_0,\psi)$ en $s=(b_{0}-1)/2$. Avec la 2e condition de ${\mathcal P}_{p}$ dans \ref{proprietedelordre} cela permet de supposer que $A\geq A_{0}$, hypoth\`ese que nous faisons dans tout ce qui suit.

\

La preuve de (ii) est partiellement commune avec celle de (i); pour pouvoir traiter les 2 en m\^eme temps, il suffit de remarquer que si la fonction $r(s,\rho,a_0,\psi)$ est holomorphe en $s=s_{0}$, l'op\'erateur $N_{\psi}(\rho,a_{0},\pi,s)$ co\"{\i}ncide en $s=s_{0}$ \`a un scalaire pr\`es avec l'op\'erateur d'entrelacement standard; ainsi l'op\'erateur d'entrelacement standard est holomorphe et donc n\'ecessairement  non nul; avec l'argument de \ref{identification} son image est l'unique sous-module irr\'eductible (cf. \ref{remarque}) de $St(\rho,a_{0})\vert\,\vert^{-(b_{0}-1)/2}\times \pi$ et en particulier est irr\'eductible. Donc (ii) est aussi d\'emontr\'e dans le cas holomorphe et m\^eme si $(\rho,a_{0},b_{0}-2)\not \in Jord(\psi)$, on continue la preuve en supposant simplement que $A\geq A_{0}$; $\pi_{>}$ est bien d\'efini comme ci-dessus.

\

On va \'etudier le diagramme ci-dessous, dont les fl\`eches, qui ne sont pas les inclusions se d\'eduisant de (1), sont des op\'erateurs d'entrelacement normalis\'es uniquement d\'etermin\'es \`a une fonction holomorphe inversible pr\`es par le choix d'un \'el\'ement du groupe de Weyl qui autorise les fl\`eches; les doubles classes de ces \'el\'ements du groupe de Weyl sont uniquement d\'etermin\'ees pour $s$ g\'en\'eral et donc pour tout $s$.

$$
\begin{matrix}St(\rho,a_{0})\vert\,\vert^{s}\times \pi_{>}&\hookrightarrow &St(\rho,a_{0})\vert\,\vert^s\times <\zeta(B+1), \cdots, \zeta(A+1)>_{\rho}\times \pi
\\ 
&&\downarrow \\&&<\zeta(B+1), \cdots, \zeta(A+1)>_{\rho}\times St(\rho,a_{0})\vert\,\vert^{s}\times \pi\\ \downarrow
&&\downarrow
\\  & &<\zeta(B+1), \cdots, \zeta(A+1)>_{\rho}\times St(\rho,a_{0})\vert\,\vert^{-s}\times \pi\\
&& \uparrow
\\
St(\rho,a_{0})\vert\,\vert^{-s}\times \pi_{>}&\hookrightarrow &St(\rho,a_{0})\vert\,\vert^{-s}\times <\zeta(B+1), \cdots, \zeta(A+1)>_{\rho}\times \pi
\end{matrix}
$$
Ce diagramme commute \`a une fonction m\'eromorphe pr\`es: en effet il faut comparer le compos\'e d'applications suivantes:$$\begin{matrix}
St(\rho,a_{0})\vert\,\vert^{s}\times <\zeta(B+1),\cdots, \zeta(A+1)>_{\rho}\times \pi  \\ \downarrow
\\
<\zeta(B+1),\cdots, \zeta(A+1)>_{\rho}\times St(\rho,a_{0})\vert\,\vert^{s}\times \pi
\\
\downarrow \\
<\zeta(B+1),\cdots, \zeta(A+1)>_{\rho}\times St(\rho,a_{0})\vert\,\vert^{-s}\times \pi
\end{matrix}
$$
au compos\'e des applications
$$
\begin{matrix}
St(\rho,a_{0})\vert\,\vert^{s}\times <\zeta(B+1),\cdots, \zeta(A+1)>_{\rho}\times \pi  \\ \downarrow
\\
St(\rho,a_{0})\vert\,\vert^{-s}\times <\zeta(B+1),\cdots, \zeta(A+1)>_{\rho}\times \pi 
\\ \downarrow\\
<\zeta(B+1),\cdots, \zeta(A+1)>_{\rho}\times St(\rho,a_{0})\vert\,\vert^{-s}\times \pi.
\end{matrix}
$$ 
Or ces compos\'es correspondent \`a 2 d\'ecompositions d'un \'el\'ement du groupe de Weyl en produit d'autres \'el\'ements mais la premi\`ere d\'ecomposition se fait avec ajout des longueurs ce qui n'est pas le cas de la deuxi\`eme. Ces 2 compos\'es sont donc \'egaux (\`a une fonction holomorphe inversible pr\`es) \`a condition de multiplier le premier compos\'e par la fonction m\'eromorphe produit des op\'erateurs d'entrelacement (dans un groupe $GL$ convenable):
$$
<\zeta(B+1), \cdots, \zeta (A+1)>_{\rho}\times St(\rho,a_{0})\vert\,\vert^{-s}\rightarrow St(\rho,a_{0})\vert\,\vert^{-s}\times <\zeta (B+1), \cdots, \zeta (A+1)>_{\rho}$$
$$\rightarrow 
<\zeta(B+1), \cdots, \zeta (A+1)>_{\rho}\times St(\rho,a_{0})\vert\,\vert^{-s}.
$$
On note $\Delta(s)$ cette fonction m\'eromophe, c'est le produit des facteurs de normalisation \`a la Langlands-Shahidi de \cite{shahidi}; gr\^ace \`a \cite{jpss} th. 8.2, on sait la calculer. Elle vaut, \`a une fonction holomorphe inversible pr\`es

si $\zeta=+$,
$$
\frac{L(St(\rho,a_{0})\times \rho, B+1+s)}{L(St(\rho,a_{0})\times \rho, A+2+s)}\frac{L(St(\rho,a_{0})\times \rho, -s-(A+1))}{L(St(\rho,a_{0})\times \rho, -s-B)},$$
ou encore, en posant $s=(b_{0}-1)/2+s'$
$$\frac
{L(\rho\times \rho, (a_{0}-1)+(b_{0}-1)/2+B+1+s')}{L(\rho\times \rho, (a_{0}-1)/2+(b_{0}-1)/2+A+2+s')}$$
$$
\times \frac
{L(\rho\times \rho, (a_{0}-1)/2-(b_{0}-1)/2-(A+1)-s')}{L(\rho\times \rho, (a_{0}-1)/2-(b_{0}-1)/2-B-s')}=
$$
$$\frac
{L(\rho\times \rho,A_{0}+B+1+s')}{L(\rho\times \rho,A_{0}+A+2+s')}\frac{L(\rho\times \rho, \zeta_{0}B_{0}-(A+1)-s')}{L(\rho\times \rho, \zeta_{0}B_{0}-B-s')}.$$
Quel que soit la valeur de $\zeta_{0}$, on a s\^urement $\zeta_{0}B_{0}-A-1<0$ car $A\geq A_{0}\geq  B_{0}$ et l'ordre de la fonction $\Delta(s)$ en $s'=0$ est 1 si $\zeta_{0}B_{0}=B$ et $0$ sinon. La fonction est donc holomorphe inversible en $s'=0$ sauf exactement si $\zeta_{0}=+$ et $B_{0}=B$ o\`u  elle a un z\'ero simple. On remarque pour la suite que l'ordre de $\Delta(s)$ en $s=s_{0}$ est exactement l'ordre de la fonction $r(s,\rho,a_0,\psi)/r(\psi_{>},\rho,a_{0},s)$: en effet ces 2 fonctions ont m\^eme ordre sauf si $(\rho,A+1,B+1,\zeta)$ participe aux p\^oles du d\'enominateur sans que $(\rho,A,B,\zeta)$ participe aux p\^oles du num\'erateur ou vice et versa. Si $\zeta_{0}=-$, comme $\zeta=+$ aucun des 2 ne participe aux p\^oles et on a le r\'esultat. Si $\zeta_{0}=+$ on a  par hypoth\`ese $A\geq A_{0}$ donc le seul cas o\`u la fonction quotient n'est pas d'ordre 0 est celui o\`u $B_{0}\geq B$ sans avoir $B_{0}\geq B+1$ o\`u elle a un z\'ero simple. Et c'est bien le cas $B_{0}=B$.

Si $\zeta=-$, par un calcul analogue, on trouve que $\Delta(s)$ \`a une fonction holomorphe inversible pr\`es vaut:
$$\frac{
L(\rho\times \rho, A_{0}-(A+1)+s')}{L(\rho\times \rho, A_{0}-B+s')}\frac{L(\rho\times \rho, \zeta_{0}B_{0}+B+1-s')}{L(\rho\times \rho, \zeta_{0}B_{0}+A+2-s')}.
$$

En $s'=0$, cette fonction est holomorphe inversible sauf exactement si:

$
\zeta_{0}=+$ et $A_{0}=B$ o\`u elle a un z\'ero d'ordre 1;

$\zeta_{0}=-$ et $B_{0}=B+1$ o\`u elle a un p\^ole d'ordre 1, ou $A_{0}=B$ o\`u elle a un 0 d'ordre 1.

On v\'erifie ici aussi que $r(s,\rho,a_0,\psi)/r(\psi_{>},\rho,a_{0},s)$ a exactement le m\^eme ordre que $\Delta(s)$ en $s'=s-(b_{0}-1)/2=0$: en effet si $\zeta_{0}=+$ et $\zeta=-$, on a s\^urement $A\geq A_{0}$ et $(\rho,A,B,\zeta)$ contribue aux p\^ole de $r(s,\rho,a_0,\psi)$ en $s=s_{0}$ si et seulement si $B\leq A_{0}$. Et on a un r\'esultat analogue pour $(\rho,A+1,B+1,\zeta)$ et $r(s,\rho,a_0,\psi)$ donc la seul diff\'erence d'ordre se produit quand $B=B_{0}$ et on a un z\'ero simple. Supposons que $\zeta_{0}=-$; comme $\zeta=-$, on a s\^urement aussi $A\geq A_{0}$ et  $(\rho,A,B,\zeta)$ contribue aux p\^oles de $r(s,\rho,a_0,\psi)$ en $s=s_{0}$, exactement quand $B\in [B_{0},A_{0}]$. Et $(\rho,A+1,B+1,\zeta)$ contribue aux p\^oles de $r(\psi_{>},\rho,a_{0},s)$, en $s=s_{0}$, si $B+1\in [B_{0},A_{0}]$. Donc si $B+1=B_{0}$, $r(\psi_{>},\rho,a_{0},s)/r(s,\rho,a_0,\psi)$ a un p\^ole simple et si $B=A_{0}$ ce quotient \`a un z\'ero simple. D'o\`u le r\'esultat annonc\'e.

On obtient donc un diagramme commutatif, \`a une fonction holomorphe inversible pr\`es en $s=s_{0}$:
$$
\begin{matrix}
St(\rho,a_{0})\vert\,\vert^{s}\times  \pi_{>} &\rightarrow &<\zeta(B+1),\cdots, \zeta(A+1)>_{\rho}\times St(\rho,a_{0})\vert\,\vert^{s}\times \pi
\\
 N_{\psi_{>}}(s,\rho,a_{0},\pi_{>})\, \downarrow &&N_{\psi}(s,\rho,a_{0},\pi)\, \downarrow\\
St(\rho,a_{0})\vert\,\vert^{-s}\times \pi_{>} &\rightarrow
&<\zeta(B+1),\cdots, \zeta(A+1)>_{\rho}\times St(\rho,a_{0})\vert\,\vert^{-s}\times \pi.
\end{matrix}
$$ 
On note $H_{1}(s)$ et $H_{2}$ les morphismes horizontaux. Ils sont le compos\'e de l'inclusion (1) de $\pi_{>}$ dans $ <\zeta(B+1),\cdots, \zeta(A+1)>_{\rho}\times \pi$ suivi pour $H_{1}(s)$ de l'op\'erateur d'entrelacement standard:
$$
M_{1}(s):=St(\rho,a_{0})\vert\,\vert^{s}\times <\zeta(B+1),\cdots, \zeta(A+1)>_{\rho}\times \pi \rightarrow$$
$$
<\zeta(B+1),\cdots, \zeta(A+1)>_{\rho}\times St(\rho,a_{0})\vert\,\vert^{s}\times <\zeta(B+1),\cdots, \zeta(A+1)>_{\rho}\times \pi$$
et pour $H_{2}(s)$ de l'op\'erateur d'entrelacement standard
$$
M_{2}(s):=St(\rho,a_{0})\vert\,\vert^{-s}\times <\zeta(B+1),\cdots, \zeta(A+1)>_{\rho}\times \pi \rightarrow$$
$$
<\zeta(B+1),\cdots, \zeta(A+1)>_{\rho}\times St(\rho,a_{0})\vert\,\vert^{-s}\times <\zeta(B+1),\cdots, \zeta(A+1)>_{\rho}\times \pi.$$

Admettons que $M_{1}(s)$ et $M_{2}(s)$ sont holomorphes en $s=s_{0}$, ce que l'on d\'emontrera \`a la fin de la preuve. On v\'erifie l'\'egalit\'e suivante:
$$
Jac_{\zeta (B+1), \cdots, \zeta (A+1)}Im\, N_{\psi}(s_{0},\rho,a_{0},\pi)\circ H_{1}(s_{0})= N_{\psi}(s,\rho,a_{0},\pi)(St(\rho,a_{0})\vert\,\vert^{s_{0}}\times \pi). \eqno(2)
$$
En effet, on montre d'abord que pour tout $x\in [\zeta(B+1),\zeta(A+1)]$, on a $$Jac_{x, \cdots, \zeta(A+1)}St(\rho,a_{0})\vert\,\vert^{s_{0}}\times \pi=0.\eqno(2)'$$ En effet, supposons qu'il n'en soit pas ainsi pour un $x$ bien choisi. On rappelle que  $St(\rho,a_{0})\vert\,\vert^{s_{0}}=<A_{0}, \cdots, -\zeta_{0}B_{0}>_{\rho}$. On a $A_{0}<A+1$, par hypoth\`ese et il existe donc n\'ecessairement $ y\in [\zeta (B+1),\zeta (A+1)]$ avec $Jac_{y, \cdots, \zeta (A+1)}\pi\neq 0$. Mais  ceci est impossible avec \ref{proprietedujac} puisque $Jord(\psi)$ ne contient pas d'\'el\'ement $(\rho,A',B',\zeta')$ avec $\zeta'=\zeta$ et $B'\geq B+1>B$ et $A'\geq A+1>A$ car un tel \'el\'ement est s\^urement strictement sup\'erieur \`a $(\rho,A,B,\zeta)$ et v\'erifie donc, d'apr\`es nos hypoth\`eses g\'en\'erales, $B'>>A$.

On a donc montr\'e que $$Jac_{\zeta(B+1), \cdots, \zeta (A+1)} N_{\psi}(s,\rho,a_{0},\pi)\biggl(<\zeta(B+1), \cdots, \zeta(A+1)>_{\rho}\times St(\rho,a_{0})\vert\,\vert^{s_{0}}\times \pi\biggr)
$$
$$
=N_{\psi}(s,\rho,a_{0},\pi)\biggl( St(\rho,a_{0})\vert\,\vert^{s_{0}}\times \pi\biggr). \eqno(3)
$$
On v\'erifie aussi que 
$$
Jac_{\zeta (B+1), \cdots \zeta (A+1)} Im\, M_{1}(s_{0})=St(\rho,a_{0})\vert\,\vert^{s_{0}}\times \pi:
$$en effet dans l'op\'erateur d'entrelacement,
il n'y a que le groupe $GL$ portant les repr\'esentations induite $St(\rho,a_{0})\vert\,\vert^{s_{0}}\times <\zeta (B+1), \cdots, \zeta (A+1)>_{\rho}$ qui intervient et $M_{1}(s_{0})$ \'echange les 2 facteurs. Donc en tenant compte de (2)' l'\'egalit\'e r\'esulte du calcul des modules de Jacquet des induites.

D'o\`u, par exactitude du foncteur de Jacquet, $Jac_{\zeta (B+1), \cdots, \zeta (A+1)}Im\, N_{\psi}(s,\rho,a_{0},\pi)\circ M_{1}(s_{0})$ vaut (2).

On v\'erifie par les formules calculant le module de Jacquet des induites que $$Jac_{\zeta(B+1), \cdots , \zeta(A+1)}St(\rho,a_{0})\vert\,\vert^s\times \pi_{>}$$ contient certainement l'induite $St(\rho,a_{0})\vert\,\vert^{s_{0}}\times Jac_{\zeta (B+1), \cdots, \zeta (A+1)}\pi_{>}$ vu comme \'el\'ement du bon groupe de Grothendieck, c'est-\`a-dire l'induite $St(\rho,a_{0})\vert\,\vert^{s_{0}}\times \pi$. Donc par exactitude du foncteur de Jacquet:
$$
Jac_{\zeta (B+1), \cdots, \zeta (A+1)}\biggl(St(\rho,a_{0})\vert\,\vert^{s_{0}}\times <\zeta (B+1), \cdots, \zeta (A+1))>_{\rho}\times \pi/St(\rho,a_{0})\vert\,\vert^{s_{0}}\times \pi_{>}\biggr)=0.
$$
Et cela prouve (2), encore par exactitude du foncteur de Jacquet.

La conclusion est que, sous l'hypoth\`ese que $M_{1}(s_{0})$ et $M_{2}(s_{0})$ sont holomorphes, on a montr\'e, par la commutativit\'e du diagramme:
$$
N_{\psi}(s_{0},\rho,a_{0},\pi)\biggl(St(\rho,a_{0})\vert\,\vert^{s_{0}}\times \pi\biggr)= 
$$
$$
Jac_{\zeta(B+1), \cdots, \zeta (A+1)}M_{2}(s_{0})N_{\psi_{>}}(s_{0},\rho,a_{0},\pi_{>})\biggl(St(\rho,a_{0})\vert\,\vert^{s_{0}}\times \pi_{>}\biggr).
$$
On montre comme ci-dessus que l'on peut enlever $M_{2}(s_{0})$ dans la formule pr\'ec\'edente et finalement $$
N_{\psi}(s_{0},\rho,a_{0},\pi)\biggl(St(\rho,a_{0})\vert\,\vert^{s_{0}}\times \pi\biggr)=$$
$$
Jac_{\zeta(B+1), \cdots, \zeta (A+1)}N_{\psi_{>}}(s_{0},\rho,a_{0},\pi_{>})\biggl(St(\rho,a_{0})\vert\,\vert^{s_{0}}\times \pi_{>}\biggr).\eqno(4)$$

Concluons pour obtenir (i) et (ii) sous  les hypoth\`eses d'holomorphie faites pour $M_{i}(s)$ pour $i=1,2$: si $\psi^+$ est d\'efini comme dans (i), alors $\psi^+_{>}$ est aussi d\'efini. On a donc d\'efini $\pi^+$ et $\pi^+_{>}$. Par d\'efinition $\pi^+=Jac_{\zeta (B+1), \cdots, \zeta (A+1)}\pi^+_{>}$. Donc si le th\'eor\`eme est vrai pour $N_{\psi_{>}}(s_{0},\rho,a_{0},\pi_{>})$, le membre de droite de (4) est exactement $\pi^+$ et (i)  est donc vrai pour $N_{\psi}(s_{0},\rho,a_{0},\pi)$ et cela permet de terminer la r\'ecurrence de (i).

Supposons maintenant que $(\rho,a_{0},b_{0}-2)\notin Jord(\psi)$; on veut, dans ce cas, simplement montrer que le membre de gauche de (4) est nul ou est une repr\'esentation irr\'eductible, en admettant le m\^eme r\'esultat pour
$$
\tau:=N_{\psi_{>}}(s_{0},\rho,a_{0},\pi_{>})\biggl(St(\rho,a_{0})\vert\,\vert^{s_{0}}\times \pi_{>}\biggr).
$$
Si $\tau=0$, il n'y a \'evidemment pas de difficult\'e. Supposons qu'il n'en soit pas ainsi et il faut montrer que $Jac_{\zeta (B+1), \cdots, \zeta (A+1)}\tau$ est soit nul soit est une repr\'esentation irr\'eductible. On suppose que le module de Jacquet est non nul et par r\'eciprocit\'e de Frobenius on sait qu'il existe une repr\'esentation irr\'eductible convenable $\sigma$ et une inclusion 
$$
\tau \hookrightarrow \times_{c\in [\zeta (B+1),\zeta (A+1)]}\rho\vert\,\vert^{c}\times \sigma.\eqno(5)
$$
On sait aussi que $\tau$ est un sous-module de $St(\rho,a_{0})\vert\,\vert^{s_{0}}\times \pi_{>}$. Donc en particulier pour tout $y\in ]\zeta (B+1),\zeta (A+1)]$ la non nullit\'e de $Jac_{y, \cdots, \zeta (A+1)}\tau$ entra\^{\i}ne la non nullit\'e, pour un bon choix de  $y'\in [y, \zeta (A+1)]$ de $Jac_{y', \cdots, \zeta (A+1)}\pi_{>}$ (par le fait que $A\geq A_{0}$). Ceci est impossible car $\zeta y'>B+1$ (cf. \ref{proprietedujac}). Cela montre que (5) se factorise n\'ecessairement via un sous-module
$$
\tau \hookrightarrow <\zeta (B+1), \cdots, \zeta (A+1)>_{\rho}\times \sigma. \eqno(6)
$$
On calcule $Jac_{\zeta (B+1), \cdots, \zeta (A+1)}\tau$ en utilisant cette inclusion. Le r\'esultat est r\'eduit \`a $\sigma$ si pour tout $x\in [\zeta (B+1), \zeta (A+1)]$ on a $Jac_{x, \cdots, \zeta (A+1)}\sigma=0$. Mais s'il n'en est pas ainsi, par r\'eciprocit\'e de Frobenius pour un bon choix de repr\'esentation irr\'eductible $\sigma'$, on a encore une inclusion
$$
\tau \hookrightarrow <\zeta (B+1), \cdots, \zeta (A+1)>_{\rho}\times_{c'\in [x,\zeta (A+1)]}\rho\vert\,\vert^{c'}\times \sigma'
$$
$$
\simeq \times_{c'\in [x,\zeta (A+1)]}\rho\vert\,\vert^{c'}\times <\zeta (B+1), \cdots, \zeta (A+1)>_{\rho}\times \sigma'.
$$
Et on aurait encore $Jac_{x, \cdots, \zeta (A+1)}\tau \neq 0$ ce qui a \'et\'e exclu. D'o\`u l'irr\'eductibilit\'e cherch\'ee qui conclut la r\'ecurrrence pour (ii).

\

Il reste donc \`a montrer les propri\'et\'es d'holomorphie pour $M_{i}(s)$ en $s=s_{0}$ pour $i=1,2$.

On commence par le cas o\`u $\zeta=-$; dans ce cas $<-(B+1), \cdots, -(A+1)>_{\rho}$ est une s\'erie discr\`ete tordue par le caract\`ere $\vert\, \vert^{-(B+A)/2-1}$. En particulier $M_{1}(s)$ est holomorphe en $s=s_{0}$ par positivit\'e: $(b_{0}-1)/2>-(B+A)/2-1$. Le m\^eme argument vaut pour $M_{2}(s)$: en effet $$-(b_{0}-1)/2 \geq -(a_{0}-1)/2-(b_{0}-1)/2=-A_{0}\geq -A\geq -(A+B)/2>-(A+B)/2-1.$$

On consid\`ere maintenant les cas o\`u $\zeta=+$; ici $<\zeta (B+1), \cdots, \zeta (A+1)>_{\rho}$ est une repr\'esentation de Speh. Donc l'induite
$$
St(\rho,a_{0})\vert\,\vert^{\pm s_{0}}\times <(B+1), \cdots,  (A+1)>_{\rho}
$$
est l'induite d'une s\'erie discr\`ete tordue avec une repr\'esentation de Speh; cela n'entre pas directement dans le cadre \'etudi\'e par Zelevinsky mais on a \'etudie cette situation dans \cite{mw} I.
Si $A_{0}<B+1$, $M_{1}(s)$ est holomorphe en $s=s_{0}$ pour $i=1,2$ simplement parceque les op\'erateur d'entrelacement standard  $\rho\vert\,\vert^{x}\times \rho\vert\,\vert^{y}\rightarrow \rho\vert\,\vert^{y}\times \rho\vert\,\vert^{x}$ (avec $x,y$ au voisinage de ${\mathbb R}$) n'ont de p\^oles que sur la droite $x-y=0$. On suppose donc que $A_{0}\geq (B+1)$. On consid\`ere l'op\'erateur normalis\'e \`a la Langlands-Shahidi comme en \cite{mw} (on calculera le facteur de normalisation ci-dessous) associ\'e \`a l'entrelacement $M_{1}(s)$; on le note $N_{1}(s)$. D'apr\`es \cite{mw} I.8, cet op\'erateur est holomorphe car $A_{0}\geq B+1$ par hypoth\`ese.
Le facteur de normalisation est $L(St(\rho,a_{0})\times \rho,s-(A+1))/L(St(\rho,a_{0})\times \rho, s-B)$. En posant encore $s=(b_{0}-1)/2+s'$ cela devient $L(\rho\times \rho,A_{0}-(A+1)+s')/L(\rho\times \rho, A_{0}-B)$. Avec nos hypoth\`eses ce facteur de normalisation n'a ni p\^ole ni z\'ero en $s=s_{0}$. 

Consid\'erons maintenant le cas de $M_{2}(s)$; ici on a directement l'holomorphie si $\zeta_{0}=-$ puisque $St(\rho,a_{0})\vert\,\vert^{-s_{0}}$ est alors la repr\'esentation $<-B_{0}, \cdots, -A_{0}>_{\rho}$. Supposons que  $\zeta_{0}=+$, ici $St(\rho,a_{0})\vert\,\vert^{-s_{0}}$ $=$ $<B_{0}, \cdots, -A_{0}>_{\rho}$. On a donc directement l'holomorphie cherch\'ee si $B_{0}<B+1$. Si $B_{0}\geq B+1$, on a l'holomorphie de l'op\'erateur d'entrelacement normalis\'e avec la m\^eme r\'ef\'erence. Le facteur de normalisation est ici
$$\frac{
L(St(\rho,a_{0})\times \rho,-s-(A+1))}{L(St(\rho,a_{0})\times \rho,-s-B)}=\frac{L(\rho\times \rho, B_{0}-(A+1)-s')}{L(\rho\times \rho, B_{0}-B-s')}.
$$
Or on a suppos\'e que l'on est dans le cas o\`u $B<B_{0}$ et $A_{0}\leq A$; de plus comme $\zeta_{0}=+$, on a aussi $a_{0}\geq b_{0}\geq 2$ d'o\`u aussi $B_{0}<A_{0}$, d'o\`u $B<B_{0}<A_{0}\leq A$. Le facteur de normalisation est donc encore holomorphe inversible en $s'=0$. Cela termine la preuve.
\section{G\'en\'eralisation}
Il y a plusieurs g\'en\'eralisations \`a faire; enlever l'hypoth\`ese de bonne parit\'e pour $\psi$, prendre aussi en compte le fait que la conjecture de Ramanujan n'est pas connue et remplacer la s\'erie discr\`ete $St(\rho,a_{0})$ consid\'er\'ee jusqu'\`a pr\'esent par une composante en la place p-adique fix\'ee d'une repr\'esentation cuspidale unitaire quelconque.

D'abord on consid\`ere un morphisme $\psi$ g\'en\'eral, c'est-\`a-dire contenant des sous-repr\'esenta\-tions irr\'eductibles de $W_{F}\times SL(2,{\mathbb C})\times SL(2,{\mathbb C})$ qui soit ne sont pas autoduales soit ne se factorisent pas par un groupe de m\^eme type que le groupe dual de $G$ et on va aussi ajouter des repr\'esentations non unitaires qui se doivent d'intervenir si la conjecture de Ramanujan n'est pas v\'erifi\'ee. On note $\psi_{bp}$ la somme des sous-repr\'esentations de $\psi$ qui sont autoduales et se factorisent par un sous-groupe de m\^eme type que le dual de $G$, $\psi_{mp}$ la somme des repr\'esentations irr\'eductibles unitaires qui ne se factorisent pas par un groupe de m\^eme type que $G$ et $\psi_{nu}$ la somme des repr\'esentations non unitaires. On peut pr\'eciser un peu $\psi_{nu}$. Pour tout $x\in {\mathbb R}$, on note $\psi_{x}$ la somme des sous-repr\'esentations irr\'eductibles incluses dans $\psi$ sur lesquelles $W_{F}$ agit par une repr\'esentation de d\'eterminant $\vert\, \vert^x$. Pour les paquets d'Arthur, on peut  (et on le fait) se limiter aux $\psi_{nu}$ tels que $\psi_{x}=0$ sauf \'eventuellement pour $\vert x\vert \in ]0,1/2[$. De plus $\psi$ \'etant autodual, on a n\'ecessairement $\psi_{x} \simeq (\psi_{-x})^*$. On note $\psi_{nu,>0}:=\sum_{x\in ]0,1/2[}\psi_{x}$. On d\'ecoupe aussi $\psi_{mp}$ en une somme de sous-repr\'esentations $\psi_{1/2,mp}\oplus \psi_{-1/2,mp}$ avec $\psi_{1/2,mp}\simeq \psi_{-1/2,mp}^*$ (cf. \ref{descriptionpaquets}). On a d\'ej\`a expliqu\'e comment les travaux de Bernstein-Zelevinsky et la preuve de la conjecture de Langlands permettent d'associer une repr\'esentation irr\'eductible d'un groupe $GL$ convenable \`a toute repr\'esentation irr\'eductible de $W_{F}\times SL(2,{\mathbb C})\times SL(2,{\mathbb C})$ en loc. cite. En induisant pour passer au cas des repr\'esentations semi-simples, on associe une repr\'esentation d'un $GL$ convenable \`a $\psi_{1/2,mp}$ et $\psi_{nu,>0}$ que l'on note $\pi^{GL}(\psi_{1/2,mp})$ et $\pi^{GL}(\psi_{nu,>0})$.

Le paquet de repr\'esentations associ\'e \`a $\psi$ est alors exactement l'ensemble des sous-quotients irr\'eductibles inclus dans les induite $\pi^{GL}(\psi_{nu,>0})\times \pi^{GL}(\psi_{1/2,mp})\times \pi_{bp}$ o\`u $\psi_{bp}$ parcourt l'ensemble des repr\'esentations dans le paquet associ\'e \`a $\psi_{bp}$. On a montr\'e en \cite{pourshahidi} 3.2 que les induites $\pi^{GL}(\psi_{1/2,mp})\times \pi_{bp}$ sont irr\'eductibles pour les groupes que l'on consid\`ere. Et on a repr\'esent\'e cette induite comme sous-module de Langlands convenable.
On g\'en\'eralisera ces r\'esultats \`a l'induction par $\pi^{GL}(\psi_{nu,>0})$. 

On d\'efinit encore $Jord(\psi)$ en utilisant une d\'ecomposition de $\psi$ en sous-repr\'esentations irr\'eductibles; pour ne pas modifier les notations d\'ej\`a utilis\'ees, on \'ecrit $Jord(\psi)$ sous forme de triplets $(\rho',a',b')_{x}$ o\`u $x\in {\mathbb R}$ et o\`u $\rho'$ est une repr\'esentation unitaire de $W_{F}$ et $a',b'$ sont des entiers repr\'esentant la dimension des repr\'esentations irr\'eductibles $sp_{a'}$, $sp_{b'}$ de $SL(2,{\mathbb C})$. La repr\'esentation correspondante est $\rho'\vert\, \vert^{x}\otimes sp_{a'}\otimes sp_{b'}$. Et on peut oublier $x$ si $x=0$.

Il faut aussi g\'en\'eraliser $St(\rho,a_{0})$ en une repr\'esentation irr\'eductible $\tau$ qui est unitaire, a un mod\`ele de Whittaker et s'\'ecrit donc comme induite $\times_{(\rho',a',x')\in {\mathcal J} }St(\rho',a')\vert\,\vert^{x'}$ et on impose aux $x'$ intervenant de v\'erifier $x'\in ]-1/2,1/2[$. Ce sont des propri\'et\'es qu'ont toutes les composantes locales des repr\'esentations cuspidales des groupes $GL$.

\subsection{Description des repr\'esentations dans les paquets d'Arthur g\'en\'eraux\label{descriptiongenerale}}
Soit $\psi$ un morphisme de $W_{F}\times SL(2,{\mathbb C})\times SL(2,{\mathbb C})$ dans le groupe dual de $G$; on le suppose semi-simple mais on ne le suppose plus tout \`a fait unitaire, on fait l'hypoth\`ese de l'introduction de cette section. On reprend les notations $\psi_{nu,>0},\psi_{nu,<0}, \psi_{mp}, \psi_{\pm 1/2,mp},\psi_{bp}$ que l'on a d\'ej\`a introduites.

En utilisant la correspondance de Langlands, on sait d\'efinir les repr\'esentations $\pi^{GL}(\psi)$, $\pi^{GL}(\psi_{nu,>0})$, $\cdots$ et l'on a
$$
\pi^{GL}(\psi)\simeq \pi^{GL}(\psi_{nu},>0)\times \pi^{GL}(\psi_{nu},<0)\times \pi^{GL}(\psi_{1/2,mp})\times \pi^{GL}(\psi_{-1/2,mp})\times \pi^{GL}(\psi_{bp}),\eqno(1)
$$
cette induite \'etant irr\'eductible.

Pour $(\rho,a,b)$ un triplet d\'efinissant une repr\'esentation irr\'eductible de $W_{F}\times SL(2,{\mathbb C})\times SL(2,{\mathbb C})$, on introduit encore une notation, pour tout $y$ tel que $y+(b-1)/2\in {\mathbb N}$, $J(St(\rho,a),-(b-1)/2,y)$ est l'unique sous-module irr\'eductible de l'induite, pour le $GL$ convenable:
$$
St(\rho,a)\vert\,\vert^{-(b-1)/2} \times \cdots \times St(\rho,a_{0})\vert\,\vert^{y};
$$
c'est le sous-module de Langlands de l'induite ci-dessus. Pour $b$ un entier, on pose $\delta_{b}=0$ si $b$ est impair et $1/2$ si $b$ est pair et $\delta'_{b}=1$ si $b$ est impair et $1/2$ si $b$ est pair. Dans tous les cas, on a $[-(b-1)/2,(b-1)/2]=[-(b-1)/2,-\delta_{b}]\cup [\delta'_{b},(b-1)/2]$ et c'est pour avoir ce d\'ecoupage que l'on a introduit ces notations compliqu\'ees.

On rappelle que l'on a montr\'e en \cite{pourshahidi} 3.2 que l'induite $\pi^{GL}(\psi_{1/2,mp})\time \pi_{bp}$ est irr\'eductible pour toute repr\'esentation irr\'eductible $\pi_{bp}$ dans le paquet de repr\'esentations associ\'e \`a $\psi_{bp}$. 
\begin{prop} Avec les notations ci-dessus, pour tout $\pi_{bp}$ dans le paquet de repr\'esentations associ\'ees \`a $\psi_{bp}$ l'induite $\pi^{GL}(\psi_{nu,>0})\times \pi^{GL}(\psi_{1/2,mp})\times \pi$ est irr\'eductible et c'est l'unique sous-module irr\'eductible de l'induite:
$$
\times_{(\rho,a,b)_{x}\in Jord(\psi_{nu,<0}),b>1}J(St(\rho,a),-(b-1)/2,-\delta'_{b})\vert\,\vert^{x}\times J(St(\rho^*,a),-(b-1)/2,-\delta_{b}))\vert\,\vert^{-x}$$
$$\times_{(\rho,a,b)_{x}\in Jord(\psi_{nu,<0}); b\equiv 1[2]}St(\rho,a)\vert\,\vert^{x}\times \pi^{GL}(\psi_{1/2,mp})\times \pi_{bp}
;
$$
\end{prop}
On commence par expliquer pourquoi on n'a pas mis d'ordre dans les inductions ci-dessus: l'induite  
$$
\times_{(\rho,a,b)_{x}\in Jord(\psi_{nu,<0}),b>1}J(St(\rho,a),-(b-1)/2,-\delta_{b})\vert\,\vert^{x}\times J(St(\rho,a,-(b-1)/2,-\delta'_{b}))\vert\,\vert^{-x}$$
pour le $GL$ convenable est irr\'eductible: en  \cite{mw} I.9 on a d\'emontr\'e que l'induite $$J(St(\rho,a),-(b-1)/2,y)\vert\,\vert^z\times J(St(\rho',a'),-(b'-1)/2,-y')\vert\,\vert^{z'}$$ est irr\'eductible si l'une des 3 conditions suivantes est r\'ealis\'ee:
$\rho\not\simeq \rho'$, $(a+b)/2+x-(a'+b')/2-x' \notin {\mathbb Z}$ ou $\vert y+z-y'-z'\vert <1$. Or on a, avec les notations introduites $-\delta_{b}+x\in ]-1,0[$ car $x\in ]-1/2,0[$ et $-\delta_{b}$ vaut $0$ ou $-1$ et $-\delta'_{b}-x\in ]-1,0[$ car $-x\in ]0,1/2[$ et $-\delta'_{b}=-1$ ou $-1/2$. Donc l'\'equivalent de $y+z-y'-z'$ est de la forme $t-t'$ avec $t,t'\in ]-1,0[$, d'o\`u $\vert t-t'\vert\in [0,1[$ et la troisi\`eme condition est toujours v\'erifi\'ee.

On note $\pi:=\pi^{GL}(\psi_{1/2,mp})\times \pi_{bp}$, on sait que c'est une repr\'esentation irr\'eductible comme on vient de le rappeler.  On montre d'abord que $\pi^{GL}(\psi_{nu,>0})\times \pi$ est irr\'eductible. On montre qu'une telle repr\'esentation a un unique sous-module irr\'eductible, que l'on note $\sigma$ et que $\sigma$ intervient avec multiplicit\'e 1 comme sous-quotient de l'induite: en effet, on note $d$ la dimension de la repr\'esentation $\psi_{nu,>0}$ et on calcule le module de Jacquet de l'induite $\pi^{GL}(\psi_{nu,>0})\times \pi$ le long du radical unipotent d'un parabolique standard de Levi $GL(d)\times G'$ o\`u $G'$ est un groupe classique de m\^eme type que $G$. Dans ce module de Jacquet, on ne conserve que les sous-repr\'esentations irr\'eductibles dont le support cuspidal pour l'action de $GL(d)$ est le m\^eme que celui de la repr\'esentation $\pi^{GL}(\psi_{nu,>0})$. Il  ne reste plus que la repr\'esentation irr\'eductible $\pi^{GL}(\psi_{nu,>0})\otimes \pi$. Soit $\sigma$ une sous-repr\'esentation irr\'eductible de l'induite. Par r\'eciprocit\'e de Frobenius le module de Jacquet de $\sigma$ doit contenir cette repr\'esentation irr\'eductible et comme cette repr\'esentation n'intervient qu'avec multiplicit\'e 1 dans le module de Jacquet de toute l'induite, $\sigma$ est unique et a multiplicit\'e 1 comme sous-quotient irr\'eductible. On d\'emontre de la m\^eme fa\c{c}on que l'induite a un unique quotient irr\'eductible et qu'un tel quotient intervient avec multiplicit\'e 1 comme sous-quotient irr\'eductible de l'induite: ici il suffit de dualiser pour se ramener \`a une assertion o\`u quotient est remplac\'e par sous-module. En dualisant on change $\psi_{nu,>0}$ en $\psi_{nu,<0}$ et  on change \'eventuellement   $\pi$ mais par un automorphisme que l'on contr\^ole (cf. \cite{mvw}, II.1) et la m\'ethode est la m\^eme. On note $\tau$ ce quotient irr\'eductible. On remarque pour la suite que le module de Jacquet de $\tau$ contient n\'ecessairement $\pi^{GL}(\psi_{nu,<0})\otimes \pi$.

On montre que l'induite $\pi^{GL}(\psi_{nu,>0})\times \pi$ est isomorphe \`a l'induite $\pi^{GL}(\psi_{nu,<0})\times \pi$ en construisant l'isomorphisme avec des op\'erateurs d'entrelacement: en \'ecrivant $\pi^{GL}(\psi_{nu,>0})$ sous la forme $\pi_{1}\times \pi_{\ell}$, chaque $\pi_{i}$ \'etant une repr\'esentation de la forme $Speh(St(\rho',a'),b')\vert\,\vert^{x'}$ et $\pi^{GL}(\psi_{nu,<0})$ est \'ecrit $\pi_{\ell}^*\times \cdots \times \pi_{1}^*$ et l'op\'erateur d'entrelacement est celui qui correspond \`a l'\'el\'ement du groupe de Weyl de longueur minimal qui conjugue les repr\'esentations de la fa\c{c}on sugg\'er\'ee par les notations. On va v\'erifier que pour tout $i\in [1,\ell]$ les op\'erateurs d'entrelacement standard $\pi_{i}\times \pi\rightarrow \pi_{i}^*\times \pi$ sont holomorphes bijectifs. Avec les r\'ef\'erences d\'ej\`a donn\'ee, on sait que pour $i,j\in [1,\ell]$ les induites $\pi_{j}\times \pi_{i}^*$ sont irr\'eductibles et on aura ainsi l'isomorphisme cherch\'e. Pour v\'erifier l'affirmation pr\'ec\'edente, on \'ecrit $\pi_{i}$ \`a l'aide des segments de Zelevinsky:
$$
\pi_{i}=\bigg< \begin{matrix} (a'-b')/2+x' & \cdots & (a'+b')/2-1+x'\\
\vdots&\vdots&\vdots&\\
-(a'+b')/2+1+x' &\cdots &-(a'-b')/2+x'
\end{matrix}\bigg>_{\rho'}
$$
o\`u les lignes sont des segments croissants et les colonnes des segments d\'ecroissants (cf. \ref{correspondance}). On note ${\mathcal E}$ l'ensemble des coefficients de la matrice \'ecrite, ensemble que l'on ordonne en lisant d'abord de gauche \`a droite puis de haut en bas. Soit $z\in {\mathcal E}$; on note ${\mathcal E}_{<z}$ les \'el\'ements de ${\mathcal E}$ \`a la gauche de $z$ (pour l'ordre) et ${\mathcal E}_{>z}$ ceux \`a la droite de $z$; aucun de ces ensembles ne contient $z$. Dans la notation ci-dessous $opp$ veut dire que l'on inverse l'ordre et l'on regarde l'op\'erateur
$$
\times_{z'\in {\mathcal E}_{<z}}\rho'\vert\,\vert^{z'}\times \rho'\vert\,\vert^{z}\times_{z''\in {\mathcal E}^{opp}_{>z}}\rho^{*' }\vert\,\vert^{-z''}\times \pi \rightarrow \times_{z'\in {\mathcal E}_{<z}}\rho'\times_{z''\in {\mathcal E}^{opp}_{>z}}\rho^{*' }\vert\,\vert^{-z''}\times \rho'\vert\,\vert^{z}\times \pi;
$$
il est holomorphe et inversible; c'est quasi automatique si $\rho\not\simeq \rho^*$ mais l'argument g\'en\'eral, qui vaut dans tous les cas, est que pour tout $z''\in {\mathcal E}_{>z}$ $z+z''\notin 1/2 {\mathbb Z}$. Ensuite on v\'erifie que l'op\'erateur $\rho'\vert\,\vert^{z}\times \pi \rightarrow \rho^{*'}\vert\,\vert^{-z}\times \pi$ est holomorphe inversible et ici c'est le fait que $z$ n'est pas un demi-entier qui sert (cf. la remarque ci-dessous). Finalement on montre que l'op\'erateur
$$
\times_{z'\in {\mathcal E}_{<z}}\rho'\vert\,\vert^{z'}\times \rho'\vert\,\vert^{z}\times_{z''\in {\mathcal E}^{opp}_{>z}}\rho^{*' }\vert\,\vert^{-z''}\times \pi \rightarrow \times_{z'\in {\mathcal E}_{<z}}\rho'\times_{z''\in {\mathcal E}^{opp}_{>z}}\rho^{*' }\vert\,\vert^{-z''}\times \rho^{'*}\vert\,\vert^{-z}\times \pi
$$
est holomorphe inversible. Ainsi en restreignant cet op\'erateur \`a $\pi_{i}\times \pi$ on obtient l'assertion cherch\'ee. Ensuite on conclut facilement \`a l'isomorphisme cherch\'e. Ainsi $\sigma$ est aussi un sous-module irr\'eductible de $\pi^{GL}(\psi_{\nu,<0})\times \pi$ et son module de Jacquet contient donc aussi $\pi^{GL}(\psi_{\nu,<0})\otimes \pi$. Comme ce module de Jacquet \'etait n\'ecessairement dans celui de $\tau$ et comme il n'intervient qu'avec multiplicit\'e 1 dans le module de Jacquet de toute l'induite, n\'ecessairement $\tau\simeq \sigma$. D'o\`u l'irr\'eductibilit\'e cherch\'ee.

Montrons maintenant la 2e partie de la proposition. Soit $(\rho,a,b)_{x}\in Jord(\psi_{nu,<0})$. On vient sait par d\'efinition et avec les notations introduites avant l'\'enonc\'e que $$Speh(St(\rho,a),b)\vert\,\vert^{x}\hookrightarrow J(St(\rho,a),-(b-1)/2-\delta_{b})\times J(St(\rho,a),\delta'_{b},(b-1)/2).
$$
Comme ci-dessus, on v\'erifie que $$J(St(\rho,a),\delta'_{b},(b-1)/2)\vert\,\vert^{x}\times \pi\simeq
J(St(\rho^*,a),-(b-1)/2,-\delta'_{b})\vert\,\vert^{-x}\times \pi.
$$
Soit $(\rho',a',b')_{x'}\in Jord(\psi_{nu,<0})$, par d\'efinition $Speh(St(\rho',a'),b')\vert\,\vert^{x'}=J(St(\rho',a'),-(b'-1)/2,(b'-1)/2)\vert\,\vert^{x'}$. On a introduit une condition  de segments li\'es de \cite{mw} 1.7 qui g\'en\'eralise celle de Zelevinsky pour ce genre de repr\'esentation et qui assure l'irr\'eductibilit\'e de l'induite $Speh(St(\rho',a'),b')\vert\,\vert^{x'}\times J(St(\rho,a),-(b-1)/2,y)\vert\,\vert^{\tilde{x}}$. Supposons que $b'\geq b$, on a alors $-(b'-1)/2\leq -(b-1)$ et $y\leq (b'-1)/2$ et encore $y+\tilde{x}< (b'-1)/2+x'+\vert (a-a')/2\vert+1$ tandis que $-(b'-1)/2+x' <-(b-1)/2+\tilde{x}+\vert (a-a')/2\vert+1$ et les segments ne sont donc pas li\'es. Donc ainsi
$$
Speh(St(\rho,a'),b')\vert\,\vert^{x'}\times Speh(St(\rho,a),b)\vert\,\vert^{x}\times \pi\hookrightarrow$$
$$
Speh(St(\rho,a'),b')\vert\,\vert^{x'}\times J(St(\rho,a),-(b-1)/2,-\delta_{b})\vert\,\vert^{x}\times J(St(\rho^*,a),-(b-1)/2,-\delta'_{b})\vert\,\vert^{-x}\times \pi
$$
$$
\simeq 
J(St(\rho,a),-(b-1)/2,-\delta_{b})\vert\,\vert^{x}\times J(St(\rho^*,a),-(b-1)/2,-\delta'_{b})\vert\,\vert^{-x}\times
Speh(St(\rho,a'),b')\vert\,\vert^{x'}\times \pi
$$
$$
\hookrightarrow 
J(St(\rho,a),-(b-1)/2,-\delta_{b})\vert\,\vert^{x}\times J(St(\rho^*,a),-(b-1)/2,-\delta'_{b})\vert\,\vert^{-x}$$
$$\times 
J(St(\rho',a'),-(b'-1)/2,-\delta'_{b'})\vert\,\vert^{x'}\times J(St(\rho^{'*},a'),-(b'-1)/2,-\delta'_{b'})\vert\,\vert^{-x'}\times \pi.
$$
On peut ensuite continuer avec plus de 2 facteurs et obtenir l'inclusion de l'\'enonc\'e. Le fait que l'induite de droite a un unique sous-module irr\'eductible vient du fait que l'on peut remplacer l'induite des repr\'esentations $J(\cdots)$ par une induite de repr\'esentation de Steinberg tordu, les exposants \'etant dans la chambre de Weyl ferm\'e n\'egative; c'est  l'argument de  \cite{manuscripta} 3.6, lemme 1 qui repose encore une fois sur \cite{mw} I.9. Cela termine la preuve.

\begin{rem} {\bf 1} Avec les notations pr\'ec\'edentes, pour $\rho'$ une repr\'esentation cuspidale unitaire d'un groupe $GL$ et pour $x\in {\mathbb R}-1/2 {\mathbb Z}$, l'op\'erateur d'entrelacement standard $$M(\rho,x):\rho'\vert\,\vert^{x}\times \pi^{GL}(\psi_{1/2,mp})\times \pi_{bp} \rightarrow \rho^{'*}\vert\,\vert^{-x}\times \pi^{GL}(\psi_{1/2,mp})\times \pi_{bp}$$ est holomorphe et le compos\'e $M(\rho',-x)\circ M(\rho,x)$ est un scalaire non nul. En particulier ces induites sont irr\'eductibles.
\end{rem}
Le support cuspidal de $\times \pi^{GL}(\psi_{1/2,mp})\times \pi_{bp}$ est de la forme une collection de repr\'esentations cuspidales de $GL$, de la forme  $\rho\vert\,\vert^{c}$ avec $\rho$ une repr\'esentation cuspidale unitaire d'un $GL$ et $c$ un \'el\'ement de $1/2{\mathbb Z}$ et d'une repr\'esentation cuspidale $\pi_{cusp}$ d'un groupe classique de m\^eme type que $G$. On peut \'ecrire
$$
 \pi^{GL}(\psi_{1/2,mp})\times \pi_{bp}\hookrightarrow \times_{(\rho,c)\in {\mathcal E}}\rho\vert\,\vert^{c}\times \pi_{0}, \eqno(1)
$$
o\`u ${\mathcal E}$ est un ensemble totalement ordonn\'e convenable. On sait que pour tout $(\rho,c)\in {\mathcal E}$ l'op\'erateur d'entrelacement standard qui \'echange les 2 facteurs de l'induite $\rho'\vert\,\vert^{x}\times \rho \vert\,\vert^{c}$ est holomorphe pour $x-c\neq 0$ ainsi que  l'op\'erateur d'entrelacement standard $$\rho'\vert\,\vert^{x}\times \pi_{0}\rightarrow \rho^{'*}\vert\,\vert^{-x}\times \pi_{0}$$ en $x\neq 0$ et l'op\'erateur qui \'echange les 2 copies de $\rho\vert\,\vert^{c}\times \rho^{'*}\vert\,\vert^{-x}$ en $x+c\neq 0$ par les r\'esultats g\'en\'eraux d'Harish-Chandra. D'o\`u de fa\c{c}on \'evidente l'holomorphie de l'op\'erateur d'entrelacement standard d\'ecrit dans l'\'enonc\'e. On sait aussi que $M(\rho^{'*},-x)\circ M(\rho,x)$ est un scalaire.

Pour le calculer on peut encore utiliser l'inclusion (1) et d\'ecomposer en op\'erateur \'el\'emen\-taire; on v\'erifie que chaque induite $\rho'\vert\,\vert^{x}\times \rho \vert\,\vert^{c}$, $\rho'\vert\,\vert^{x}\times \pi_{0}$, $\rho\vert\,\vert^{c}\times \rho^{'*}\vert\,\vert^{-x}$ est irr\'eductible car $x\notin 1/2{\mathbb Z}$. Le compos\'e des op\'erateurs d'entrelacement d\'ecrit est donc une bijection et le scalaire ne peut \^etre 0. Cela entra\^{\i}ne que l'op\'erateur $M(\rho',x)$ est bijectif, d'inverse $M(\rho^{'*},-x)$. Par un calcul de module de Jacquet, on v\'erifie que l'induite $\rho'\vert\,\vert^{x}\times \pi^{GL}(\psi_{1/2,mp})\times \pi_{pb}$ a un unique sous-module irr\'eductible: le terme $\rho'\vert\,\vert^{x}\otimes \biggl(\pi^{GL}(\psi_{1/2,mp})\times \pi_{pb}\biggr)$ arrive avec multiplicit\'e 1 dans le module de Jacquet de l'induite et que ce sous-module irr\'eductible a multiplicit\'e 1 comme sous-quotient irr\'eductible. En dualisant, on v\'erifie que ce sous-module irr\'eductible est aussi l'unique quotient irr\'eductible de l'induite $\rho^{'*}\vert\,\vert^{-x}\times \pi^{GL}(\psi_{1/2,mp})\times \pi_{pb}$. Par la bijectivit\'e de $M(\rho',x)$, cela entra\^{\i}ne l'irr\'eductibilit\'e cherch\'ee et la remarque.
\begin{rem} {\bf 2} Avec les notations pr\'ec\'edentes, l'induite 
$$
\times_{(\rho,a,b)\in Jord(\psi_{1/2,mp})} J(St(\rho,a),-(b-1)/2,-\delta_{b})\times J(St(\rho,a),-(b-1)/2,-\delta'_{b})\times \pi_{bp}
$$
a un unique sous-module irr\'eductible et ce sous-module est l'induite $\pi^{GL}(\psi_{1/2,mp})\times \pi_{bp}$.
\end{rem}
En \cite{pourshahidi} 3.2, on a \'etabli l'irr\'eductibilit\'e de $\pi^{GL}(\psi_{1/2,mp})\times \pi_{bp}$. Comme le module de Jacquet de cette induite contient le terme
$$
\biggl(\times_{(\rho,a,b)\in Jord(\psi_{1/2,mp})} J(St(\rho,a),-(b-1)/2,-\delta_{b})\times J(St(\rho,a),-(b-1)/2,-\delta'_{b}))\otimes \pi_{bp}
$$
la repr\'esentation $\pi^{GL}(\psi_{1/2,mp})\times \pi_{bp}$ est un sous-module irr\'eductible de l'induite \'ecrite. Montrons que cette induite a un unique sous-module irr\'eductible; avec \cite{mw} I.8, on peut \'echanger les facteurs comme on veut et on \'ecrit cette induite sous la forme
$$
\times_{(\rho,a,b)\in Jord(\psi_{1/2,mp}); b\equiv 0[2]}J(St(\rho,a),-(b-1)/2,-1/2)\times J(St(\rho,a),-(b-1)/2,-1/2))$$
$$
\times_{(\rho,a,b)\in Jord(\psi_{1/2,mp}),b\equiv 1[2]}J(St(\rho,a),-(b-1)/2,-1)$$
$$
\times_{(\rho,a,b)\in Jord(\psi_{1/2,mp}); b\equiv 0[2]} J(St(\rho,a),-(b-1)/2,0)\times \pi_{bp}.
$$
Or $J(St(\rho,a),-(b-1)/2,0)\hookrightarrow J(St(\rho,a),-(b-1)/2,-1)\times St(\rho,a)$; en le faisant progressivement, on peut encore remplacer la ligne ci-dessus par une inclusion dans l'induite
$$
\times_{(\rho,a,b)\in Jord(\psi_{1/2,mp}),b\equiv 1[2]}J(St(\rho,a),-(b-1)/2,-1)\times_{(\rho,a,b)\in Jord(\psi_{1/2,mp}); b\equiv 1[2]}St(\rho,a)\times \pi_{bp}.
$$
Cela inclut l'induite de l'\'enonc\'e dans l'induite
$$
\times_{(\rho,a,b)\in Jord(\psi_{1/2,mp})}J(St(\rho,a),-(b-1)/2,-\delta'_{b})
\times_{(\rho,a,b)\in Jord(\psi_{1/2,mp}),b\equiv 1[2]}St(\rho,a) \times \pi_{bp}.
$$
L'induite $\times_{(\rho,a,b)\in Jord(\psi_{1/2,mp}),b\equiv 1[2]}St(\rho,a) \times \pi_{bp}$ est irr\'eductible par \cite{pourshahidi} 3.2 et toute l'induite \'ecrite ci-dessus contient dans son module de Jacquet le terme 
$$
\biggl(\times_{(\rho,a,b)\in Jord(\psi_{1/2,mp})}J(St(\rho,a),-(b-1)/2,-\delta'_{b})\biggr)\otimes
\biggl(\times_{(\rho,a,b)\in Jord(\psi_{1/2,mp}),b\equiv 1[2]}St(\rho,a) \times \pi_{bp}\biggr)
$$
avec multiplicit\'e 1. Cela prouve l'unicit\'e du sous-module irr\'eductible cherch\'e; ce sous-module irr\'eductible est le sous-quotient de Langlands (bien d\'efini, quand on remplace $\pi_{bp}$ par l'induite avec les param\`etres de Langlands de $\pi_{bp}$).

\subsection{Propri\'et\'es d'holomorphie des op\'erateurs d'entrelacement normalis\'es\label{generalisation}}
On fixe $\psi$ un morphisme de $W_{F}\times SL(2,{\mathbb C})\times SL(2,{\mathbb C})$ comme dans \ref{descriptiongenerale}, c'est-\`a-dire que $\psi= \psi_{nu}\oplus \psi_{mp}\oplus \psi_{bp}$. On fixe aussi une repr\'esentation $\tau$ qui joue le r\^ole de la composante locale d'une forme automorphe cuspidale, autoduale  d'un groupe $GL$. On suppose donc  que $\tau$ est unitaire et a un mod\`ele de Whittaker avec la propri\'et\'e de l'introduction de cette section.

Pour $s\in {\mathbb C}$, on a d\'efini la fonction m\'eromorphe, qui utilise des facteurs $L$ pour des groupes lin\'eaires  $$r(s,\tau,\psi):=\frac{ L(\tau\times \pi^{GL}(\psi),s)}{L(\tau\times \pi^{GL}(\psi),s+1)
}\frac{L(\tau,r_{G},2s)}{L(\tau,r_{G},2s+1)}.$$

Soit $\pi$ dans le paquet de repr\'esentations associ\'e \`a $\psi$.
On consid\`ere l'op\'erateur d'entrelacement standard, pour tout $s\in {\mathbb C}$:
$$
M(s,\tau,\pi): \quad \tau\vert\,\vert^{s}\times \pi\rightarrow \tau^{*}\vert\,\vert^{-s}\times \pi;
$$
il est associ\'e \`a l'\'el\'ement du groupe de Weyl de $G$, de longueur minimale dans sa double classe modulo le groupe de Weyl du parabolique maximal. Et on pose $N_{\psi}(s,\tau,\pi):=r(s,\tau,\psi)^{-1}\times M(s,\tau,\pi)$. 

\begin{prop}L'op\'erateur $N_{\psi}(s,\tau,\pi)$ est holomorphe en tout $s\in {\mathbb R}_{> 0}$.
\end{prop}

On fixe $\pi_{bp}$ dans le paquet associ\'e \`a $\psi_{bp}$ avec les notations de \ref{descriptiongenerale} de tel sorte que $\pi\simeq \pi^{GL}(\psi_{nu,>0})\times \pi^{GL}(\psi_{1/2,mp})\times \pi_{bp}$. On d\'ecompose l'op\'erateur d'entrelacement en le produit des op\'erateurs d'entrelacement suivant:
$$
\tau\vert\,\vert^{s}\times \pi^{GL}(\psi_{nu,>0})\times \pi^{GL}(\psi_{1/2,mp})\rightarrow
\pi^{GL}(\psi_{nu,>0})\times \pi^{GL}(\psi_{1/2,mp}) \times \tau\vert\,\vert^{s};\eqno(1)
$$
$$
\tau\vert\,\vert^{s}\times \pi_{bp}\rightarrow \tau^*\vert\,\vert^{-s}\times \pi_{bp};\eqno(2)
$$
$$
\pi^{GL}(\psi_{nu,>0})\times \pi^{GL}(\psi_{1/2,mp})\times \tau^*\vert\,\vert^{-s}\rightarrow
\tau^*\vert\,\vert^{-s}\times \pi^{GL}(\psi_{nu,>0})\times \pi^{GL}(\psi_{1/2,mp}); \eqno(3)
$$
o\`u le 1e et le 3e op\'erateur sont situ\'es dans des groupes $GL$ convenables. On pose $\sigma:=\pi^{GL}(\psi_{nu,>0})\times \pi^{GL}(\psi_{1/2,mp})$ et on sait que $\sigma^*\simeq \pi^{GL}(\psi_{nu,<0})\times \pi^{GL}(\psi_{-1/2,mp})$ puisque $\psi$ est \`a valeurs dans le groupe dual de $G$.

On utilise la d\'ecomposition de \ref{descriptiongenerale} (1) pour d\'ecomposer $L(\tau\times \pi^{GL}(\psi),s)/L(\tau\times \pi^{GL}(\psi),s+1)$ en le produit des 3 facteurs ci-dessous. Pour simplifier l'\'ecriture on pose $\sigma:=\pi^{GL}(\psi_{nu,>0})\times \pi^{GL}(\psi_{1/2,mp})$ et on utilise le fait que $\sigma^*\simeq \pi^{GL}(\psi_{nu,<0})\times \pi^{GL}(\psi_{-1/2,mp})$. On a donc:
$
L(\tau\times \pi^{GL}(\psi),s)/L(\tau\times \pi^{GL}(\psi),s+1)=$
$$
L(\tau\times \sigma^*,s)/L(\tau\times \sigma^*,s+1)\eqno(4)
$$
$$
\times L(\tau\times \pi^{GL}(\psi_{bp}),s)/L(\tau,\times \pi^{GL}(\psi_{bp}),s+1)\eqno(5)
$$
$$
\times L(\tau\times \sigma,s)/L(\tau\times \sigma,s+1).\eqno(6)
$$
L'op\'erateur d'entrelacement (1) normalis\'e par  la fonction (4), est exactement l'op\'erateur d'entrelacement \'etudi\'e dans \cite{mw} I.2 et il en est de m\^eme de (3) normalis\'e par (6). On v\'erifie que l'hypoth\`ese de \cite{mw} I.8 est satisfaite de fa\c{c}on \`a tirer de cette r\'ef\'erence l'holomorphie: on se place en $s\in {\mathbb R}_{>0}$ et l'hypoth\`ese en question est que les modules de Speh qui d\'efinissent la repr\'esentation $\tau\vert\,\vert^s$ (resp. $\sigma$) dominent ceux qui d\'efinissent la repr\'esentation $\sigma$ (resp. $\tau^*\vert\,\vert^{-s}$) pour l'op\'erateur (1) (resp. (3)) ou ne sont pas li\'es. Faisons le pour l'op\'erateur (1): la repr\'esentation $\tau$ est une induite de $St(\rho'',c)\vert\,\vert^{x+s}$ avec $x\in ]-1/2,1/2[$ donc de la forme $Speh(St(\rho'',c),1)\vert\,\vert^{x+s}$. Ceux qui d\'efinissent $\sigma$ sont de la forme $Speh (St(\rho''',a),b)\vert\,\vert^{z}$ avec $z\in ]-1/2,1/2[$; donc on a $x+s-z>-1$ et si les segments sont li\'es, ce nombre est $\geq 0$ ainsi le premier domine le second (cf. loc.cite). Pour l'op\'erateur (3), on a $z+s-x$ qui a la m\^eme propri\'et\'e. D'o\`u l'holomorphie des op\'erateurs (1) et (3) normalis\'es.

Il reste l'op\'erateur (2) normalis\'e par la fonction que l'on est en droit d'\'ecrire $r(s, \tau,\psi_{bp})$. Ici on \'ecrit $\tau\simeq \times_{(\rho',a',x')\in {\mathcal J}}St(\rho',a')\vert\,\vert^{x'}$ o\`u ${\mathcal J}$ est un ensemble, totalement ordonn\'e, de triplets form\'es d'une repr\'esentation cuspidale unitaire, d'un entier et d'un nombre r\'eel. La fonction $r(s,\tau,\psi_{bp})$ se d\'ecompose de fa\c{c}on analogue en produit. L'op\'erateur d'entrelacement s'\'ecrit ainsi, apr\`es avoir fix\'e l'isomorphisme de $\tau$ avec l'induite \'ecrite:
$$
\times_{(\rho',a',x')\in {\mathcal J}}St(\rho',a')\vert\,\vert^{x'+s} \times \pi_{bp} \rightarrow \times_{(\rho',a',x')\in {\mathcal J}^{opp}}St(\rho^{'*},a)\vert\,\vert^{-x'-s}\times \pi_{bp},
$$
o\`u l'exposant $opp$ signifie que l'ordre sur ${\mathcal J}$ a \'et\'e invers\'e. Donc il faut donc d\'emontrer l'assertion quand ${\mathcal J}$ est r\'eduit \`a un \'el\'ement c'est-\`a-dire quand $\tau\simeq St(\rho',a')\vert\,\vert^{x'}$ avec $x'\in ]-1/2,1/2[$. On a d\'emontr\'e l'assertion d'holomorphie en \cite{holomorphie} dans le cas o\`u $s+x'\geq 0$ et $\rho'$ est autodual. Si  si $s+x'\notin 1/2{\mathbb Z}$, on a prouv\'e l'holomorphie dans la remarque de \ref{descriptiongenerale}; la preuve donn\'ee en loc.cite s'applique aussi au cas o\`u $\rho'$ n'est pas autoduale  puisqu'ici la repr\'esentation du groupe classique est $\pi_{bp}$, associ\'e \`a des triplets de bonne parit\'e donc avec une repr\'esentation cuspidale qui, elle, est autoduale. On couvre ainsi tous les cas puisque pour $s\geq 0$, $s+x'>-1/2$ par hypoth\`ese sur $x'$. Cela termine la preuve de l'holomorphie.
\subsection{Image des op\'erateurs d'entrelacement normalis\'es\label{imagecageneral}}
On reprend les notations $\tau,\psi,\pi$ de \ref{generalisation}. On a d\'efini en loc.cite l'op\'erateur d'entrelacement $N_{\psi}(s,\tau,\pi)$ et montr\'e son holomorphie pour tout $s\in {\mathbb R}_{>0}$.

On va d\'ecrire l'image; dans le cas g\'en\'eral le r\'esultat est purement qualitatif, on veut montrer que l'image est nulle ou irr\'eductible; on le fera sous l'hypoth\`ese que  $s\geq 1/2$, vue nos hypoth\`eses cette condition est n\'ecessaire pour avoir l'irr\'eductibilit\'e: en effet, fixons $\rho$ tel que l'induite $\rho\times \pi$ est r\'eductible (il y a de nombreux exemples) et posons $\tau=\rho\vert\,\vert^{-x}\times \rho\vert\,\vert^{x}$ avec $x\in ]0,1/2[$. Alors en $s=x$ l'op\'erateur d'entrelacement n'a aucun raison d'avoir une image irr\'eductible.

\subsubsection{Description qualitative de l'image\label{qualitatifgeneral}}
\begin{prop}L'image de $N_{\psi}(s,\tau,\pi)$ en $s=s_{0}\in {\mathbb R}_{\geq 1/2}$ est soit identiquement 0, soit est  une repr\'esentation irr\'eductible.
\end{prop}
On \'ecrit $\tau\vert\,\vert^{s_{0}}$ comme induite irr\'eductible et pour fixer les notations les plus simples possibles, on \'ecrit:
$$
\tau\vert\,\vert^{s_{0}}\simeq \times_{(\rho',a',s')\in {\mathcal J}}St(\rho',a')\vert\,\vert^{s'};
$$
ici on sait que $s'>0$. On ordonne ${\mathcal J}$ de telle sorte que les premiers facteurs sont ceux pour lesquels $s'\notin 1/2{\mathbb Z}$ et les derniers ont la propri\'et\'e oppos\'ee. D'o\`u en regroupant ces 2 types de facteurs, $\tau\vert\,\vert^{s_{0}}\simeq \tau_{ne}\times \tau_{e}$, $ne$ pour non demi-entier et $e$ pour demi-entier. On \'ecrit $\pi=\pi^{GL}(\psi_{nu,>0})\times \pi^{GL}(\psi_{1/2,mp})\times \pi_{bp}$ avec les notations de \ref{descriptiongenerale}. On note $\pi_{0}$ l'image de l'op\'erateur d'entrelacement normalis\'e:
$$
\tau_{e}\times \pi^{GL}(\psi_{1/2,mp})\times \pi_{bp}\rightarrow \tau_{e}^*\times \pi^{GL}(\psi_{1/2,mp})\times \pi_{bp}.
$$
Comme le support cuspidale de $\pi^{GL}(\psi_{nu,<0})$ et celui de $\pi_{e}$ n'ont pas les m\^emes propri\'et\'es d'int\'egralit\'e les op\'erateurs d'entrelacement qui \'echange ces facteurs sont des isomorphismes qu'on les normalise ou pas. Cela entra\^{\i}ne que l'image cherch\'ee n'est autre que l'image de l'op\'erateur d'entrelacement normalis\'e par $r(s, \tau_{ne},\psi)$ de l'op\'erateur d'entrelacement:
$$
\tau_{ne}\times \pi^{GL}(\psi_{nu,>0})\times \pi_{0}\rightarrow \tau_{ne}^*\times \pi^{GL}(\psi_{nu,>0})\times \pi_{0}.
$$En particulier si $\pi_{0}=0$ l'image cherch\'ee est nulle; supposons donc que $\pi_{0}\neq 0$; on admet pour le moment que $\pi_{0}$ est soit nul soit  irr\'eductible. On conclut avec l'hypoth\`ese que $\pi_{0}$ est irr\'eductible.

L'argument, dont on donne les d\'etails ci-dessous,   est que le terme de gauche a un unique quotient irr\'eductible et que celui de droite a un unique sous-module irr\'eductible que ces 2 repr\'esentations sont isomorphes et interviennent avec multiplicit\'e 1 comme sous-quotient irr\'eductible des induites \'ecrites. Cela force le fait que tout morphisme non nul entre les 2 induites est d'image irr\'eductible. Et prouve l'assertion. La d\'emonstration sera r\'eutilis\'ee en \ref{parametregeneralisation} car elle donne une bonne description de l'image une fois la repr\'esentation $\pi_{0}$ connue.

En g\'en\'eralisant \ref{descriptiongenerale}, on montre que $\pi^{GL}(\psi_{nu,>0})\times \pi_{0}$ est l'unique sous-module irr\'eductible de l'induite (avec les notations de loc. cite)
$$
\times_{(\rho,a,b)_{x}\in Jord(\psi_{nu,<0})}J(St(\rho,a),-(b-1)/2,-\delta_{b})\vert\,\vert^{x}\times J(St(\rho^*,a),-(b-1)/2,-\delta'_{b})\vert\,\vert^{-x}\times \pi_{0}.
$$
On v\'erifie  que $\tau^*_{ne}\times_{(\rho,a,b)_{x}\in Jord(\psi_{nu,<0})}J(St(\rho,a),-(b-1)/2,-\delta_{b})\vert\,\vert^{x}\times J(St(\rho,a),-(b-1)/2,-\delta'_{b})\vert\,\vert^{-x}$ est le sous-module de Langlands de l'induite (cf. \cite{manuscripta} 3.6 lemme 1, o\`u l'argument a \'et\'e d\'etaill\'e). Avec cela il est simple de montrer l'unicit\'e du sous-module irr\'eductible de l'induite comme annonc\'e, ainsi que sa multiplicit\'e 1 comme sous-quotient irr\'eductible; d'apr\`es cette description, on a les param\`etres de Langlands de ce sous-module en fonction de ceux de $\pi_{0}$. 

Avec les m\^emes arguments, ou plus simplement en dualisant on montre que $\pi$ est l'unique quotient irr\'eductible de l'induite
$$
\times_{(\rho,a,b)_{x}\in Jord(\psi_{nu,>0})}J(St(\rho,a),\delta_{b},(b-1)/2)\vert\,\vert^{x}\times J(St(\rho^*,a),\delta'_{b},(b-1)/2)\vert\,\vert^{-x}\times \pi_{0}.
$$
Et que $\tau_{ne}\times \pi$ est le quotient de Langlands de l'induite $\tau_{ne}\times$ l'induite que l'on vient d'\'ecrire. D'o\`u l'unicit\'e du quotient irr\'eductible, c'est le quotient de Langlands et il a les m\^emes param\`etres de Langlands que le sous-module d\'ecrit pr\'ec\'edemment. Cela suffit pour la proposition.

Pour la suite, on veut v\'erifier que l'op\'erateur est non nul si $\pi_{0}\neq 0$: il suffit de montrer que l'op\'erateur d'entrelacement (les $(\rho,a,b)_{x}$ parcourt $Jord(\psi_{nu,<0})$ ce que l'on n'a pas la place d'\'ecrire ci-dessous)
$$
\tau_{ne}\vert\,\vert^{s'} \times_{(\rho,a,b)_{x}}J(St(\rho,a),-(b-1)/2,-\delta_{b})\vert\,\vert^{x}\times J(St(\rho^*,a),-(b-1)/2,-\delta'_{b})\vert\,\vert^{-x}\times \pi_{0}$$
$$
\rightarrow
$$
$$
\tau_{ne}^*\vert\,\vert^{-s'}\times_{(\rho,a,b)_{x}}J(St(\rho,a),-(b-1)/2,-\delta_{b})\vert\,\vert^{x}\times J(St(\rho^*,a),-(b-1)/2,-\delta'_{b})\vert\,\vert^{-x}\times \pi_{0}
$$
normalis\'e par la fonction $r(\psi,\tau_{ne},s')$ est non nul en $s'=0$ (on a d\'ej\`a pris en compte $s_{0}$ pour d\'efinir $\tau_{ne}$). On peut encore remplacer $r(s',\tau_{ne},\psi)$ par $r(s',\tau_{ne},\psi_{nu})$ car ces fonctions ont le m\^eme ordre en $s'=0$. On  d\'ecompose cet op\'erateur en le produit de l'op\'erateur d'entrelacement standard (qui est holomorphe par positivit\'e)
$$
\tau_{ne}\vert\,\vert^{s'} \times_{(\rho,a,b)_{x}}J(St(\rho,a),-(b-1)/2,-\delta_{b})\vert\,\vert^{x}\times J(St(\rho^*,a),-(b-1)/2,-\delta'_{b})\vert\,\vert^{-x}\times \pi_{0}$$
$$
\rightarrow
$$
$$
\times_{(\rho,a,b)_{x}}J(St(\rho,a),-(b-1)/2,-\delta_{b})\vert\,\vert^{x}\times J(St(\rho^*,a),-(b-1)/2,-\delta'_{b})\vert\,\vert^{-x}\times \tau_{ne}^*\vert\,\vert^{-s'}\times \pi_{0}
$$
avec l'op\'erateur d'entrelacement
$$
\times_{(\rho,a,b)_{x}}J(St(\rho,a),-(b-1)/2,-\delta_{b})\vert\,\vert^{x}\times J(St(\rho^*,a),-(b-1)/2,-\delta'_{b})\vert\,\vert^{-x}\times \tau_{ne}^*\vert\,\vert^{-s'}\times \pi_{0}$$
$$\rightarrow
$$
$$\tau_{ne}^*\vert\,\vert^{-s'}
\times_{(\rho,a,b)_{x}}J(St(\rho,a),-(b-1)/2,-\delta_{b})\vert\,\vert^{x}\times J(St(\rho^*,a),-(b-1)/2,-\delta'_{b})\vert\,\vert^{-x}\times  \pi_{0}
$$
Le premier est non nul car c'est un op\'erateur d'entrelacement standard holomorphe et son image n'annule pas le sous-quotient de Langlands. Donc il suffit de montrer que le 2e op\'erateur d'entrelacement n'annule pas non plus le sous-quotient de Langlands et \`a ce point, il suffit de montrer qu'il est non  nul puisque le sous-quotient de Langlands est l'unique sous-module irr\'eductible de l'induite d'arriv\'ee. Or cet op\'erateur est purement un op\'erateur dans un $GL$, c'est l'op\'erateur
$$
\times_{(\rho,a,b)_{x}}J(St(\rho,a),-(b-1)/2,-\delta_{b})\vert\,\vert^{x}\times J(St(\rho^*,a),-(b-1)/2,-\delta'_{b})\vert\,\vert^{-x}\times \tau^*_{ne}\vert\,\vert^{-s'}$$
$$ \rightarrow$$
$$ \tau^*_{ne}\vert\,\vert^{-s'}\times
\times_{(\rho,a,b)_{x}}J(St(\rho,a),-(b-1)/2,-\delta_{b})\vert\,\vert^{x}\times J(St(\rho^*,a),-(b-1)/2,-\delta'_{b})\vert\,\vert^{-x}$$
normalis\'e comme expliqu\'e. Il n'y a plus qu'\`a remarquer que les p\^oles de notre normalisation sont ceux de la fonction $$\times _{(\rho,a,b)_{x}\in Jord(\psi)_{nu,<0}} \biggl(L(\tau_{ne}\times St(\rho,a),s'-(b-1)/2-x)L(\tau_{ne}\times St(\rho^*,a),s'-(b-1)/2+x)\biggr);$$
ce sont aussi les p\^oles de la fonction de la normalisation utilis\'ee en \cite{mw} I.8 pour cet op\'erateur et on y a d\'emontr\'e la non nullit\'e. D'o\`u l'assertion cherch\'ee.

On a donc montr\'e que l'image de $N_{\psi}(s_{0},\tau,\pi)$ est l'unique sous-module irr\'eductible de l'induite:
$$
\tau_{ne}^*\times_{(\rho,a,b)_{x}\in Jord(\psi_{nu,<0})}J(St(\rho,a),-(b-1)/2,-\delta_{b})\vert\,\vert^{x}\times J(St(\rho^*,a),-(b-1)/2,-\delta'_{b})\vert\,\vert^{-x}\times \pi_{0}.\eqno(1)
$$
Il faut maintenant montrer l'irr\'eductibilit\'e de $\pi_{0}$. D\'ecrivons d'abord $\tau_{e}$: on note $s^+_{0}$ et $s^-_{0}$ les demi-entiers cons\'ecutifs tels $s_{0}\in [s^-_{0},s^+_{0}]$ bien d\'efinis si $s_{0}$ n'est pas un demi-entier et on pose $s^{\pm}_{0}=s_{0}$ si $s_{0}$ est un demi-entier. Par hypoth\`ese, $s_{0}\geq 1/2$ et donc $s^{\pm}_{0}\geq 1/2$. Quand on \'ecrit $\tau_{e}$ sous forme d'induite de repr\'esentations de Steinberg tordues par un caract\`ere non unitaire, et ce caract\`ere non unitaire est soit $\vert\,\vert^{s^{+}_{0}}$ soit $\vert\,\vert^{s^-_{0}}$. Ou encore, pour fixer les notations:
$$
\tau_{e}\simeq \times_{(\rho',a',\zeta')\in {\mathcal J}'}St(\rho',a')\vert\,\vert^{s^{\zeta' }_{0}}.
$$
Si $s_{0}\in 1/2 {\mathbb N}$, dans tout ce qui suit $s^-_{0}$ n'intervient pas. En g\'en\'eral on pose pour $\zeta=\pm$,  ${\mathcal J}^{',\zeta}:=\{(\rho',a',\zeta')\in {\mathcal J}'; \zeta'=\zeta\}$ et on d\'ecompose ${\mathcal J}^{',\zeta}$ en ${\mathcal J}^{',\zeta}_{bp}\cup {\mathcal J}^{',\zeta}_{mp}$ o\`u l'indice $bp$ indique regroupe les triplets $(\rho',a',\zeta)$ tel que la repr\'esentation de $W_{F}\times SL(2,{\mathbb C})\times SL(2,{\mathbb C})$ associ\'ee au triplet $(\rho',a', 2s^{\zeta}_{0}+1)$ est de bonne parit\'e . D'o\`u une d\'ecomposition en induites:
$$
\tau\vert\,\vert^{s_{0}}\simeq \tau_{bp}^{+}\vert\,\vert^{s^+_{0}}\times \tau_{bp}^-\vert\,\vert^{s^-_{0}}\times \tau_{mp}^{+}\vert\,\vert^{s^+_{0}}\times \tau_{mp}^-\vert\,\vert^{s^-_{0}}.
$$
Pour $\zeta=\pm$ et $\zeta'=\pm$, les exposants intervenant dans le support cuspidal de  $\tau_{bp}^{\zeta}\vert\,\vert^{\pm s^\zeta_{0}}$ et de $\tau_{mp}^{\zeta'}\vert\,\vert^{\pm s^{\zeta'}_{0}}$ n'ont pas les 
m\^eme propri\'et\'e d'int\'egralit\'e: avec l'indice $bp$ le support cuspidal est de la forme $\rho'\vert\,\vert^{x}$ avec $\rho'$ necessairement autodual et $x$ est soit demi-entier non entier, soit entier, cela d\'epend de $\rho'$ et avec l'indice $mp$ une de ces conditions est n\'ecessairement oppos\'ee. Cela donne toutes les irr\'eductibilit\'es possible entre induite
$$
\tau_{bp}^{\zeta}\vert\,\vert^{\pm s^\zeta_{0}}\times \tau_{mp}\vert\,\vert^{s^{\zeta'}_{0}}, \quad 
\tau_{bp}^{\zeta}\vert\,\vert^{\pm s^\zeta_{0}}\times \tau_{mp}^*\vert\,\vert^{-s^{\zeta'}_{0}}.
$$
Au passage, on remarque que $\tau_{bp}^{\zeta}$ est n\'ecessairement autoduale. 
On consid\`ere l'op\'erateur d'entrelacement 
$$
\tau^+_{bp}\vert\,\vert^{s^+}\times \tau^-_{bp}\vert\,\vert^{s^-} \pi_{bp}\rightarrow \tau^-_{bp}\vert\,\vert^{-s^-}\times \tau^+_{bp}\vert\,\vert^{-s^+}\times \pi_{bp},
$$normalis\'e par le produit des fonctions $r(s^\zeta, \tau_{bp}^\zeta, \psi)$. Calcul\'e en $s^\zeta=s^\zeta_{0}$, il est holomorphe (cf. \ref{generalisation}).
On note $\pi_{0,bp}$ son image en $s^{\zeta}=s^\zeta_{0}$ et on admet momentan\'ement que cette image est nulle ou irr\'eductible et on conclut la preuve de la proposition. On v\'erifie que $\pi_{0}$ est l'image de l'op\'erateur d'entrelacement standard
$$
\tau_{mp}^+\vert\,\vert^{s^+}\times \tau_{mp}^-\vert\,\vert^{s^-} \times \pi^{GL}(\psi_{1/2,mp})\times \pi_{0,bp}\rightarrow \tau_{mp}^{*-}\vert\,\vert^{-s^-}\times \tau_{mp}^{*+}\vert\,\vert^{-s^+} \times \pi^{GL}(\psi_{1/2,mp})\times \pi_{0,bp},
$$
normalis\'e par le produit des fonctions $r(s^\zeta,\tau_{mp}^\zeta,\psi)$, pour $\zeta=\pm$ et calcul\'e en $s^\zeta=s^\zeta_{0}$. Le cas de mauvaise parit\'e est analogue au cas non entier (on utilise la remarque 2 de \ref{descriptiongenerale} pour avoir cette similitude): l'induite de gauche \`a un unique quotient irr\'eductible qui est aussi l'unique sous-module irr\'eductible  de l'induite de droite et cette repr\'esentation irr\'eductible intervient avec multiplicit\'e exactement 1 dans les induites \'ecrites. L'op\'erateur d'entrelacement a donc son image nulle ou co\"{\i}ncidant avec cette sous-repr\'esentation. On caract\'erise encore cette sous-repr\'esentation irr\'eductible comme le sous-module de Langlands de l'induite
$$
\tau_{mp}^*\vert\,\vert^{-s_{0}}\times_{(\rho,a,b)\in Jord(\psi_{1/2,mp})}J(St(\rho,a),-(b-1)/2,-\delta_{b})\times J(St(\rho,a),-(b-1)/2,-\delta'_{b})\times \pi_{0,bp}. \eqno(2)
$$
Il faut donc d\'emontrer que l'op\'erateur d'entrelacement est d'image non nulle. On remarque que les fonctions $r(s^\zeta,\tau_{mp}^\zeta,\psi_{bp})$ n'ont ni z\'ero ni p\^oles. Comme $\psi_{mp}=\psi_{1/2,mp}^+\oplus \psi_{1/2,mp}$, le facteur de normalisation est encore le produit de facteurs de normalisation \`a la Langlands Shahidi pour les op\'erateurs d'entrelacement, dans les GL convenables:
$$
\tau_{mp}^\zeta\vert\,\vert^{s^\zeta} \times \pi^{GL}(\psi_{1/2,mp}) \rightarrow \pi^{GL}(\psi_{1/2,mp})\times \tau_{mp}^\zeta\vert\,\vert^{s^\zeta}\eqno(3)$$
$$
 \pi^{GL}(\psi_{1/2,mp})\times \tau_{mp}^{*\zeta}\vert\,\vert^{-s^\zeta}\rightarrow \tau_{mp}^{*\zeta}\vert\,\vert^{-s^\zeta}\times \pi^{GL}(\psi_{1/2,mp}).\eqno(4)
 $$
 L'op\'erateur d'entrelacement cherch\'ee est donc le produit de (3), avec l'entrelacement $$\tau_{mp}^\zeta\vert\,\vert^{s^\zeta}\times \pi_{0,bp}\rightarrow \tau_{mp}^{*\zeta}\vert\,\vert^{-s^\zeta}\times \pi_{0,bp}$$
 qui est l'op\'erateur d'entrelacement standard et l'op\'erateur d'entrelacement (4). Tous ces op\'erateurs d'entrelacement n'annulent pas le sous-quotient de Langlands. On a donc montr\'e que $\pi_{0}$ est le sous-module de Langlands de (2).

\

La nullit\'e ou l'irr\'eductibilit\'e de $\pi_{0,bp}$ est une g\'en\'eralisation  de \ref{imagebonneparite}(ii) mais il faut reprendre la d\'emonstration; remarquons que dans les cas qui nous int\'eressent vraiment $s'_{0}$ est un demi-entier et $s^+_{0}=s^-_{0}$ avec les notations ci-dessus. Pour simplifier on suppose donc que $s'_{0}$ est un demi-entier que l'on note $s_{0}$ pour \^etre compatible aux notations d\'ej\`a utilis\'ees.

La d\'emonstration de \ref{imagebonneparite} est bas\'ee sur des diagrammes commutatifs, o\`u on peut tranquillement remplacer $St(\rho,a_{0})\vert\,\vert^{s_{0}}$ par $\tau_{bp}\vert\,\vert^{s_{0}}$, et des propri\'et\'es d'holomorphie et de module de Jacquet o\`u on ne peut pas remplacer $St(\rho,a_{0})\vert\,\vert^{s_{0}}$ car on a l'hypoth\`ese $A\geq A_{0}$ (avec les notations de loc.cite) alors qu'ici $A_{0}$ n'est pas uniquement d\'efini. On fait une r\'ecurrence sur le nombre de repr\'esentations distincts intervenant dans ${\mathcal J}_{bp}$; on initialise, trivialement, la r\'ecurrence avec le cas o\`u ${\mathcal J}_{bp}$ est vide. On fixe $\rho,a_{0}$ tel que $\tau_{bp}\simeq St(\rho,a_{0})\times \cdots \times St(\rho,a_{0})\times \tau'$, o\`u les $\cdots$ cachent des copies de $St(\rho,a_{0})$ et o\`u $\tau'$ est une induite de repr\'esentations de Steinberg de la forme $St(\rho'',a'')$ o\`u pour touts les $a''$ intervenant $a''<a_{0}$. Ainsi pour un tel $a''$, on a certainement $(a''-1)/2+s_{0}<(a_{0}-1)/2+s_{0}$. On montre que l'on peut reprendre la d\'emonstration de \ref{imagebonneparite}. On montre d'abord  que l'induite
$$
St(\rho,a_{0})\vert\,\vert^{-s_{0}}\times \cdots \times St(\rho,a_{0})\vert\,\vert^{-s_{0}}\times N_{\psi}(s_{0},\tau',\pi_{bp})(\tau'\vert\,\vert^{s_{0}}\times \pi_{bp})$$a un unique sous-module irr\'eductible \`a condition que pour tout $(\rho,a,b)\in Jord(\psi_{bp})$ l'in\'egalit\'e $(a+b)/2-1\geq (a_{0}-1)/2+s_{0}$ force $\vert a-b\vert/2>>0$ et ce sous-module irr\'eductible intervient, alors, avec multiplicit\'e 1 comme sous-quotient irr\'eductible. C'est un calcul de module de Jacquet: d'abord on a la propri\'et\'e de nullit\'e suivante:

pour tout $x\in [(a_{0}-1)/2-s_{0}, -(a_{0}-1)/2-s_{0}]$, $Jac_{x, \cdots, -(a_{0}-1)/2-s_{0}}N_{\psi}(s_{0},\tau',\pi_{bp})(\tau'\vert\,\vert^{s_{0}}\times \pi_{bp})=0$.  La repr\'esentation $\rho\vert\,\vert^{\pm (a_{0}-1)/2+s_{0})}$ n'est pas dans le support cuspidal de $\tau'\vert\,\vert^{s_{0}}$ par l'hypoth\`ese faite ci-dessus $a''<a_{0}$. Ainsi la non nullit\'e du module de ce module de Jacquet force l'existence de $x'\in [x,-(a_{0}-1)/2-s_{0}]$ tel que $Jac_{x',\cdots,-(a_{0}-1)/2-s_{0}}\pi_{bp}\neq 0$; mais cela est exclu par l'hypoth\`ese faite sur $Jord(\psi_{bp})$ et \ref{proprietedujac}.

Ainsi la repr\'esentation $$
\biggl(St(\rho,a_{0})\vert\,\vert^{-s_{0}}\times \cdots \times St(\rho,a_{0})\vert\,\vert^{-s_{0}}\biggr)\otimes \biggl(N_{\psi}(s_{0},\tau',\pi_{bp})(\tau'\vert\,\vert^{s_{0}}\times \pi_{bp})\biggr)
$$intervient avec multiplicit\'e 1 dans le module de Jacquet de l'induite consid\'er\'e. Or cette repr\'esentation est irr\'eductible par l'hypoth\`ese de r\'ecurrence d'o\`u notre assertion par r\'eciprocit\'e de Frobenius et exactitude du foncteur de Jacquet. 
On conclut \`a l'irr\'eductibilit\'e de l'image de tout morphisme non nul de $$St(\rho,a_{0})\vert\,\vert^{s_{0}}\times \cdots \times St(\rho,a_{0})\vert\,\vert^{s_{0}}\times N_{\psi}(s_{0},\tau',\pi_{bp})(\tau'\vert\,\vert^{s_{0}}\times \pi_{bp})$$
dans l'induite \'ecrite. On a donc le r\'esultat dans ce cas. Il faut maintenant enlever l'hypoth\`ese faite.

On peut donc reprendre la m\'ethode de descente de \ref{imagebonneparite} qui fait passer, avec les notations de loc.cite de $\pi_{>}$ \`a $\pi$; ici on fait $\pi_{bp,>}$ \`a $\pi_{bp}$ en fixant $(\rho,A,B,\zeta)\in Jord(\psi_{bp})$ mais on peut maintenant supposer que $A\geq (a_{0}-1)/2+s_{0}$. Donc on a aussi $A\geq (a'-1)/2+s_{0}$ pour tout $(\rho,a')\in {\mathcal J}_{bp}$. La preuve de loc.cite s'\'etend alors pour donner l'assertion.

\subsubsection{Description pr\'ecise de l'image dans certains cas\label {parametregeneralisation}}
On garde les notations pr\'ec\'edentes, $\psi$, $\tau$ et $\pi$ dans le paquet associ\'e \`a $\psi$. On note $\psi_{\tau}$ la repr\'esentation de $W_{F}\times SL(2,{\mathbb C})$ qui correspond \`a $\tau$ par la param\'etrisation de Langlands; c'est une repr\'esentation \`a valeurs dans un groupe $GL(d_{\tau},{\mathbb C})$ pour $d_{\tau}$ convenable. On fixe encore un r\'eel, $s_{0}$, qui ici est un demi-entier donc de la forme $(b_{0}-1)/2$ avec $b_{0}\geq 2$. 

On suppose que $\psi_{\tau}$ tensoris\'e par la repr\'esentation de dimension $b_{0}$ de $SL(2,{\mathbb C})$, ce qui donne donc une repr\'esentation de $W_{F}\times SL(2,{\mathbb C})\times SL(2,{\mathbb C})$,  est  autoduale et \`a valeurs dans un groupe de m\^eme type que le groupe dual de $G$  (c'est-\`a-dire est de bonne parit\'e) et on suppose aussi que

\

$b_{0}=2$ ou $b_{0}>2$ et  $\psi$ contient la sous-repr\'esentation $\psi_{\tau}\otimes sp_{b_{0}-2}$ de $W_{F}\times SL(2,{\mathbb C})\times {\mathbb C}$ o\`u $sp_{b_{0}-2}$ est la repr\'esentation de $SL(2,{\mathbb C})$ de dimension $b_{0}-2$. 

\

Remarquons que ces hypoth\`eses entra\^{\i}nent que $\tau$ est autoduale. On note $\psi^+$ le morphisme qui s'obtient en ajoutant \`a $\psi$ la repr\'esentation $\psi_{\tau}\otimes sp_{2}$ si $b_{0}=2$ et si $b_{0}>2$, en rempla\c{c}ant $\psi_{\tau}\otimes sp_{b_{0}-2}$ par $\psi_{\tau}\otimes sp_{b_{0}}$. Le morphisme $\psi^+$ est de m\^eme nature que le morphisme $\psi$.

On note encore $\psi_{bp}$ la somme des sous-repr\'esentations irr\'eductibles de $\psi$ ayant bonne parit\'e et $\psi^+_{bp}$ est d\'efini de fa\c{c}on analogue. Le passage de $\psi_{bp}$ \`a $\psi^+_{bp}$ se d\'ecrit ainsi: on note $\psi_{\tau,bp}$ la somme des sous-repr\'esentations irr\'eductible de $\psi_{\tau}$ telles que tensoris\'ees avec la repr\'esentation de dimension $b_{0}$ de $SL(2,{\mathbb C})$ elles sont de bonne parit\'e en tant que repr\'esentation de $W_{F}\times SL(2,{\mathbb C})\times SL(2,{\mathbb C})$. On passe de $\psi_{bp}$ \`a $\psi_{bp}^+$ en enlevant \`a $\psi_{bp}$ la repr\'esentation  $\psi_{\tau,bp}\otimes sp_{b_{0}-2}$ (quand $b_{0}>2$) et en lui ajoutant $\psi_{\tau,bp}\otimes sp_{b_{0}}$. Si $\psi_{\tau,bp}$ est irr\'eductible c'est le processus \'etudi\'e en \ref{imagebonneparite}; dans ce cas pour toute repr\'esentation $\pi_{bp}$ dans le paquet associ\'e \`a $\psi_{bp}$, on a d\'efini une repr\'esentation $\pi_{bp}^+$. Dans le cas g\'en\'eral, on a une succession de situation \'el\'ementaire \'etudi\'ee en loc. cite; on d\'efinit successivement le passage de $\pi_{bp}$ \`a $\pi_{bp}^+$; a priori cela d\'epend du choix de l'ordre dans lequel on fait ces passages mais comme on va caract\'eriser $\pi_{bp}^+$ intrins\`equement, cela prouvera l'ind\'ependance des choix. On pourrait donner les param\`etres $\underline{t}^+, \underline{\eta}^+$ comme on l'a fait en loc. cite mais il y a une difficult\'e dans le cas o\`u $\psi_{\tau,bp}$ a de la multiplicit\'e et on ne le fait pas ici.

On a donc d\'efini le passage $\pi_{bp}\mapsto \pi_{bp}^+$ du paquet de repr\'esentations associ\'e \`a $\psi_{bp}$ dans le paquet de repr\'esentations associ\'e \`a $\pi_{bp}^+$. Cela d\'efinit une application $\pi\in \Pi(\psi)\mapsto \pi^+\in \Pi(\psi^+)$ avec les notations de \ref{descriptiongenerale}:
$$
\pi\simeq \pi^{GL}{(\psi_{nu,<0})}\times \pi^{GL}(\psi_{1/2,mp})\times \pi_{bp}\mapsto \pi^{GL}(\psi^+_{\nu,<0})\times \pi^{GL}(\psi_{1/2,mp}^+)\times \pi^+_{bp}.
$$
Dans les notations ci-dessus, il n'est pas n\'ecessaire que $\psi^+_{1/2,mp}$ soit li\'e \`a $\psi_{1/2,mp}$ puisque les induites ne d\'ependent pas du d\'ecoupage de $\psi^+_{mp}$ et de $\psi_{mp}$. 
\begin{thm}Pour toute repr\'esentation (non nulle) $\pi$ dans le paquet associ\'e \`a $\psi$, l'image de l'op\'erateur d'entrelacement normalis\'e $N_{\psi}(s,\tau,\pi)$ en $s=(b_{0}-1)/2$ est exactement $\pi^+$; c'est-\`a-dire est nul si $\pi^+=0$ et est irr\'eductible isomorphe \`a $\pi^+$ sinon.
\end{thm}
La d\'emonstration de \ref{qualitatifgeneral} et en particulier (1) et (2) de cette r\'ef\'erence, va permettre de ramener le th\'eor\`eme au cas o\`u $\psi=\psi_{bp}$. Montrons d\'ej\`a cela. On admet ici que le $\pi_{0,bp}$ de cette preuve est $\pi^+_{bp}$ et on montre le th\'eor\`eme. Puisqu'ici on suppose d\`es le d\'epart que l'on calcule l'image en un point demi-entier, il est plus simple de d\'ecrire $\tau_{ne}$ de loc.cite. On note $\tau_{nu}$ les repr\'esentations de Steinberg tordu par un caract\`ere de la forme $\vert\,\vert^{x}$ avec $x\in ]-1/2,1/2[$, $x\neq 0$ tel que $\tau$ soit l'induite de $\tau_{nu}$ avec une induite de repr\'esentation de Steinberg. On d\'ecompose encore cette induite avec des notations que l'on esp\`ere claires:
$$
\tau\simeq \tau_{nu}\times \tau_{mp}\times \tau_{bp}.
$$
On rappelle pour \'eviter les confusions que la d\'efinition de bonne parit\'e tient compte de $b_{0}$ (de sa parit\'e) et ne d\'epend pas uniquement de $\tau$. Le $\tau_{ne}$ de loc.cite est ici $\tau_{nu}\vert\,\vert^{(b_{0}-1)/2}$ et $\tau_{mp}$ est le m\^eme.

Avec les r\'ef\'erences que l'on vient de donner, on a montr\'e que l'image cherch\'ee est l'unique sous-module irr\'eductible de l'induite (les notations sont celles de loc. cite) $$ \tau_{nu}\vert\,\vert^{-(b_{0}-1)/2}\times \eqno(1)$$
$$
\times_{(\rho,a,b)_{x}\in Jord(\psi_{nu,<0})} \biggl(J(St(\rho,a),-(b-1)/2,-\delta_{b})\vert\,\vert^{x}\times J(St(\rho^*,a),-(b-1)/2,-\delta'_{b})\vert\,\vert^{-x}\biggr)
$$
$$\times \tau_{mp}\vert\,\vert^{-(b_{0}-1)/2}\times \eqno(2)$$
$$
\times_{(\rho,a,b)\in Jord(\psi_{1/2,mp})}\biggl(J(St(\rho,a),-(b-1)/2,-\delta_{b})\times J(St(\rho^*,a),-(b-1)/2,-\delta'_{b})\biggr)$$
$$
\times \pi_{bp}^+.
$$
On r\'ecrit $\tau_{nu}$ comme induite de repr\'esentation de Steinberg tordu par des caract\`eres non unitaires en utilisant l'autodualit\'e sous la forme $$\tau_{nu}\simeq \times_{(\rho',a')_{z'}\in Jord(\psi_{\tau,nu,<0})}\biggl(St(\rho',a')\vert\,\vert^{z'}\times St(\rho^{'*},a')\vert\,\vert^{-z'}\biggr).
$$
Si $b_{0}=2$, on a donc $\delta_{b_{0}}=\delta'_{b_{0}}=1/2$ et avec des notations un peu compliqu\'ee (mais g\'en\'erales)
$
\tau_{nu}\vert\,\vert^{-(b_{0}-1)/2}\simeq$
$$ \times_{(\rho',a')_{z'}\in Jord(\psi_{\tau,nu,<0})}J(St(\rho',a'),-(b_{0}-1)/2,-\delta_{b_{0}})\times J(St(\rho',a'),-(b_{0}-1)/2,-\delta'_{b_{0}}).\eqno(3)
$$
Si $b_{0}>2$ par hypoth\`ese pour tout $(\rho',a')_{z'}\in Jord(\psi_{\tau,nu,<0})$ il existe $(\rho',a',b_{0}-2)_{z'}\in Jord(\psi_{nu,<0})$ (avec la m\^eme multiplicit\'e). On peut donc remplacer la ligne (1)
ci-dessus induite avec la ligne suivante par le sous-module irr\'eductible o\`u on remplace, pour chaque $(\rho',a')_{z'}\in Jord(\psi_{\tau,nu,<0})$
$$
St(\rho',a')\vert\,\vert^{-(b_{0}-1)/2+z'}\times St(\rho^{'*},a')\vert\,\vert^{-(b_{0}-1)/2-z'}\times $$
$$J(St(\rho',a'),-(b_{0}-3)/2,-\delta_{b_{0}})\vert\,\vert^{z'}\times J(St(\rho^{'*},a'),-(b_{0}-3)/2,-\delta'_{b_{0}})\vert\,\vert^{-z'}$$
par son sous-module irr\'eductible (c'est le sous-module de Langlands de la situation)
$$
J(St(\rho',a'),-(b_{0}-1)/2,-\delta_{b_{0}})\vert\,\vert^{z'}\times J(St(\rho^{'*},a'),-(b_{0}-1)/2,-\delta'_{b_{0}})\vert\,\vert^{-z'}\eqno(4)$$
On fait la m\^eme op\'eration avec $\tau_{mp}$ mais on commence par faire un choix de d\'ecomposition $\tau_{mp}\simeq \tau_{1/2,mp}\times \tau^*_{1/2,mp}$. Pour simplifier la d\'emonstration, on fait ce choix de telle sorte que, si $b_{0}>2$, pour tout $(\rho',a')\in Jord(\psi_{\tau,1/2,mp})$ (avec des notations \'evidentes), $(\rho',a',b_{0}-2)\in Jord(\psi_{1/2,mp})$  avec une multiplicit\'e sup\'erieure ou \'egale \`a celle de $(\rho',a')\in Jord(\psi_{\tau,1/2, mp})$. Ensuite on conclut comme ci-dessus que soit $b_{0}=2$ et que l'on \'ecrit $\tau_{mp}\vert\,\vert ^{-(b_{0}-1)/2}$ de fa\c{c}on analogue \`a (3) sans les torsions non unitaires soit $b_{0}>2$ et on peut remplacer l'induite de la ligne (2) ci-dessus avec la ligne qui la suit par l'analogue de (4) sans les torsions non unitaires. 

Le th\'eor\`eme r\'esulte alors de la proposition de  \ref{descriptiongenerale} et de la remarque 2 de loc.cite. Il nous reste donc \`a traiter le cas o\`u $\psi=\psi_{bp}$ qui va \^etre une cons\'equence simple de \ref{imagebonneparite} et des vertus des normalisations.
Ici \'ecrivons $\tau_{bp}\times_{j\in [1,v]}St(\rho_{j},a_{j})$. On consid\`ere l'op\'erateur d'entrelacement
$$
\times_{j\in [1,v]}St(\rho_{j},a_{j})\vert\,\vert^{s_{0}}\times \pi_{bp}\rightarrow \times_{j\in [v,1]}St(\rho_{j},a_{j})\vert\,\vert^{-s_{0}}\times \pi_{bp}.
$$
On fait d\'ependre cet op\'erateur de nombres complexes $s'_{j}$ pour $j\in [1,v]$ voisins de 0 et que l'on \'evaluera en $0$:
$$
\times_{j\in [1,v]}St(\rho_{j},a_{j})\vert\,\vert^{s'_{j}+(b_{0}-1)/2}\times \pi_{bp}\rightarrow \times_{j\in [v,1]}St(\rho_{j},a_{j})\vert\,\vert^{-s'_{j}-(b_{0}-1)/2}\times \pi_{bp}.\eqno(1)
$$
En \ref{imagebonneparite}, on a \'etudi\'e successivement les op\'erateurs (la diff\'erence entre (1) et (2) et le changement de $[v,1]$ en $[1,v]$ dans le membre de droite):
$$
\times_{j\in [1,v]}St(\rho_{j},a_{j})\vert\,\vert^{s'_{j}+(b_{0}-1)/2}\times \pi_{bp}\rightarrow
\times_{j\in [1,v]}St(\rho_{j},a_{j})\vert\,\vert^{-s'_{j}-(b_{0}-1)/2}\times \pi_{bp}\eqno(2)
$$
avec une normalisation  diff\'erente: en effet \`a chaque \'etape $\psi_{bp}$ est modifi\'e on remplace progressivement $(\rho_{j},a_{j},b_{0}-2)$ par $(\rho_{j},a_{j},b_{0})$; pour tout $\ell\in [1,v]$, on note $\psi_{bp}^{+_{\ell}}$ le morphisme o\`u on a effectu\'e ce changement pour tout $j\in ]\ell,v]$, donc $\psi_{bp}^{+_{v}}:=\psi_{bp}$. On a donc normalis\'e par $\times_{\ell \in [1,v]}r(\psi_{bp}^{+_{\ell}},St(\rho_{\ell},a_{\ell}),s'_{\ell}+(b_{0}-1)/2)$. Ainsi la contribution de $(\rho,a_{j},b_{0}-2)$ \`a $r(\psi_{bp},St(\rho,a_{\ell}),s'_{\ell}+(b_{0}-1)/2)$ qui \'etait  $$\frac{L(St(\rho,a_{\ell})\times St(\rho,a_{j}),s'_{\ell}+(b_{0}-1)/2-(b_{0}-3)/2)}{L(St(\rho,a_{\ell})\times St(\rho,a_{j}),s'_{\ell}+(b_{0}-1)/2+(b_{0}-1)/2)}=$$
$$\frac{
L(St(\rho,a_{\ell})\times St(\rho,a_{j}),s'_{\ell}+1)}{L(St(\rho,a_{\ell})\times St(\rho,a_{j}),s'_{\ell}+b_{0}-1)}$$
devient $L(St(\rho,a_{\ell})\times St(\rho,a_{j}),s'_{\ell})/L(St(\rho,a_{\ell})\times St(\rho,a_{j}),s'_{\ell}+b_{0}-1)$. Aux d\'enominateurs pr\`es qui n'introduisent ni z\'ero ni p\^ole on  passe donc de $r(\psi_{bp},St(\rho,a_{\ell}),s'_{j}+(b_{0}-1)/2)$ \`a $r(\psi_{bp}^{+_{\ell}},St(\rho_{\ell},a_{\ell}),s'_{\ell}+(b_{0}-1)/2)$ en multipliant par
$$
\times_{j\in ]\ell,v]}L(St(\rho_{\ell},a_{\ell})\times St(\rho_{j},a_{j}),s'_{\ell})/L(St(\rho_{\ell},a_{\ell})\times St(\rho_{j},a_{j}),s'_{\ell}+1).
$$
Ceci est exactement la normalisation de Langlands Shahidi qui permet de passer du deuxi\`eme membre de (1) \`a celui de (2) par des op\'erateurs d'entrelacement holomorphes bijectifs (\cite{mw}I.9). Ainsi l'image de l'op\'erateur (1) et exactement l'image de l'op\'erateur (2). En particulier, elle est nulle ou irr\'eductible et se calcule par les formules explicites de \ref{imagebonneparite}.

\section{Places archim\'ediennes, les repr\'esentations ayant de la cohomologie}
Dans ce paragraphe on suppose que le corps de base est archim\'edien et on va m\^eme se limiter au cas o\`u $F$, le corps de base est r\'eel. Et on ne consid\`ere que des repr\'esentations ayant de la cohomologie et qui sont composantes locales de formes automorphes de carr\'e int\'egrable donc le $\pi$ de tout ce travail est cohomologique, pour un bon syst\`eme de coefficients et est unitaire. Quand on consid\'erera des induites, on se place en des points o\`u le caract\`ere infinit\'esimal est entier r\'egulier; c'est surtout r\'egulier qui nous importe. Ces restrictions sont \'evidemment g\^enantes; elles sont d\^ues au fait que l'on ne conna\^{\i}t pas les paquets d'Arthur en g\'en\'eral dans le cas archim\'edien et aussi au fait que les op\'erateurs d'entrelacement sont plus difficile \`a \'etudier aux places archim\'ediennes car il y a plus de r\'eductibilit\'e qu'aux places p-adiques.

\subsection{Les op\'erateurs d'entrelacement standards\label{descriptionstandard}}
Pr\'ecis\'ement, on fixe une repr\'esentation irr\'eductible unitaire $\pi$ ayant de la cohomologie pour un bon syst\`eme de coefficients et une induite de s\'eries discr\`etes, $\tau$, de $GL(n,F)$; c'est-\`a-dire que $\tau$ est une repr\'esentation temp\'er\'ee. On consid\`ere l'op\'erateur d'entrelacement standard $$M(s,\tau, \pi): \tau\vert\,\vert^{s}\times \pi \rightarrow  \tau^*\vert\,\vert^{-s}\times \pi.$$
\begin{prop} On fixe $s_{0}$ un r\'eel strictement positif tel que le caract\`ere infinit\'esimal de l'induite $\delta\vert\,\vert^{s_{0}}\times \pi$ soit r\'egulier. Alors $M(s,\tau,\pi)$ est holomorphe en $s=s_{0}$. La repr\'esentation induite $\tau\vert\,\vert^{s_{0}}\times \pi$ a un unique quotient irr\'eductible et l'image de $M(s_{0},\tau, \pi)$ est irr\'eductible isomorphe \`a ce quotient irr\'eductible.
\end{prop}
C'est \'evidemment l'hypoth\`ese de r\'egularit\'e du caract\`ere infinit\'esimal coupl\'e au fait que gr\^ace \`a \cite{vz} on conna\^{\i}t les param\`etres de Langlands de $\pi$ qui rend tout tr\`es simple. Cette proposition est d\'emontr\'ee \`a l'issue des sous-sections suivantes.

\subsubsection{Caract\`ere infinit\'esimal\label{caractereinfinitesimal}}

On \'ecrit d'abord le caract\`ere infinit\'esimal d'une repr\'esentation ayant de la cohomologie; cela est fait pour les groupes unitaires et des groupes orthogonaux en \cite{bc}. On reprend \`a peu pr\`es le m\^eme point de vue. Ici on voit vraiment les repr\'esentations dans des paquets, ceux d'Adams-Johnsons, dont il n'est pas prouv\'e que ce sont les paquets d'Arthur, sauf dans le cas des groupes unitaires (\cite{aj} et \cite{johnson}); ces paquets sont param\'etr\'es par des morphismes $\psi$ de $W_{F}\times SL(2,{\mathbb C})$ dans le $L$-groupe de G; ici on consid\`ere bien tout le $L$-groupe et pas seulement sa composante alg\'ebrique. On peut encore d\'ecomposer ces morphismes en morphismes \'el\'ementaires, un morphisme \'el\'ementaire param\'etrant un couple form\'e d'une s\'erie discr\`ete, $\delta$, d'un groupe $GL(m,{\mathbb R})$ et une repr\'esentation irr\'eductible de dimension finie de $SL(2,{\mathbb C})$ param\'etr\'e par sa dimension $b$; ici n\'ecessairement $m=1$ ou $2$. De plus, on sait que $\delta^*\simeq \delta$ et il y a une condition de parit\'e sur $b$ qui d\'epend de $\delta$. Si $m=1$ $\delta$ est l'un des 2 caract\'er\`es quadratiques de ${\mathbb R}$; le caract\`ere infinit\'esimal du carat\`ere de $GL(b,{\mathbb R})$ induit par $\delta$ est \'evidemment d\'etermin\'e par la collection d'entiers: $((b-1)/2, \cdots, -(b-1)/2)$; c'est-\`a-dire le segment $[(b-1)/2,-(b-1)/2]$. Si $m=2$, $\delta$ est une s\'erie discr\`ete unitaire de $GL(2,{\mathbb C})$ et on note $a_{\delta}$ l'entier tel que le caract\`ere infinit\'esimal de $\delta$ soit $(a_{\delta}-1)/2,-(a_{\delta}-1)/2$. Le caract\`ere infinit\'esimal de la repr\'esentation $Speh(\delta,b)$ de $GL(2b,{\mathbb R})$ est alors la collection de demi-entiers, form\'ee de l'union de 2 segments, le segment $[(b-1)/2,-(b-1)/2]$ translat\'e par $(a_{\delta}-1)/2$ et le m\^eme segment translat\'e par $-(a_{\delta}-1)/2$, c'est-\`a-dire:
$$
\{(a_{\delta}-1)/2+[(b-1)/2,-(b-1)/2] \}\cup \{-(a_{\delta}-1)/2+[(b-1)/2,-(b-1)/2]\}.
$$
On revient \`a $\psi$ qui est la somme de ces morphismes \'el\'ementaires et la repr\'esentation de $\pi^{GL}(\psi)$ du $GL$ convenable associ\'e \`a $\psi$ est l'induite des repr\'esentations $Speh(\delta,b)$ comme pr\'ec\'edemment; son caract\`ere infinit\'esimal est l'union des segments que l'on vient de d\'ecrire. On transf\`ere de fa\c{c}on compl\`etement \'evidente ce caract\`ere du centre de l'alg\`ebre enveloppante du $GL$ en un caract\`ere du centre de l'alg\`ebre enveloppante de $G$ et toute repr\'esentation dans le paquet associ\'e par \cite{aj} \`a $\psi$ a ce caract\`ere infinit\'esimal. La r\'egularit\'e du caract\`ere infinit\'esimal assure que ces segments ne se coupent pas. En particulier pour au plus un couple intervenant dans $\psi$, $\delta$ est un caract\`ere. 
\subsubsection{Pr\'esentation du quotient de Langlands\label{proprietelanglands}}
On exploite la r\'egularit\'e du caract\`ere infinit\'esimal sous la forme du lemme suivant: pour tout $(\delta,b)$ intervenant dans $\psi$, on suppose donn\'e $\ell_{\delta,b}\in [0,b]$ un entier de m\^eme parit\'e que $b$.
On note ici, pour $d_{1}\geq d_{2}$, $Speh(\delta,d_{1},d_{2})$ le quotient de Langlands de l'induite $\delta\vert\,\vert^{d_{1}}, \cdots \delta \vert\,\vert^{d_{2}}$ pour le GL convenable. On suppose aussi que $\pi_{0}$ est une s\'erie discr\`ete d'un groupe de m\^eme type que $G$ de rang plus petit et que l'induite:
$$
\pi':=\times_{(\delta,b)} Speh(\delta,(b-1)/2,(\ell_{\delta,b}+1)/2)\times \pi_{0}$$
a le caract\`ere infinit\'esimal associ\'e \`a $\psi$ tel que d\'ecrit pr\'ec\'edemment et que le caract\`ere infinit\'esimal de l'induite $\tau\vert\,\vert^{s_{0}}\times \pi'$ est r\'egulier. 
\begin{lem} La repr\'esentation $\pi'$ et l'induite $\tau\vert\,\vert^{s_{0}}\times \pi'$ ont toutes deux un unique quotient irr\'eductible et c'est le sous-quotient de Langlands de ces repr\'esentations.
\end{lem}
Pour pr\'eciser, la notion de sous-quotient de Langlands, disons qu'en rempla\c{c}ant les repr\'e\-sen\-tations irr\'eductibles $Speh(\delta,(b-1)/2,(\ell_{b}+1)/2)$ par une induite de s\'erie discr\`ete tordu et en r\'eordonnant pour \^etre dans la chambre de Weyl positive, on obtient une repr\'esentation telle que $\tau\vert\,\vert^{s_{0}}\times \pi'$ en soit un quotient.

Soient $(\delta,b),(\delta',b')$ intervenant dans $\psi$ et soit $d<b$ (resp. $d'<b'$) un entier de m\^eme parit\'e que $b$ (resp. $b'$).  Supposons d'abord que $\delta'$ soit un caract\`ere. Alors $\delta$ n'en est pas un et on a s\^urement
$$
(a_{\delta}-1)/2-(b-1)/2> (b'-1)/2; \quad -(b'-1)/2>-(a_{\delta}-1)/2+(b-1)/2.
$$
Pour tout $d_{1}\leq d_{2}\in [0,b]$ et $d'_{1}<d'_{2}\in [0,b']$, on a l'irr\'eductibilit\'e de $\tilde{\tau}:=Speh(\delta,(d_{1}-1)/2,(d'_{1}-1)/2)\times Speh (\delta',(d'_{1}-1)/2,(d'_{2}-1)/2)$ d'apr\`es \cite{mw} I.9. Supposons maintenant que ni $\delta$ ni $\delta'$ ne sont des caract\`eres. Par sym\'etrie on suppose que $(a_{\delta}-1)/2+(b-1)/2>(a_{\delta'}-1)/2+(b'-1)/2$, d'o\`u aussi $-(a_{\delta}-1)/2-(b-1)/2< -(a_{\delta'}-1)/2-(b'-1)/2$. On a alors par r\'egularit\'e:
$$
(a_{\delta}-1)/2-(b-1)/2>(a_{\delta'}-1)/2+(b'-1)/2;  -(a_{\delta'}-1)/2-(b'-1)/2 >-(a_{\delta}-1)/2+(b'-1)/2.
$$
Et on peut encore appliquer \cite{mw} I.9 pour avoir l'irr\'eductibilit\'e de l'analogue de $\tilde{\tau}$. 

On consid\`ere $\tau$ et $s_{0}$ comme dans l'\'enonc\'e mais on suppose que $\tau$ est une s\'erie discr\`ete ou un caract\`ere quadratique et on note $a_{\tau}$ le param\`etre de $\tau$. Supposons que $\tau$ soit une s\'erie discr\`ete et non un caract\`ere; on reprend les notations $\delta,b,d_{1},d_{2}$ d\'ej\`a introduites et on suppose aussi que $\delta$ est une s\'erie discr\`ete. Ici, on a soit l'irr\'eductibilit\'e de $\tau\vert\,\vert^{s_{0}}\times Speh(\delta,(d_{1}-1)/2,(d_{2}-1)/2)$ parce que soit $a_{\tau}+b+2s_{0}+1$ n'est pas un entier ou parcequ'il n'a pas la bonne parit\'e et sinon on regarde la place des \'el\'ements du couple
$$
(a_{\tau}-1)/2+s_{0},-(a_{\tau}-1)/2+s_{0}\eqno(1)
$$
par rapport au 2 segments:
$$
(a_{\delta}-1)/2+[(b-1)/2,-(b-1)/2]; -(a_{\delta}-1)/2+[(b-1)/2,-(b-1)/2].\eqno(2)
$$
Une position n'est pas possible: $(a_{\tau}-1)/2+s_{0}<(a_{\delta}-1)/2-(b-1)/2$ et $-(a_{\tau}-1)/2+s_{0}<-(a_{\delta}-1)/2-(b-1)/2$ car en ajoutant les 2 in\'egalit\'es, on aurait $2s_{0}<b-1$ et $2s_{0}\leq 0$, car $s_{0}$ est maintenant un demi-entier, ce qui est contradictoire avec l'hypoth\`ese $s_{0}>0$. Il reste 2 possibilit\'es: soit les 2 segments de (2) sont \`a l'int\'erieur du segment $[(a_{\tau}-1)/2+s_{0},-(a_{\tau}-1)/2+s_{0}]$ ou \`a l'inverse le segment $[(a_{\tau}-1)/2+s_{0},-(a_{\tau}-1)/2+s_{0}]$ est \`a l'int\'erieur du segment $[(a_{\delta}-1)/2-(b-1)/2,-(a_{\delta}-1)/2+(b-1)/2]$ soit 
$$
(a_{\tau}-1)/2+s_{0}>(a_{\delta}-1)/2+(b-1)/2; -(a_{\tau}-1)/2+s_{0}>-(a_{\delta}-1)/2+(b-1)/2.
$$
Dans cette derni\`ere \'eventualit\'e, en ajoutant les in\'egalit\'es, on trouve $2s_{0}>(b-1)$ ou encore $s_{0}>(b-1)/2$. Dans le premier cas les induites
$
\tau\vert\,\vert^{s_{0}}\times Speh(\delta,d_{1},d_{2})$
sont irr\'eductibles d'apr\`es \cite{mw} I.9 et dans le 2e cas ces induites sont avec des exposants dans la chambre de Weyl strictement positive. On peut donc mettre $\tau\vert\,\vert^{s_{0}}$ a la bonne place, ce qui veut dire: que $\tau\vert\,\vert^{s_{0}}\times \pi'$ est un quotient de l'induite $$\times_{(\delta,b);(\ell_{\delta,b}+1)/2> s_{0}} Speh(\delta,(b-1)/2,(\ell_{\delta,b}+1)/2)$$
$$\times_{(\delta,b);(\ell_{\delta}+1)/2\leq s_{0}}Speh(\delta,(b-1)/2,\tilde{s}_{0})\times \tau\vert\,\vert^{s_{0}}\times_{(\delta,b);(\ell_{\delta}+1)/2\leq s_{0}}Speh(\delta,\tilde{s}_{0}-1,\ell_{\delta,b})\times
$$
$$
\times_{(\delta,b); (b-1)/2 \leq s_{0}}Speh(\sigma,(b-1)/2,(\ell_{\delta,b}+1)/2)\times \pi_{0},
$$
o\`u $\tilde{s}_{0}$ d\'epend de $(\delta,b)$ et est le plus petit demi-entier sup\'erieur strictement \`a $s_{0}$ et tel que $2\tilde{s}_{0}\equiv b[2]$.

Si $\tau$ n'est pas une s\'erie discr\`ete mais une repr\'esentation temp\'er\'ee plus g\'en\'erale, on proc\`ede de la m\^eme fa\c{c}on pour d\'eplacer les constituants de $\tau$ \'ecrite comme induite. Cela permet de conclure puisque l'on s'est essentiellement mis dans la fermeture de chambre de Weil positive et la m\'ethode d\'ecrite permet de se mettre exactement dans la fermeture de la chambre de Weyl positive. Par r\'egularit\'e du caract\`ere infinit\'esimal, on a l'unicit\'e du quotient de Langlands dans cette situation. Cela termine la preuve du lemme.

\subsubsection{Param\`etres de Langlands des repr\'esentations cohomologiques d'apr\`es \cite{vz}\label{descriptioncoho}}
On reprend toutes les notations et les hypoth\`eses de \ref{proprietelanglands}. En particulier la repr\'esentation $\tau$ est une induite de s\'erie discr\`ete avec un caract\`ere infinit\'esimal entier r\'egulier.

\begin{rem} {\bf 1}. Les repr\'esentations ayant de la cohomologie dans le paquet associ\'e \`a $\psi$ sont toutes des quotients de Langlands des repr\'esentations $\pi'$ d\'ecrites en \ref{proprietelanglands} quand la famille des $\ell_{\delta,b}$ et $\pi_{0}$ varient.
\end{rem}
L'\'enonc\'e ne dit pas que tous les choix sont possibles car il me semble qu'il y a  un rapport entre $\pi_{0}$ et les $(\delta,b,\ell_{\delta,b})$ plus fort que celui mis sur le caract\`ere infinit\'esimal.

C'est une cons\'equence de \cite{vz} 6.16 avec comme difficult\'e que cette r\'ef\'erence ne s'applique qu'aux groupes connexes; dans le cas des groupes orthogonaux, il faut utiliser la description en termes d'induites cohomologies de \cite{aj} 2.2 et le calcul des param\`etres de Langlands de telles induites donn\'e en \cite{vogan} 6.6.15.

\

\begin{rem} {\bf 2}. On suppose qu'il existe une repr\'esentation temp\'er\'ee $\tau'$ d'un GL convenable et un demi-entier $s'$ et une repr\'esentation irr\'eductible $\pi^-$ d'un groupe de m\^eme type que $G$ tel que $\pi$ soit un quotient de $\tau\vert\,\vert^{s'}\times \pi^-$. Alors $\pi^-$ a de la cohomologie.
\end{rem}
On compare la param\'etrisation de Langlands de $\pi^-$ avec celle de $\pi$ et cela montre que n\'ecessairement, il existe un sous-ensemble ${\mathcal E}$ de $Jord(\psi)$ tel que pour tout $(\delta,b)\in {\mathcal E}$ on ait $(b-1)/2=s'$, $\tau\simeq_{(\delta,b)\in {\mathcal E}}\delta$ et avec nos notations $\ell_{\delta,b}<b$. Cela donne aussi que $\pi^-$ est le quotient de Langlands de l'analogue de $\pi'$ en rempla\c{c}ant $(\delta,b)$ par $(\delta,b-2)$ dans $Jord(\psi)$ sans changer $\ell_{\delta,b}$. En utilisant \cite{aj} paragraphe 3, on voit que $\pi^-$ est dans un des paquets d'Adams-Johnson et a donc de la cohomologie.

\

\begin{rem} {\bf 3}. On pourrait tr\`es certainement aussi d\'ecrire les couples $\tau,s_{0}$ tel que le quotient de Langlands de \ref{proprietelanglands} ait de la cohomologie et la r\'eponse est sans doute  que l'analogue de $\psi^+$ existe. 
\end{rem}
Pour pouvoir utiliser une telle remarque, il faudrait supposer que les paquets d'Adams-Johnson sont ceux d'Arthur donc mieux vaut mettre directement l'hypoth\`ese d'existence de cohomologie les rares fois o\`u cela est n\'ecessaire.

\subsubsection{Holomorphie des op\'erateurs d'entrelacement standard\label{holomorphiestandard}}
On consid\`ere ici l'op\'erateur d'entrelacement standard d\'ej\`a introduit $M(s,\tau,\pi)$ en $s=s_{0}>0$. Avec \ref{descriptioncoho} on fixe $\pi'$ tel que $\pi$ soit le quotient de Langlands de $\pi'$ et on consid\`ere l'op\'erateur d'entrelacement standard $M(s,\tau,\pi')$. Il suffit de d\'emontrer que $M(s,\tau,\pi')$ est holomorphe en $s=s_{0}$ pour avoir le m\^eme r\'esultat pour $M(s,\tau,\pi)$. On d\'emontre donc que $M(s,\tau,\pi')$ est holomorphe en $s=s_{0}$ avec les hypoth\`eses de r\'egularit\'e faites sur le caract\`ere infinit\'esimal de $\tau\vert\,\vert^{s_{0}}\times \pi'$. On d\'ecompose $M(s,\tau,\pi')$ en le produit des 3 op\'erateur d'entrelacement standard, le premier et le dernier \'etant dans des groupes de type GL et celui du milieu dans un groupe de m\^eme type que $G$:
$$
\tau\vert\,\vert^{s}\times_{(\delta,b)}Speh(\delta,(b-1)/2,(\ell_{\delta,b}+1)/2) \rightarrow
\times_{(\delta,b)}Speh(\delta,(b-1)/2,(\ell_{\delta,b}+1)/2)\times \tau\vert\,\vert^{s};\eqno(1)
$$
$$
\tau\vert\,\vert^s\times \tau' \rightarrow \tau^*\vert\,\vert^{-s}\times \tau';\eqno(2)
$$
$$
\times_{(\delta,b)}Speh(\delta,(b-1)/2,(\ell_{\delta,b}+1)/2)\times \tau^*\vert\,\vert^{-s}\rightarrow
\tau^*\vert\,\vert^{-s}\times_{(\delta,b)}Speh(\delta,(b-1)/2,(\ell_{\delta,b}+1)/2).\eqno(3)
$$
L'op\'erateur d'entrelacement (2) est holomorphe en $s=s_{0}$ gr\^ace aux travaux d'Harish-Chandra car on suppose $s_{0}>0$; on rappelle que $\tau'$ est une s\'erie discr\`ete et que $\tau$ est temp\'er\'ee.

On a d\'emontr\'e l'holomorphie des op\'erateurs (1) et (3) en \cite{mw}I.8 mais \`a condition de les normaliser \`a la Langlands-Shahidi. Il suffit donc de montre que les facteurs de normalisations n'ont ni z\'ero ni p\^ole.  Les facteurs de normalisation sont les m\^emes pour (1) et (3) (peut-\^etre aux facteurs $\epsilon$ pr\`es que l'on ne consid\`ere par ici) et ont m\^eme comportement que
$$
\prod_{(\delta,b)}\frac{L(\tau\times \delta,s-(b-1)/2))}{L(\tau\times \delta,s+(\ell_{\delta,b}+1)/2+1))}.
$$
Ce qui compte est donc le num\'erateur; on \'ecrit explicitement $L(\tau\times \delta,s-(b-1)/2))$ gr\^ace \`a \cite{shahididuke} 1.4 repris en \cite{kasten} 1.5. Avec la formule de produit, on traite le cas o\`u $\tau$ est une s\'erie discr\`ete ou un caract\`ere unitaire et on va m\^eme supposer que $\tau$ et $\delta$ sont des s\'eries discr\`etes (le cas des caract\`eres s'obtient en faisant $a_{\tau}=1$ ou $a_{\delta}=1$ ci-dessous et en ne gardant qu'un seul des 2 facteurs qui deviennent \'egaux)
$$L(\tau\times \delta,s-(b-1)/2))=
\Gamma(s-(b-1)/2+\vert a_{\tau}-a_{\delta}\vert/2)\Gamma(s-(b-1)/2+(a_{\tau}+a_{\delta})/2-1).
$$
On ne peut pas avoir $s_{0}+(a_{\tau}-1)/2\leq (b-1)/2-(a_{\delta}-1)/2$ car par r\'egularit\'e, on aurait aussi $s_{0}+(a_{\tau}-1)/2<-(b-1)/2-(a_{\delta}-1)/2$ ce qui contredit la positivit\'e de $s_{0}$. Ainsi le 2e facteur n'a pas de p\^ole. Etudions le 1e facteur en supposant d'abord que $a_{\tau}\geq a_{\delta}$; dans ce cas, si $s_{0}+(a_{\tau}-1)/2\leq (b-1)/2+(a_{\delta}-1)/2$ par r\'egularit\'e on a aussi 
$s_{0}+(a_{\tau}-1)/2\leq -(b-1)/2+(a_{\delta}-1)/2$ et $$(a_{\tau}-1)/2<s_{0}+(a_{\tau}-1)/2\leq -(b-1)/2+(a_{\delta}-1)/2\leq (a_{\delta}-1)/2$$
d'o\`u $a_{\tau}<a_{\delta}$ ce qui est contradictoire. D'o\`u encore l'absence de p\^ole si $a_{\tau}\geq a_{\delta}$. Dans le cas oppos\'e, l'in\'egalit\'e $s_{0}-(a_{\tau}-1)/2\leq -(a_{\delta}-1)/2+(b-1)/2$, force, par r\'egularit\'e, l'in\'egalit\'e $s_{0}-(a_{\tau}-1)/2\leq -(a_{\delta}-1)/2-(b-1)/2$ et on trouve encore $a_{\tau}>a_{\delta}$ ce qui est contradictoire. On a donc d\'emontr\'e l'absence de p\^ole dans tous les cas. Cela montre l'holomorphie de l'op\'erateur d'entrelacement standard associ\'e \`a (1). D'o\`u l'holomorphie cherch\'ee qui termine la preuve.
\subsubsection{Conclusion de la preuve}
On d\'emontre maintenant la proposition: on a d\'emontr\'e l'unicit\'e du quotient irr\'eductible et son identification au quotient de Langlands en \ref{proprietelanglands}; on a d\'emontr\'e l'holomorphie de l'op\'erateur d'entrelacement standard en \ref{holomorphiestandard}. On sait maintenant que l'image de l'op\'erateur d'entrelacement en $s=s_{0}$, \'etant non nulle a pour quotient irr\'eductible le quotient de Langlans. Par la m\^eme m\'ethode que celle donnant l'identification du quotient irr\'eductible de $\tau\vert\,\vert^{s_{0}}\times \pi$ avec le quotient de Langlands, on montre que l'induite $\tau^*\vert\,\vert^{-s_{0}}\times \pi$ admet cette repr\'esentation irr\'eductible comme unique sous-module irr\'eductible. La multiplicit\'e 1 force alors l'image de l'op\'erateur d'entrelacement \`a \^etre r\'eduit \`a cette repr\'esentation. Cela termine la preuve de la proposition.

\subsection{Op\'erateurs d'entrelacement normalis\'es\label{casarchimedien}}
Fixons $\psi$ param\'etrant un paquet de repr\'esentations unitaires ayant de la cohomologie de $G({\mathbb R})$. On se doute que le paquet de repr\'esentations associ\'e par Arthur \`a $\psi$ est le m\^eme que celui de \cite{aj}, mais ceci n'est pas \'ecrit et on n'en a pas besoin, comme me l'a fait remarquer H. Grobner, la seule chose qui compte pour nous est que le caract\`ere infinit\'esimal des repr\'esentations dans le paquet d'Arthur est le caract\`ere infinit\'esimal d\'ecrit en \ref{caractereinfinitesimal}. Ceci est connu. On fixe  $\pi$ dans le paquet d'Arthur associ\'e \`a $\psi$.

Soit $\tau$ une repr\'esentation temp\'er\'ee; on d\'efinit $r(s,\tau,\psi)$ comme dans le cas $p$-adique. Cela permet de poser: $N_{\psi}(s,\tau,\pi)=r(s,\tau,\psi)^{-1}M(\tau, \pi,s)$.
\begin{prop} Soit $s_{0}\in {\mathbb R}_{>0}$; on suppose que l'induite $\tau\vert\,\vert^{s_{0}}\times \pi$ a un caract\`ere infinit\'esimal r\'egulier. Alors $N_{\psi}(s,\tau,\pi)$ est holomorphe en $s=s_{0}$ et son image est la repr\'esentation irr\'eductible isomorphe \`a l'unique (cf. \ref{descriptionstandard}) quotient irr\'eductible de l'induite de $\tau\vert\,\vert^{s_{0}}\times \pi$.
\end{prop}
En tenant compte de \ref{descriptionstandard}, il suffit (il faut d'ailleurs aussi) montrer que $ordre_{s=s_{0}}r(s,\psi,\tau)=0$. Par la formule de produit, on se ram\`ene au cas o\`u $\tau$ est une s\'erie discr\`ete et on d\'ecompose $\psi$ comme en \ref{caractereinfinitesimal}. Comme dans le cas $p$-adique le facteur faisant intervenir $r_{G}$ ne donne ni z\'ero ni p\^ole car il ne d\'epend que de $\tau$ et est calcul\'e en $s_{0}>0$. Les autres termes sont: (les notations sont celles de \ref{caractereinfinitesimal})
$$
\times_{(\rho,b)}L(\tau\times \delta,s-(b-1)/2)/L(\tau\times \delta,s+(b+1)/2).
$$
Il faut donc v\'erifier que $ordre_{s=s_{0}}L(\tau\times \delta,s-(b-1)/2)=0$; mais cette d\'emonstration a \'et\'e faite en \ref{holomorphiestandard}. Cela termine la preuve.

\section{Applications aux s\'eries d'Eisenstein\label{applicationeisenstein}}
Ici on fixe un corps de nombres, $k$. On fixe aussi une repr\'esentation automorphe de carr\'e int\'egrable $\pi$ de $G$. On note $\pi^{GL}$ la repr\'esentation $\theta$-discr\`ete du groupe lin\'eaire convenable associ\'e par Arthur \`a $\pi$; on \'ecrit $\pi^{GL}\simeq \times_{(\rho',b')\in Jord(\pi^{GL})}Speh(\rho',b')$; ceci d\'efinit $Jord(\pi^{GL})$. On note encore $r_{G}$ la repr\'esentation $Sym^{2}$ (resp. $\wedge^2$) de tout groupe lin\'eaire si $G$ est un groupe symplectique ou orthogonal pair (resp. un groupe orthogonal impair). Soit $\rho$ une repr\'esentation cuspidale unitaire d'un groupe $GL$. On a l'analogue global de $r(s,\tau,\psi)$, not\'e $r(s,\rho,\pi^{GL})$ avec:
$$
r(s,\rho,\pi^{GL}):=\frac{L(\rho\times \pi^{GL},s)}{L(\rho\times \pi^{GL},s+1)}\frac{L(\rho,r_{G},2s)}{L(\rho,r_{G},2s+1)}.
$$
Comme toutes les fonctions $L$ intervenant sont des fonctions pour des groupes $GL$, on conna\^{\i}t assez leurs p\^oles en tout point $s=s_{0}\in {\mathbb R}_{\geq 1/2}$. En un tel point, on a 
$
ordre_{s=s_{0}}r(s,\rho,\pi^{GL})\geq -1$. Et cet ordre vaut $-1$ exactement quand $s_{0}$ est un demi-entier et les 2 conditions ci-dessous sont satisfaites:

(1) soit $s_{0}=1/2$ et $L(\rho,r_{G},s)$ a un p\^ole en $s=1$ soit $s_{0}\geq 1$ et $\pi^{GL}$ \'ecrit sous-forme d'induite de repr\'esentations de Speh a dans cet induction le facteur $Speh(\rho,2s_{0}-1)$;

(2) pour tout $\rho',b'$ tel que $\pi^{GL}$ \'ecrit sous-forme d'induite de repr\'esentations de Speh a dans cet induction le facteur $Speh(\rho,b')$, avec $b'=2s_{0}$, $L(\rho\times \rho',1/2)\neq 0$.

\

Pour $s\in {\mathbb R}$, on consid\`ere les s\'eries d'Eisenstein $E(\rho\times \pi,s)$. On fixe aussi $s_{0}\in {\mathbb R}_{>0}$ et on suppose que le caract\`ere infinit\'esimal de l'induite $\rho\vert\,\vert^{s_{0}}\times \pi$ est entier r\'egulier.  On a le th\'eor\`eme suivant (essentiellement d\'emontr\'e en \cite{manuscripta}):
\begin{thm} En $s=s_{0}$, les s\'eries d'Eisenstein $E(\rho\times \pi,s)$ ne sont  pas holomorphes seulement si $ordre_{s=s_{0}}r(s,\rho,\pi^{GL})=-1$. De plus si l'une de ces conditions est satisfaite, le p\^ole en $s=s_{0}$ est d'ordre au plus 1.
\end{thm}
En effet en \cite{manuscripta} 3.6, on a montr\'e que l'holomorphie des op\'erateurs d'entrelacement locaux, en les places p-adiques, normalis\'es, not\'e $N_{\psi}(s,\cdots)$ en $s=s_{0}$ entra\^{\i}ne l'holomorphie de tous les op\'erateurs d'entrelacement normalis\'es qui calculent les termes constants des s\'eries d'Eisenstein consid\'er\'ees. Gr\^ace \`a \ref{generalisation}, on a donc cette propri\'et\'e d'holomorphie pour tous les op\'erateurs d'entrelacement \`a consid\'erer.  Pour les places archim\'ediennes, l'holomorphie est facile car on a montr\'e comment remplacer la repr\'esentation $\rho_{v}\times \pi_{v}$ par une induite dans la chambre de Weyl positive (\ref{proprietelanglands}); a priori on n'est que dans la fermeture de la chambre de Weyl positive, et il faut donc normaliser les op\'erateurs d'entrelacement pour avoir l'holomorphie; mais la r\'egularit\'e du caract\`ere infinit\'esimal permet de montrer comme en \ref{casarchimedien}  que ces facteurs de normalisations n'ont ni z\'ero ni p\^ole.

Les p\^oles des s\'eries d'Eisenstein ne viennent donc que des facteurs de normalisations de \cite{manuscripta} 2.3. On a calcul\'e ces facteurs de normalisations en loc. cite, ce qui est tr\`es facile, et avec l'hypoth\`ese de r\'egularit\'e sur le caract\`ere infinit\'esimal, le p\^ole est au plus d'ordre 1. De plus si un tel facteur de normalisation a un p\^ole alors il en est de m\^eme de la fonction $r(s,\rho,\pi)$ qui est le facteur de normalisation pour l'entrelacement $\rho\vert\,\vert^{s}\times \pi\rightarrow \rho\vert\,\vert {-s}\times \pi$. D'o\`u le th\'eor\`eme.

\

Dans tout ce qui suit, on suppose que les hypoth\`eses du th\'eor\`eme sont satisfaites: $\pi$ a de la cohomologie \`a l'infini et l'induite $\rho\vert\,\vert^{s_{0}}\times \pi$ a un caract\`ere infinit\'esimal entier r\'egulier et on suppose que $ordre_{s=s_{0}}r(s,\rho,\pi^{GL})=-1$. On d\'ecrit d'abord a priori la repr\'esentation qui va d\'eterminer l'existence d'un r\'esidu non nul; pour le faire on n'a besoin que de la condition (1) ci-dessus, pas la condition (2).

Par d\'efinition, la repr\'esentation $\pi^+$ est l'image de l'op\'erateur d'entrelacement normalis\'e par la fonction $r(s,\rho,\tau)$:
$$
\rho\vert\,\vert^{s}\times \pi\rightarrow \rho\vert\,\vert {-s}\times \pi
$$
calcul\'e en $s=s_{0}$. Plus pr\'ecis\'ement:

soit $v$ une place finie, on a d\'efini $\pi^+_{v}$ en partant de $\pi_{v}$ en \ref{parametregeneralisation}; soit $v$ une place archim\'edienne, on note $\pi^+_{v}$ l'unique quotient irr\'eductible de l'induite $\rho_{v}\vert\,\vert^{s_{0}}\times \pi_{v}$. Et on pose $\pi^+:=\otimes_{v}\pi_{v}$; dans certains cas, $\pi^+=0$ et si cette repr\'esentation n'est pas nulle on ne sait pas que cette repr\'esentation est automorphe de carr\'e int\'egrale; on montre ici qu'elle est de carr\'e int\'egrable  quand les conditions (1) et (2) du th\'eor\`eme pr\'ec\'edent sont satisfaites car se r\'ealisant dans les r\'esidus. La r\'eciproque est certainement aussi vraie mais on renvoie \`a \ref{commentaires} pour expliquer ce qui manque pour la montrer.

\

Le but de cette section est de montrer que $s_{0}$ v\'erifiant les conditions n\'ecessaires du th\'eor\`eme est un p\^ole si et seulement si $\pi^{+}$ est non nul; dans le cas o\`u $\pi$ est cuspidal c'est exactement cela que l'on d\'emontre et dans le cas g\'en\'eral, pour le faire, on ajoute l'hypoth\`ese que, $\pi^+_{v}$ a de la cohomologie si $v$ est une place archim\'edienne et que les repr\'esentations automorphes irr\'eductibles ayant de la cohomologie interviennent avec multiplicit\'e au plus 1 dans le spectre discret des groupes classiques de m\^eme type que $G$. On d\'emontre alors aussi que l'espace des r\'esidus est irr\'eductible quand il n'est pas nul.

\subsection{Le cas de rang 1\label{rang1}}
On suppose ici que $\pi$ est une repr\'esentation cuspidale et ce cas est alors tr\`es simple. 
\begin{thm}L'espace des r\'esidus $\biggl((s-s_{0})E(\rho\times \pi,s)\biggr)_{s=s_{0}}$ est non nul exactement quand $\pi^+$ est non nul et alors la repr\'esentation ainsi d\'efinie est irr\'eductible isomorphe \`a $\pi^+$.
\end{thm}
On sait que l'espace des r\'esidus $\biggl((s-s_{0})E(\rho\times \pi,s)\biggr)_{s=s_{0}}$ a une projection nulle sur l'espace des formes automorphes cuspidales (\cite{bible} 1.3.4 pour la d\'efinition de cette projection). Et toute fonction non nulle dans cet espace de r\'esidu a exactement un terme constant cuspidal non nul, c'est celui qui est relatif au parabolique maximal de Levi, $GL(d_{\rho})\times G$. Et ce terme constant appartient \`a l'image de l'op\'erateur d'entrelacement $\otimes'_{v}N_{\psi}(\rho\times \pi, s_{0})$ (c'est un calcul facile fait en \cite{bible} II.1.7). On a vu que l'image de cet op\'erateur est  exactement $\pi^+$, c'est-\`a-dire est nul exactement si $\pi^+$ est nul et est isomorphe \`a cette repr\'esentation sinon. D'o\`u le th\'eor\`eme.

\subsection{Le cas g\'en\'eral\label{residugeneral}}
On veut le m\^eme r\'esultat que dans le cas de rang 1 mais j'ai besoin d'une hypoth\`ese de multiplicit\'e 1. La formule de multiplicit\'e annonc\'ee par Arthur, donne la multiplicit\'e de $\pi^+$ en fonction de multiplicit\'e locale; on sait que ces multiplicit\'es locales sont au plus 1 aux places p-adiques (\cite{discret}, \cite{pourshahidi}) et si l'on admet que $\pi^+_{\infty}$ a de la cohomologie et que les paquets d'Arthur sont aussi ceux d'Adams-Johnson dans ce cas, il r\'esulte aussi de \cite{aj} que la multiplicit\'e locale aux places archim\'ediennes est au plus 1. Ainsi la formule de multiplicit\'e d'Arthur donne au plus 1. On prend donc ici l'hypoth\`ese raisonable que l'on sait a priori que les repr\'esentations automorphes de carr\'e int\'egrable ayant de la cohomologie \`a l'infini interviennent avec multiplicit\'e 1 dans le spectre discret. 

On introduit les notations suivantes: soit $\pi'$ une repr\'esentation automorphe de carr\'e int\'egrable d'un groupe classique explicitement r\'ealis\'ee. Soit $\rho'$ une repr\'esentation cuspidale unitaire d'un $GL(d')$ (cela d\'efinit $d'$) et soit $s'\in {\mathbb R}$. On suppose que le groupe classique portant $\pi'$ admet un sous-groupe parabolique standard maximal dont le Levi est isomorphe \`a $GL(d')\times$ un groupe classique, sinon la notation est vide. On note $\pi'[\rho',s']$, la projection de l'espace des termes constants de $\pi'$ le long du radical unipotent du parabolique standard sur la repr\'esentation cuspidale $\rho'\vert\,\vert^{-s'}$ de $GL(d')$ (cf. \cite{bible}1.3.4). La non nullit\'e de $\pi'[\rho',s']$ impose des conditions fortes \`a $\rho',s'$ en particulier $s'$ est n\'ecessairement strictement positif et est un demi-entier (ici aussi on utilise les r\'esultats annonc\'es d'Arthur).

Dans le th\'eor\`eme ci-dessous, $s_{0}$ v\'erifie les conditions n\'ecessaires de \ref{applicationeisenstein} et si $\pi^+$ est non nulle, on suppose que $\pi^+$ a de la cohomologie. 

\begin{thm}(i) L'espace des r\'esidus $\biggl((s-s_{0})E(\rho\times \pi,s)\biggr)_{s=s_{0}}$ est non nul exactement quand $\pi^+$ est non nulle et  la repr\'esentation ainsi d\'efinie est alors irr\'eductible isomorphe \`a $\pi^+$.

(ii) Soit $\rho',s'$ un couple form\'e d'une repr\'esentation cuspidale unitaire d'un $GL(d_{\rho'})$ et d'un demi-entier strictement positif; on suppose que  $\pi^+[\rho',s']\neq 0$. Alors, il existe une repr\'esentation de carr\'e int\'egrable $\pi'$ d'un groupe classique de m\^eme type que $G$ tel que
$
\pi^+= \biggl((s'-s'_{0})E(\rho'\times \pi',s')\biggr)_{s'=s'_{0}}.
$
C'est une \'egalit\'e dans l'espace des formes automorphes du groupe convenable.
\end{thm}
On montre (i) et (ii) simultan\'ement par r\'ecurrence sur le nombre $S(\pi)$ d\'efini pour toute repr\'esentation de carr\'e int\'egrable, $\tilde{\pi}$ d'un groupe classique de m\^eme type que $G$ par la formule: $S(\tilde{\pi}):=\sum_{(\rho',b')\in Jord(\tilde{\pi})}(b'-1)$.

Initions la r\'ecurrence: on suppose que $S(\pi)=0$. On sait alors que $\pi$ est cuspidale (cf. par exemple les conditions n\'ecessaires \cite{manuscripta} 1.3). Alors (i) est d\'emontr\'e m\^eme dans sa version plus forte en \ref{rang1} et le (ii)  est clair car la non nullit\'e entra\^{\i}ne que $\rho'\simeq \rho$ et $s'=s_{0}$.

On fixe un entier $N>0$ et on suppose que (i) et (ii) sont d\'emontr\'ees pour toutes les repr\'esentations $\tilde{\pi}$ ayant de la cohomologie et telles que $S(\tilde{\pi})<N$. On suppose que  $S(\pi)=N$ et on d\'emontre (i) et (ii).

\

On suppose d'abord que $\pi^+\neq 0$; on calcule les termes constants $$\biggl((s-s_{0})E(\rho\times \pi,s)\biggr)_{s=s_{0}}[\rho,s_{0}]$$ avec une g\'en\'eralisation \'evidente de la notation d\'ej\`a introduite. Cet espace est exactement l'image  de l'op\'erateur d'entrelacement standard $(s-s_{0})M(\rho\vert\,\vert^{s}\times \pi)$ \'evalu\'e en $s=s_{0}$ appliqu\'e \`a l'induite $\rho\vert\,\vert^{s_{0}}\times \pi$. Cet op\'erateur n'est autre que l'op\'erateur d'entrelacement normalis\'e et l'image est donc $\pi^+$. Ainsi l'espace des r\'esidus est certainement non nul.

On suppose maintenant que l'espace des r\'esidus est non nul et on montre que $\pi^+\neq 0$ et l'identification de $\pi^+$ avec l'espace des r\'esidus. On sait  que l'espace des r\'esidus est form\'e de formes automorphes de carr\'e int\'egrable: en effet la r\'egularit\'e du caract\`ere infinit\'esimal de l'induite $\rho\vert\,\vert^{s_{0}}\times \pi$ assure que $(\rho,2s_{0}+1)\notin Jord(\pi^{GL})$, les hypoth\`eses de \cite{manuscripta} 1.2.2 sont satisfaites et cette r\'ef\'erence donne le r\'esultat.

On fixe une sous-repr\'esentation irr\'eductible, disons $\pi''$ (la notation va dispara\^{\i}tre rapidement), de cet espace et on va montrer qu'elle est isomorphe \`a $\pi^+$, avec la multiplicit\'e 1, cela suffira. Pour montrer cela, on va montrer que pour tout $\rho',s'_{0}$ tel que $\pi''[\rho',s'_{0}]\neq 0$, l'ensemble de ces termes constants vus comme sous-module de l'induite global de $\rho'\vert\,\vert^{-s'_{0}}\times {\mathcal A}(G(k)\backslash G({\mathbb A}))$ est une sous-repr\'esentation isomorphe \`a $\pi^+$.

D'apr\`es ce que l'on a vu ci-dessus, si $(\rho',s'_{0})=(\rho_{0},s_{0})$, le r\'esultat est clair. On suppose donc qu'il n'en est pas ainsi.
On sait alors, par le calcul facile des termes constants des s\'eries d'Eisenstein (cf. \cite{bible} II.1.7) que la projection des termes constants de $\pi$,  $\pi[\rho',s'_{0}]$, est non nulle; ainsi $\pi$ n'est pas cuspidal et on lui applique (ii): $\pi=\biggl((s'-s'_{0})E(\rho'\times \pi', s')\biggr)_{s'=s'_{0}}$ pour un bon choix de repr\'esentation de carr\'e int\'egrable $\pi'$. On change un peu les notations pour qu'elles soient plus explicites en rempla\c{c}ant $\rho$ par $\rho\vert\,\vert^s$ et $\rho'$ par $\rho'\vert\,\vert^{s'}$. D'o\`u:
$$
\biggl((s-s_{0})E(\rho\vert\,\vert^{s}\times \pi,s)\biggr)_{s=s_{0}}=\biggl((s-s_{0})\biggl((s'-s'_{0})E(\rho\vert\,\vert^s\times \rho'\vert\,\vert^{s'}\times \pi',s,s')\biggr)_{s'=s'_{0}}\biggr)_{s=s_{0}}. \eqno(1)
$$
On ne peut pas enlever les parenth\`eses dans le terme de droite et on doit distinguer suivant le signe de $s_{0}-s'_{0}$.

On suppose d'abord  que $s_{0}\geq s'_{0}$. Sous cette hypoth\`ese, la fonction m\'eromorphe de $s,s'$ valant $(s-s_{0})(s'-s'_{0})E(\rho\vert\,\vert^{s}\times \rho'\vert\,\vert^{s'}\times \pi',f)$ o\`u $f$ est une fonction  dans la bonne induite globale, est holomorphe au voisinage de $s=s_{0},s'=s'_{0}$. Pour calculer le terme de droite ci-dessus, on peut \'evaluer dans l'ordre que l'on veut. On  note $\sigma_{1,2}$ l'\'el\'ement du groupe de Weil de $GL(d+d')$ qui \'echange les facteurs $GL(d)\times GL(d')$ et de longueur minimale dans sa double classe. Pour $f$ comme ci-dessus $M(\sigma_{1,2},s-s')f$ est bien d\'efini (\cite{mw} 1.10 dit que les p\^oles de $M(\sigma_{1,2},s-s')$ sont ceux de $L(\rho\times \rho^{'*},s-s')/L(\rho\times \rho^{'*},s-s'+1)$ qui n'en a pas) et on a l'\'equation fonctionnelle, qui est une \'egalit\'e de fonctions m\'eromorphes:
$$(s-s_{0})(s'-s'_{0})
E(\rho\vert\,\vert^{s}\times \rho'\vert\,\vert^{s'}\times \pi',f)=$$
$$(s-s_{0})(s'-s'_{0})E(\rho'\vert\,\vert^{s'}\times \rho\vert\,\vert^{s}\times \pi',M(\sigma_{1,2},s-s')f). \eqno(2)
$$Le terme de gauche est une bonne fonction holomorphe au voisinage de $s=s_{0}, s'=s'_{0}$ alors que pour le terme de droite, il faut prendre plus de pr\'ecaution. En effet
soit $f'$ une fonction dans l'induite $\rho\vert\,\vert^{s'_{0}}\times\rho\vert\,\vert^{s_{0}}\times \pi'$; on construit la fonction m\'eromorphe
$
(s-s_{0})(s'-s'_{0})E(\rho'\vert\,\vert^{s'}\times \rho\vert\,\vert^{s}\times \pi',f')$. Cette fonction n'a aucune raison d'\^etre holomorphe en $s=s_{0},s'=s'_{0}$ par contre elle l'est si on la calcule d'abord sur l'hyperplan $s=s_{0}$ puis en $s'=s'_{0}$. En faisant cela on garde l'\'egalit\'e de (2). Avec cette \'egalit\'e, on calcule les termes constants relativement \`a $\rho',s'_{0}$ pour la s\'erie d'Eisenstein associ\'e \`a $f$ gr\^ace \`a l'analogue pour la s\'erie d'Eisenstein associ\'ee \`a la  fonction $f':=M(\sigma_{1,2},s_{0}-s'_{0})f$. Ainsi
 les termes constants  $$\biggl((s-s_{0})E(\rho\vert\,\vert^s\times \rho'\vert\,\vert^{s'}\times \pi',s,s')\biggr)_{(s=s_{0})}[\rho',s'_{0}]\hookrightarrow$$
 $$ \biggl(s'-s'_{0})\biggl((s-s_{0})E(\rho'\vert\,\vert^{s'}\times \rho\vert\,\vert^s\times \pi',s',s)\biggr)_{s=s_{0}}\biggr)_{s'=s'_{0}}[\rho',s'_{0}].$$On sait calculer le terme de droite c'est l'image par l'op\'erateur d'entrelacement convenable 
$$
\rho'\vert\,\vert^{s'_{0}}\times \biggl((s-s_{0})E( \rho\vert\,\vert^s\times \pi',s)\biggr)_{s=s_{0}}\rightarrow \rho'\vert\,\vert^{-s'_{0}}\times \biggl((s-s_{0})E( \rho\vert\,\vert^s\times \pi',s)\biggr)_{s=s_{0}}.
$$
Il faut identifier cette repr\'esentation. On consid\`ere le diagramme o\`u les op\'erateurs d'entrelacement sont les op\'erateurs d'entrelacement standard:
$$
\begin{array}{lll}
\rho\vert\,\vert^s\times \rho'\vert\,\vert^{s'}\times \pi' & \stackrel{M(\sigma_{1,2}, s-s')}{\rightarrow} &\rho'\vert\,\vert^{s'}\times \rho\vert\,\vert^{s}\times \pi'\\
M(w,s,s') \downarrow & &\downarrow M(w',s',s)\\
\rho\vert\,\vert^{-s}\times \rho'\vert\,\vert^{-s'}\times \pi' &\stackrel{M(\sigma_{2,1},s-s')}{\leftarrow}&\rho'\vert\,\vert^{-s'}\times\rho\vert\,\vert^{-s}\times \pi'.
\end{array}
$$
On note $N(w,s,s'):=(s-s_{0})(s'-s'_{0})M(w,s,s')$ et $N(w',s',s):= (s-s_{0})(s'-s'_{0})M(w',s,s')$. On a alors encore un diagramme commutatif
$$
\begin{array}{lll}
\rho\vert\,\vert^s\times \rho'\vert\,\vert^{s'}\times \pi' & \stackrel{M(\sigma_{1,2}, s-s')}{\rightarrow} &\rho'\vert\,\vert^{s'}\times \rho\vert\,\vert^{s}\times \pi'\\
N(w,s,s') \downarrow & &\downarrow N'(w',s',s)\\
\rho\vert\,\vert^{-s}\times \rho'\vert\,\vert^{-s'}\times \pi' &\stackrel{M(\sigma_{2,1},s-s')}{\leftarrow}&\rho'\vert\,\vert^{-s'}\times\rho\vert\,\vert^{-s}\times \pi'.
\end{array}
$$
Tous les op\'erateurs consid\'er\'es sont holomorphes en $s=s_{0}, s'=s'_{0}$ sauf pr\'ecis\'ement $N(w',s',s)$; cet op\'erateur doit \^etre calcul\'e d'abord sur l'hyperplan $s=s_{0}$ puis en $s'=s'_{0}$. On note $N'(w,s'_{0},s_0)$ le r\'esultat de ce calcul.  Pour les op\'erateurs $M(\sigma_{1,2},s-s')$ et $M(\sigma_{2,1},s-s')$, on utilise  \cite{mw} I. 10, comme ci-dessus. Le diagramme  est  commutatif.

Par hypoth\`ese, on sait que $N'(w,s'_{0},s_{0})\circ M(\sigma_{1,2},s_{0}-s'_{0})$ est non identiquement 0; cet op\'erateur a alors pour image la repr\'esentation irr\'eductible image de $N'(w,s'_{0},s_{0})$; on note $\pi^{'+}$ cette repr\'esentation.

Les 2 repr\'esentations sur une m\^eme colonne du diagramme sont en dualit\'e naturelle; et dans ces dualit\'e $M(\sigma_{2,1},s-s')$ est l'application duale de $M(\sigma_{1,2},s-s')$. C'est-\`a-dire, que pour $x$ dans l'induite en haut \`a gauche et $y^*$ dans l'induite en bas \`a droite, on a:
$$
<M(\sigma_{1,2},s_{0}-s'_{0})x,y^*>=<x, M(\sigma_{2,1},s_{0}-s'_{0})y^*>.
$$
Soit $y,y'$ dans l'induite en haut \`a droite (au point $s=s_{0},s'=s'_{0}$);  $$<N'(w',s'_{0},s_{0})y,y'>=0 \Leftarrow \hbox{   si $y$ ou $y' \in Ker\,N'(w',s'_{0},s_{0})$.}$$
 Ainsi l'application $(y,y')\mapsto  <N'(w',s'_{0},s_{0})y,y'>$ induit une  forme bilin\'eaire  non d\'eg\'en\'er\'ee sur  la repr\'esentation irr\'eductible, $\pi^{'+}\simeq \bigl(\rho'\vert\,\vert^{s'_{0}}\times \rho\vert\,\vert^{s_{0}}\times \pi'\bigr)/Ker\, N'(',s'_{0},s_{0})$. Soit  $x\in \rho\vert\,\vert^{s_{0}}\times \rho'\vert\,\vert^{s'_{0}}\times \pi'$ tel que $N'(w,s'_{0},s_{0})\circ M(\sigma_{1,2},s_{0}-s'_{0})(x)\neq 0$. Comme l'image de $M(\sigma_{1,2},s_{0}-s'_{0})$ n'est pas incluse dans le noyau de $N'(w',s'_{0},s_{0})$, 	avec les propri\'et\'es que l'on vient de voir, il existe $x' \in \rho\vert\,\vert^{s_{0}}\times \rho'\vert\,\vert^{s'_{0}}\times \pi'$  tel que:
$$
<M(\sigma_{1,2},s_{0}-s'_{0})x',N'(w,s'_{0},s_{0})\circ M(\sigma_{1,2},s_{0}-s'_{0})(x)>\neq 0.
$$
Cela vaut aussi
$$
<x,M(\sigma_{2,1},s_{0}-s'_{0})N'(w,s'_{0},s_{0})\circ M(\sigma_{1,2},s_{0}-s'_{0})(x)>=<x,N(w,s_{0},s'_{0})x>\neq 0.
$$
Ainsi l'op\'erateur $N(w,s_{0},s'_{0})$ est non nul et son image $\pi^{+}$ est aussi non nulle. Au passage, on a montr\'e que $N'(w',s'_{0},s_{0})$ calcul\'e comme expliqu\'e est non nul et son image qui est une repr\'esentation irr\'eductible est isomorphe \`a $\pi^+$ par la commutativit\'e du diagramme. On a la suite d'application
$$
\biggl((s-s_{0})E(\rho\times \pi,s))\biggr)_{s=s_{0}}[\rho',s'_{0}]=\biggl((s-s_{0})(s'-s'_{0})E(\rho\times \rho'\times \pi',s,s')\biggr)_{s'=s'_{0},s=s_{0}}[\rho',s'_{0}]$$
$$
\hookrightarrow \rho'\vert\,\vert^{-s'_{0}}\times \biggl((s-s_{0})E(\rho\times \pi',s)_{s=s_{0}}\biggr) \rightarrow \rho'\vert\,\vert^{-s'_{0}}\times \biggl((s-s_{0})E(\rho\times \pi',s)_{s=s_{0}}\biggr) [\rho,s_{0}].
$$
La derni\`ere fl\`eche est non nulle et elle est alors n\'ecessairement injective car par l'hypoth\`ese de r\'ecurrence $\times \bigl((s-s_{0})E(\rho\times \pi',s)\bigr)_{s=s_{0}} $ est irr\'eductible. Le compos\'e de ces applications repr\'esente la projection des termes constant de $ \bigl((s-s_{0})E(\rho\times \pi,s)\bigr)_{s=s_{0}}$  sur $\rho'\vert\,\vert^{-s'_{0}}\otimes \rho\vert\,\vert^{-s_{0}}$; cet espace n'est autre que l'image de $$Im\, N(w',s'_{0},s_{0})\circ M(\sigma_{1,2},s_{0}-s'_{0}) \rho\vert\,\vert^{s_{0}}\times \rho'\vert\,\vert^{s'_{0}}\times \pi'.$$
On vient de voir que cette image est isomorphe \`a $\pi^+$ d'o\`u l'isomorphisme cherch\'ee pour terminer la preuve de (i) dans le cas o\`u $s_{0}\geq s'_{0}$. On prouve aussi (ii) dans ce cas: l'\'egalit\'e (2) (qui est  une \'equation fonctionnelle) donne
$$
\biggl((s-s_{0})(s'-s'_{0})E(\rho\times \rho'\times \pi',s,s',f)\biggr)_{s=s_{0},s'=s'_{0}}=$$
$$
\biggl((s-s_{0})\bigg((s'-s'_{0})E(\rho'\times \rho\times \pi',s,s',M(\sigma_{1,2},s_{0}-s'_{0})f)
\biggr)_{s'=s'_{0}}
\biggr)_{s=s_{0}}.
$$
Ceci inverse $\rho$ et $\rho'$ et prouve (ii) dans ce cas.

\

Il faut maintenant consid\'erer le cas o\`u $s_{0}\leq s'_{0}$; le cas $s_{0}=s'_{0}$ a d\'ej\`a \'et\'e trait\'e mais on le retrouve aussi ici. On consid\`ere le compos\'e des entrelacements: 
$$N(w,s,s'):
\rho\vert\,\vert^{s}\times \rho'\vert\,\vert^{s'}\times \pi' \rightarrow \rho\vert\,\vert^{-s}\times \rho'\vert\,\vert^{-s'}\times \pi'$$ (c'est la m\^eme notation que ci dessus) avec $M'(\sigma_{1,2},s'-s)$ qui correspond \`a l'\'echange de $\rho\vert\,\vert^{-s}\times \rho'\vert\,\vert^{-s'}\times \pi'\rightarrow \rho'\vert\,\vert^{-s'}\times \rho\vert\,\vert^{-s}\times \pi'$. Ici  $M'(\sigma_{1,2},s'-s)$  est holomorphe au voisinage de $s=s_{0}$ et $s'=s'_{0}$ par l'hypoth\`ese $s_{0}\leq s'_{0}$; $N(w,s,s')$ doit \^etre calcul\'e d'abord sur $s'=s'_{0}$ puis sur $s=s_{0}$ pour avoir l'holomorphie en $(s,s')=(s_{0},s'_{0})$. On revient \`a (1) ci-dessus; on calcule les termes constants du membre de gauche relativement \`a $\rho',s'_{0}$ et on trouve un sous-module de l'induite $\rho'\vert\,\vert^{-s'_{0}}\times \bigl((s-s_{0})E(\rho\times \pi',s)\bigr)_{s=s_{0}}$. Cet espace est donc non nul, a fortiori. On applique (i) par r\'ecurrence \`a $\pi'$ qui donne l'irr\'eductibilit\'e de l'espace $\bigl((s-s_{0})E(\rho\times \pi',s)\bigr)_{s=s_{0}}$. Cette irr\'eductibilit\'e assure que l'application qui a un r\'esidu associe sont terme constant relatif \`a $\rho\vert\,\vert^{-s_{0}}$ est injective et que l'espace ainsi engendr\'e est l'image de l'op\'erateur d'entrelacement normalis\'e \'evident. On calcule alors les termes constant du 2e membre de (1) relativement \`a $\rho'\vert\,\vert^{-s'_{0}}\otimes \rho\vert\,\vert^{-s_{0}}$ et on obtient un espace non nul. Ces termes constant sont l'image de $ M'(\sigma_{1,2},s'_{0}-s_{0})\circ N(w,s_{0},s'_{0})$. Donc en particulier $N(w,s_{0},s'_{0})$ est non nul; son image est $\pi^+$ ce qui donne la non nullit\'e cherch\'ee de $\pi^+$. De plus les termes constants cherch\'es, sont l'image de $M'(\sigma_{1,2},s'_{0}-s_{0})\circ N(w,s_{0},s'_{0})$ et forment une repr\'esentation isomorphe \`a $\pi^+$. Cela termine la preuve de (i). Montrons encore (ii) dans ce cas.

On a  ici le diagramme commutatif
$$
\begin{array}{lll}
\rho\vert\,\vert^s\times \rho'\vert\,\vert^{s'}\times \pi' & \stackrel{M'(\sigma_{1,2}, s'-s)}{\leftarrow} &\rho'\vert\,\vert^{s'}\times \rho\vert\,\vert^{s}\times \pi'\\
N(w,s,s') \downarrow & &\downarrow N(w',s',s)\\
\rho\vert\,\vert^{-s}\times \rho'\vert\,\vert^{-s'}\times \pi' &\stackrel{M'(\sigma_{2,1},s'-s)}{\rightarrow}&\rho'\vert\,\vert^{-s'}\times\rho\vert\,\vert^{-s}\times \pi'.
\end{array}
$$Toutes les fl\`eches sont holomorphe en $s=s_{0},s'=s'_{0}$ sauf $N(w,s,s')$ qu'il faut calculer d'abord sur $s'=s'_{0}$ puis sur $s=s_{0}$. Le diagramme est toujours commutatif apr\`es cette \'evaluation.
On sait que $M'(\sigma_{2,1},s'_{0}-s_{0})\circ N(w,s_{0},s'_{0})\not\equiv 0$. Ceci est encore \'equivalent \`a ce que  $\pi^+$ ne soit pas inclus dans le noyau de $M'(\sigma_{2,1},s'_{0}-s_{0})$. Par dualit\'e comme on l'a vu pr\'ec\'edemment, ceci montre que l'image de $M'(\sigma_{1,2},s'_{0}-s_{0})$ n'est pas incluse dans le noyau de $N(w,s'_{0},s'_{0})$. D'o\`u aussi
$$
N(w',s'_{0},s_{0})=M'(\sigma_{2,1},s'_{0}-s_{0})\circ N(w,s_{0},s'_{0})\circ M'(\sigma_{1,2},s'_{0}-s_{0})\not\equiv 0.
$$
Comme l'image de $N(w',s'_{0},s_{0})$ est irr\'eductible quand elle n'est pas nulle, cette image est n\'ec\'essairement isomorphe \`a $\pi^+$. Ainsi $\pi^+$ se r\'ealise dans les r\'esidus:
$$
\biggl((s'-s'_{0})(s-s_{0}) E(\rho'\times \rho\times \pi',s',s)_{{s=s_{0},s'=s'_{0}}}\biggr).
$$
Le terme entre parenth\`ese peut se calculer d'abord sur $s=s_{0}$ puis sur $s'=s'_{0}$ par holomorphie. On utilise ici encore le fait que $\pi^+$ est \'egal \`a l'espace des r\'esidus $\bigl((s-s_{0})E(\rho\times \pi,s)_{s=s_{0}}\bigr)$ par l'irriductibilit\'e de (i) que l'on vient de prouver. On obtient donc (ii). Cela termine la preuve.

\subsection{Remarque sur les p\^oles des s\'eries d'Eisenstein}
Dans cet remarque on donne un exemple d'une s\'erie d'Eisenstein form\'ee avec un parabolique maximal et une repr\'esentation non de carr\'e int\'egrable du Levi de ce parabolique ayant un p\^ole d'ordre 1 et dont l'espace des r\'esidus est de carr\'e int\'egrable.

\

Soit $M$ un sous-groupe de Levi de $G$ de la forme $GL(d)\times GL(d')\times G'$ pour $G'$ un sous-groupe de $G$ de m\^eme type que $G$; on fixe $\rho, \rho'$ des repr\'esentations cuspidales irr\'eductibles et autoduales de $GL(d)$ et $GL(d')$ respectivement et $\pi'$ une repr\'esentation de carr\'e int\'egrable irr\'eductible de $G'$; on d\'efinit donc, en suivant Arthur, la repr\'esentation $\pi^{'GL}$, $\theta$-discr\`ete du groupe lin\'eaire qui param\'etrise le paquet de repr\'esentations de carr\'e int\'egrable contenant $\pi'$.

Soit $s_{0},s'_{0}$ des demi-entiers positifs strictement et on suppose que $s_{0}= s'_{0}+1/2$. 
On suppose  que les quotients de Langlands des repr\'esentations $\rho_{v}\vert\,\vert^{s_{0}}\times \rho'_{v}\vert\,\vert^{s'_{0}}\times \pi'_{v}$ aux places archim\'ediennes $v$ ont de la cohomologie pour un bon syst\`eme de coefficients.

On suppose aussi que  $ordre_{s=s_{0}}r(s,\rho,\pi^{'GL})=-1$. On note $t:=ordre_{s=1/2}L(\rho\times \rho',s)$ et on suppose que $ordre_{s'=s'_{0}}r(s,\rho',\pi^{'GL})=t-1$; cela veut dire que $Jord(\pi^{'GL})$ contient $(\rho',2s'_{0}-1)$ et $(\rho,2s_{0}-1)=(\rho,2s'_{0})$ et que pour tout autre \'el\'ement de la forme $(\rho'',b'')$, si $b''=2s'_{0}$, $L(\rho\times \rho'',1/2)\neq 0$ et si $b''=2s_{0}$, $L(\rho\times \rho'',1/2)\neq 0$.

En termes plus explicites, \'ecrivons $\pi^{'GL}=\times_{(\rho'',b'')\in Jord(\pi^{'GL})}Speh(\rho'',b'')$; le fait que $ordre_{s=s_{0}}r(s,\rho,\pi^{'GL})=-1$ assure que $Jord(\pi^{'GL})$ contient $(\rho,2s_{0}-1)$. On rappelle que par hypoth\`ese $2s_{0}-1=2s'_{0}$. Ainsi, la fonction $r(s',\rho',\pi{'GL})$ contient le facteur $$
\frac{L(\rho\times \rho',s'-s_{0}+1)}{L(\rho\times \rho',s'+s_{0})}=\frac{L(\rho\times \rho',s'-s'_{0}+1/2)}{L(\rho\times \rho',s'+s'_{0}+1/2)}$$ et  les autres termes de cette fonction ne fournissent aucun z\'ero mais fournissent un p\^ole.

\begin{rem}{\bf 1} On suppose que $\pi'$ est cuspidale et  que $L(\rho\times \rho',1/2)=0$, c'est-\`a-dire que $t\geq 1$. Alors les s\'eries d'Eisenstein $E(\rho'\times \pi',s')$ sont holomorphes en $s'=s'_{0}$ et on note $\tilde{\pi}'$ la repr\'esentation automorphe qu'elles engendrent. Les s\'eries d'Eisenstein $E(\rho\times \tilde{\pi}',s)$ sont holomorphes en $s=s_{0}$ si $t>1$ et si $t=1$ elles peuvent avoir un p\^ole en ce point. Dans ce cas l'espace des r\'esidus $
\bigl((s-s_{0})E(\rho\times \tilde{\pi}',s)_{s=s_{0}}\bigr)$ est inclus dans l'ensemble des formes automorphes de carr\'e int\'egrable.
\end{rem}
La premi\`ere assertion concernant l'holomorphie des s\'eries d'Eisenstein est un cas particulier du th\'eor\`eme \ref{applicationeisenstein} puisque $r(s',\rho',\pi^{'GL})$ n'a pas de p\^ole en $s'=s'_{0}$. On \'etudie les s\'eries d'Eisenstein d\'ependant des variables $s,s'$ au voisinage de $s_{0},s'_{0}$:
$
E(\rho\times \rho'\times\pi',s,s')$.

On calcule tous les termes constants; ceci sont donn\'es par des op\'erateurs d'entrelacement partant de l'induite $\rho\vert\,\vert^{s}\times \rho'\vert\,\vert^{s'}\times \pi'$. Il y en a 8:

l'identit\'e qui est clairement holomorphe;

$\rho\vert\,\vert^{s}\times \rho'\vert\,\vert^{s'}\times \pi' \rightarrow \rho'\vert\,\vert^{s'}\times \rho\vert\,\vert^{s}\times \pi'$ holomorphe car $s>s'$ dans le voisinage consid\'er\'e et n\'ecessairement $\rho\not\simeq \rho'$ par l'assertion d'existence de cohomologie; on remarque pour la suite que le facteur global fait appara\^{\i}tre $L(\rho\times \rho',s-s')$ qui a un z\'ero en $s-s'=1/2$;

$\rho\vert\,\vert^{s}\times \rho'\vert\,\vert^{s'}\times \pi \rightarrow \rho\vert\,\vert^{s}\times \rho'\vert\,\vert^{-s'}\times \pi$ qui est holomorphe en $s=s_{0}$ et $s'=s'_{0}$ car la fonction $r(s',\rho',\pi^{'GL})$ n'a pas de p\^ole en $s'=s'_{0}$;

en composant avec un op\'erateur n'ayant clairement pas de p\^oles on a aussi l'holomorphie en $s=s_{0},s'=s'_{0}$ de l'op\'erateur d'entrelacement $\rho\vert\,\vert^{s}\times \rho'\vert\,\vert^{s'}\times \pi \rightarrow \rho'\vert\,\vert^{-s'}\times \rho\vert\,\vert^{s}\times \pi$;

$\rho\vert\,\vert^{s}\times \rho'\vert\,\vert^{s'}\times \pi \rightarrow \rho'\vert\,\vert^{s'}\times \rho\vert\,\vert^{-s}\times \pi$, a un p\^ole sur $(s-s_{0})$ (en g\'en\'eral) d'ordre 1qui vient de la fonction $r(s,\rho,\pi^{'GL})$ en $s=s_{0}$ mais on a aussi un facteur $L(\rho\times \rho',s-s')$ qui vient de l'\'echange des 2 premier facteur; on a donc holomorphie sur  l'hyperplan $s'=s'_{0}$;

en composant avec un op\'erateur n'ayant clairement pas de p\^oles on a aussi l'holomorphie sur l'hyperplan $s'=s'_{0}$ de l'op\'erateur $\rho\vert\,\vert^{s}\times \rho'\vert\,\vert^{s'}\times \pi \rightarrow \rho\vert\,\vert^{-s}\times \rho'\vert\,\vert^{s'}\times \pi$;

$\rho\vert\,\vert^{s}\times \rho'\vert\,\vert^{s'}\times \pi \rightarrow \rho\vert\,\vert^{-s}\times \rho'\vert\,\vert^{-s'}\times \pi$ est holomorphe sur l'hyperplan $s'=s'_{0}$ par exemple, parce qu'il s'obtient en combinant le pr\'ec\'edent avec l'op\'erateur $\rho'\vert\,\vert^{s'}\times \pi' \rightarrow \rho'\vert\,\vert^{-s'}\times \pi'$ que l'on a d\'ej\`a \'etudi\'e (c'est le 3e ci-dessus);

par contre l'op\'erateur 
$\rho\vert\,\vert^{s}\times \rho'\vert\,\vert^{s'}\times \pi \rightarrow \rho'\vert\,\vert^{-s'}\times \rho\vert\,\vert^{-s}\times \pi$ est identiquement 0 sur $s'=s'_{0}$ si $t>1$ \`a cause de l'ordre de la fonction $r(s',\rho',\pi^{'GL})$ sur cet hyperplan. Par contre si $t=1$ cette fonction global est inversible et la fonction $r(s,\rho,\pi^{',GL})$ fourni alors un p\^ole d'ordre 1 en $s=s_{0}$. Pour \^etre s\^ur qu'il y ait bien un p\^ole, il faut que les op\'erateurs d'entrelacement locaux n'aient pas de z\'ero apr\`es normalisation par le produit des fonctions $r(s',\rho',\pi^{'GL})r(s,\rho,\pi^{'GL})$. Mais rien ne force ces op\'erateurs \`a avoir un z\'ero: un cas o\`u on n'a s\^urement pas de z\'ero est le cas o\`u les fonctions locales normalisant les op\'erateurs d'entrelacement sont holomorphes et pour cela il suffit que pour tout $(\rho'',b'')\in Jord(\pi^{'GL})$ on ait $s'_{0}>(b''-1)/2$ et ceci est possible car les seuls $(\rho'',b'')$ impos\'es sont $(\rho,2s_{0}-1)=(\rho,2s'_{0})$ et $(\rho',2s'_{0}-1)$.

Ainsi si on calcule les fonctions $(s-s_{0})E(\rho\vert\,\vert^{s}\times \rho'\vert\,\vert^{s'}\times \pi',f)$ (pour $f$ une section de l'induite) sur l'hyperplan $s'=s'_{0}$ puis en $s=s_{0}$, le r\'esultat est une forme automorphe dont le seul termes constant non nul est dans l'induite $\rho'\vert\,\vert^{-s'_{0}}\times \rho\vert\,\vert^{-s_{0}}\times \pi'$; elle est donc de carr\'e int\'egrable et n\'ecessairement irr\'eductible. 

La preuve de la remarque est termin\'e. On va encore montrer comment l'espace des r\'esidus, quand il est de carr\'e int\'egrable s'obtient autrement et de fa\c{c}on conforme \`a nos r\'esultats. Ici pour pourvoir utiliser \ref{residugeneral}, on suppose que  l'on sait a priori que les repr\'esentations de carr\'e int\'egrable irr\'eductible ayant de la cohomologie interviennent avec multiplicit\'e 1 dans le spectre discret.

\begin{rem}{\bf 2} On garde les hypoth\`eses et notations pr\'ec\'edentes et on note $\sigma$ la repr\'esentation de carr\'e int\'egrable d\'efinie par les r\'esidus de la remarque 1; on suppose que $\sigma$ est  non nulle. Alors la repr\'esentation $\pi^{'+}:=\bigl((s-s_{0})E(\rho\vert\,\vert^{s}\times \pi',s)\bigr)_{s=s_{0}}$ est de carr\'e int\'egrable et 
$$
\sigma=\biggl((s'-s'_{0})E(\rho'\vert\,\vert^{s'}\times \pi^{'+},s')\biggr)_{s'=s'_{0}}.
$$
\end{rem}
Le fait que $\pi^{'+}$ est de carr\'e int\'egrable est \'evident (cf. \ref{rang1}).
On a le diagramme commutatif d'op\'erateur d'entrelacement globaux:
$$\begin{matrix}
\rho\vert\,\vert^{s}\times \rho'\vert\,\vert^{s'}\times \pi' &\rightarrow& \rho'\vert\,\vert^{s'}\times \rho\vert\,\vert^{s}\times \pi'\\
\downarrow&&\downarrow\\
\rho'\vert\,\vert^{-s'}\times \rho\vert\,\vert^{-s}\times \pi' &\simeq &\rho'\vert\,\vert^{-s'}\times \rho\vert\,\vert^{-s}\times \pi'.
\end{matrix}
$$
On regarde en toute place pour garder la commutativit\'e, on ajoute les normalisations: pour la ligne du haut c'est la normalisation de Langlands-Shahidi et pour les colonnes ce sont celles que l'on a utilis\'ees ici. On sait que la fl\`eche de droite a une image nulle ou irr\'eductible et on identifie alors le produit tensoriel en toute place de ces images \`a 
$$
\biggl((s'-s'_{0})\biggl((s-s_{0})E(\rho'\times \rho\times \pi',s',s)_{s=s_{0}}\biggr)_{s'=s'_{0}}\biggr);
$$pour cela on utilise d'abord \ref{rang1} puis \ref{residugeneral} (i).
Avec la commutativit\'e du diagramme ci-dessus, l'espace des r\'esidus de la remarque est inclus dans celui que l'on vient d'\'ecrire; cela est conforme \`a \ref{residugeneral}(ii) et donne la remarque
\subsection{Commentaires sur les formes automorphes de carr\'e int\'egrable non cuspidales\label{commentaires}}
Notre but, expliqu\'e en \cite{manuscripta}, est d'obtenir des conditions n\'ecessaires et suffisantes sur les param\`etres d'Arthur pour savoir quand une repr\'esentation automorphe irr\'eductible de carr\'e int\'egrable n'est pas cuspidale. Le travail fait ici, dans le cas des repr\'esentations automorphes ayant de la cohomologie, cl\^ot presque la question et on va expliquer ce qui manque pour conclure. On admet ici que toute repr\'esentation automorphe ayant de la cohomologie intervient avec multiplicit\'e au plus 1 dans le spectre discret des groupes de m\^eme type que $G$.

Soit $\sigma$ une forme automorphe de carr\'e int\'egrable ayant de la cohomologie. On a montr\'e en \cite{manuscripta} que si $\sigma$ n'est pas cuspidale, il existe une repr\'esentation de carr\'e int\'egrable $\pi$ d'un groupe de m\^eme type que $G$, une repr\'esentation cuspidale autoduale $\rho$ d'un groupe $GL$ convenable et un demi-entier $s_{0}$ tel que 
$$
\sigma\hookrightarrow \biggl((s-s_{0})E(\rho\times \pi,s)_{s=s_{0}}\biggr).\eqno(1)
$$
On a en plus pr\'ecis\'e comment s'obtient $\pi^{GL}$ le param\`etre du paquet de $\pi$ \`a partir de $\sigma^{GL}$ le param\`etre du paquet de $\sigma$.

Avec les notations introduites dans cet article et en particulier en \ref{residugeneral} on sait que (1) est \'equivalent \`a ce que $\sigma=\pi^+$; on a s\^urement la non nullit\'e. Cela veut dire qu'en toute place $v$, $\sigma_{v}=\pi^+_{v}$. En \ref{parametregeneralisation} on a construit $\pi^+_{v}$ en fonction de $\pi_{v}$, ici on fait l'inverse, en ayant $\sigma_{v}$ on construit des param\`etres pour une repr\'esentation dans le paquet associ\'e \`a $\pi^{GL}$. Il y a une condition pour que ce soit possible: aux places archim\'ediennes, il faut que $\sigma_{v}$ soit un quotient de Langlands convenable et aux places $p$-adiques il faut que les param\`etres de $\sigma_{v}$ soit dans l'image de \ref{parametregeneralisation} et la non surjectivit\'e est celle de \ref{descriptionparametre}: cette r\'ef\'erence s'applique par \'etage comme expliqu\'e en \ref{parametregeneralisation} et c'est dans les cas o\`u il faut $t_{0}^+=t_{0}+1$ qu'il y a un d\'efaut de surjectivit\'e (c'est le m\^eme type de condition que dans le cas archim\'edien mais dans le cas archim\'edien les hypoth\`eses sont tellement simplificatrices que l'on ne voit pas la difficult\'e). M\^eme sans les avoir rendu totalement explicite on voit comment le fait que $\sigma_{v}$ doit \^etre de la forme $\pi^+_{v}$ donne des conditions n\'ecessaires et quand ces conditions sont remplies on a des param\`etres pour un unique \'el\'ement du paquet associ\'e \`a $\pi^{GL}$. On note cette \'el\'ement $\sigma^-_{v}$, ici elle peut \^etre 0 puisqu'elle n'est connu que par ses param\`etres. Et on doit n\'ecessairement avoir:
$$
\pi\simeq \otimes_{v}\sigma^-_{v}.\eqno(2)
$$

\begin{prop}
Pour que  $\sigma$ soit non cuspidal (on note $Jord(\sigma^{GL}))$ les couples $(\rho',b')$ tel que $\sigma^{GL}\simeq Speh(\rho',b')$), il faut et il suffit qu'il existe
$\rho$ une repr\'esentation cuspidale unitaire d'un certain $GL$ et $s_{0}$ un demi-entier tels que

(3) soit $s_{0}=1/2$ avec $L(\rho,r_{G},s)$ a un p\^ole en $s=1$, soit $Jord(\sigma^{GL})$ contient $(\rho,2s_{0}-1)$;

(4) pour tout $(\rho',b')\in Jord(\sigma^{GL})$ tel que $b'=2s_{0}$, $L(\rho\times \rho',1/2)\neq 0$;

(5) pour toute place $v$, $\sigma^-_{v}$ est d\'efini et est non nulle;

(6) la repr\'esentation $\sigma^-:=\otimes_{v}\sigma^-_{v}$ est une repr\'esentation automorphe de carr\'e int\'egrable.
\end{prop}
C'est un corollaire de \ref{residugeneral} mais ce n'est pas satisfaisant car je pense que la condition (6) est inutile et plus exactement qu'elle est cons\'equence de la construction de $\sigma^-$ \`a partir de $\sigma$; si on a en t\^ete la formule de multiplicit\'e annonc\'ee par Arthur, la condition de signe qui doit \^etre v\'erifi\'ee par $\sigma^-$ r\'esulte \`a mon avis de celle v\'erifi\'ee n\'ecessairement par $\sigma$. C'est l'int\'er\^et d'avoir des param\`etres, car le signe doit pouvoir se calculer en fonction des param\`etres. C'est ce qui reste \`a faire.

\section{Appendice}
\subsection{Propri\'et\'es des modules de Jacquet des repr\'esentations dans un paquet d'Arthur\label{proprietedujac}}

Soit $\psi$ un param\`etre pour un paquet d'Arthur et soit $\pi$ dans le paquet associ\'e. On reprend les notations $(\rho',A',B',\zeta')$ pour d\'ecrire $Jord(\psi)$. Soit $[x,y]$ un segment et on fixe $\rho$. On suppose que $\psi$ a bonne parit\'e (pour simplifier)
\begin{lem}
On suppose que $Jac_{x,\cdots,y}\pi\neq 0$, alors il existe un ensemble fini $[1,v]$ et pour tout $i\in [1,v]$, $(\rho,A_{i},B_{i},\zeta_{i})\in Jord(\psi)$ tel que $\zeta_{1}B_{1}=x$, $A_{v}\geq \vert y \vert$ et pour tout $i\in [1,v[$, $B_{i+1}\leq A_{i}+1$.
\end{lem}
On a montr\'e ce lemme en \cite{holomorphie} 2.7 sous l'hypoth\`ese suppl\'ementaire $\vert y\vert \geq \vert x\vert$ et $x\neq 0$. C'est en fait le cas le plus difficile et on va montrer les cas restants. Les cas restants sont les cas o\`u $0\in [x,y]$: en effet si $x\neq 0$, le cas non trait\'e est celui o\`u $\vert y\vert <\vert x\vert$. Supposons qu'il en soit ainsi; le simple fait que $Jac_{x}\pi\neq 0$ n\'ecessite qu'il existe $(\rho,A,B,\zeta)$ avec $\zeta B=x$. Et clairement $v=1$ et $(\rho,A_{1},B_{1},\zeta_{1})=(\rho,A,B,\zeta)$ convient. On suppose donc dans ce qui suit que $0\in [x,y]$ et que $\vert y\vert >\vert x\vert$. Remarquons qu'il suffit de montrer l'assertion suivante: il existe $(\rho,A,B,\zeta)$ tel que $\zeta B \in [x,y]$ et $A\geq \vert y \vert$: en effet on sait d\'ej\`a qu'il existe $(\rho,A_{1},B_{1},\zeta_{1})$ avec $\zeta_{1}B_{1}=x$ et on fixe $A_{1}$ maximum avec cette propri\'et\'e. Si $A_{1}\geq \vert y\vert$, l'assertion du lemme est prouv\'ee; sinon on note $\zeta$ le signe de $y$ et on remarque que $Jac_{x, \cdots, \zeta (A_{1}+1) }\pi\neq 0$. On applique l'assertion momentan\'ement admise qui assure l'existe de $(\rho,A_{2},B_{2},\zeta_{2})\in Jord(\psi)$ avec $\zeta_{2}B_{2}\in [x,\zeta (A_{1}+1)]$ et $A_{2}\geq (A_{1}+1)$; on fixe un tel choix avec $A_{2}$ maximal. On remarque que $B_{2} \leq sup( \vert x\vert, A_{1}+1)= A_{1}+1$. Si $A_{2}\geq \vert y\vert$, le lemme est d\'emontr\'e et sinon on continue: supposons que l'on ait construit une suite  d'\'el\'ements de $Jord(\psi)$ ind\'ex\'ee par $[1,v]$ v\'erifiant toutes les propri\'et\'es du lemme sauf $A_{v}\geq \vert y\vert$. On a donc certainement encore $Jac_{x, \cdots, \zeta (A_{v}+1)}\pi\neq 0$; l'assertion admise donne l'existence de $(\rho,A_{v+1},B_{v+1},\zeta_{v+1})$ avec $\zeta_{v+1} B_{v+1}\in [x, \zeta(A_{v}+1)]$ et $A_{v+1} \geq A_{v}+1$; d'o\`u encore $B_{v+1}\leq A_{v}+1$ et on la suite d'\'el\'ements index\'ee par $[1,v+1]$ v\'erifie encore toutes les conditions du lemme sauf \'eventuellement $A_{v+1}\geq \vert y\vert$. En un nombre fini d'\'etapes, on obtient aussi cette derni\`ere condition. Il reste \`a montrer l'assertion admise. On la montre d'abord pour les morphismes $\psi$ ayant la propri\'et\'e suivante: soit $(\rho',A',B',\zeta')\in Jord(\psi)$; si $\rho'\simeq \rho$, soit $\zeta' B' \in [x,y]$ soit $\vert y\vert << B'$. Sous ces hypoth\`eses, on reprend les d\'efinitions pour construire $\pi$ \`a partir d'un paquet de repr\'esentations associ\'e \`a des morphismes de restriction discr\`etes \`a la diagonale de $SL(2,{\mathbb C})\times SL(2,{\mathbb C})$; on doit faire ''descendre'' pour arriver aux \'el\'ements de $Jord(\psi)$ de la forme $(\rho,A',B',\zeta')$ avec $\zeta' B'\in [x,y]$ car on peut s'arranger pour que ces \'el\'ements soient plus petit (dans l'ordre fix\'e) que ceux ne v\'erifiant pas cette propri\'et\'e (prendre par exemple l'ordre sur les ''B''); on v\'erifie alors qu'il existe

(1) un ensemble totalement ordonn\'e de demi-entiers relatifs, ${\mathcal E}$, tel que pour $c\in {\mathcal E}$ il existe $(\rho,A',B',\zeta')\in Jord(\psi)$ avec $\zeta' B' \in [x,y]$ et $\vert c\vert \leq A'$

(2) une repr\'esentation irr\'eductible $\sigma$ d'un groupe de m\^eme type que $G$ v\'erifiant $Jac_{z}\sigma\neq 0$ seulement s'il existe $(\rho,A',B',\zeta')\in Jord(\psi)$ avec $\zeta' B'\notin [x,y]$ et $z=\zeta' B'$ et en particulier $z>>\vert y\vert$;

(3) et une inclusion $\pi\hookrightarrow \times_{c\in {\mathcal E}}\rho\vert\,\vert^c\times \sigma$.

La non nullit\'e de $Jac_{x,\cdots, y}\pi$ entra\^{\i}ne la non nullit\'e du m\^eme module de Jacquet pour l'induite que l'on vient d'\'ecrire en (3). Ainsi, avec (2), il existe $c\in {\mathcal E}$ avec $\vert c\vert=\vert y \vert$. En appliquant (1), on obtient l'assertion cherch\'ee.

On enl\`eve maintenant l'hypoth\`ese sur $\psi$; en d'autres termes il faut consid\'erer la situation suivante; soit $(\rho,A',B',\zeta')\in Jord(\psi)$ et supposons qu'il existe $\psi'$ dominant $\psi$ tel que $Jord(\psi')$ diff\`ere de $Jord(\psi)$ simplement en rempla\c{c}ant $(\rho,A',B',\zeta')\in Jord(\psi)$ par $(\rho,A'+1,B'+1,\zeta')\in Jord(\psi')$ et soit $\pi'$ dans le paquet associ\'e \`a $\psi'$ tel que
$$
\pi'\hookrightarrow \rho'\vert\,\vert^{\zeta'(B'+1)}\times \cdots \rho'\vert\,\vert^{\zeta'(A'+1)}\times \pi.
$$
Dans cette situation, on suppose que le  lemme est connu pour $\pi'$ et que $\zeta' B'\notin [x,y]$ (en particulier $B'\neq 0$) et il faut d\'emontrer le lemme pour $\pi$. Par hypoth\`ese $Jac_{x, \cdots, y}\pi\neq 0$. V\'erifions que $[\zeta' B', \zeta' A']\cap [x,y]=\emptyset$: en effet si cette intersection est non vide, il y a a priori 2 possibilit\'e, l'un des segments est inclus dans l'autre ou l'intersection  est un segment qui contient une des extr\'emit\'es de $[\zeta' B',\zeta' A']$ et une extr\'emit\'e de $[x,y]$. Le segment $[x,y]$ contient $0$ et ne peut donc \^etre inclus dans $[\zeta' B', \zeta' A']$, le segment $[\zeta' B', \zeta' A']$ n'est pas inclus dans $[x,y]$ car $\zeta' B'$ n'y est pas. Donc il reste la possibilit\'e o\`u soit $\zeta' B'$ soit $\zeta' A'$ est dans $[x,y]$; on sait que $\zeta' B' \notin [x,y]$ par hypoth\`ese; on a donc soit $B'\geq \vert y\vert$ et a fortiori $\zeta' A' \notin [x,y]$ soit $\zeta' B'$ est de signe oppos\'e \`a $y$ avec $B' > \vert x\vert$. Dans ce dernier cas $\zeta' A'$ est aussi de signe oppos\'e \`a $y$ et a le signe de $x$ en v\'erifiant aussi $A'>\vert x\vert$. On sait alors que pour tout $z \in [x,y]$ et pour tout $C\in [\zeta' (B'+1),\zeta' (A'+1)]$ l'induite $\rho\vert\,\vert^{z}\times \rho\vert\,\vert^{C}$ pour le GL convenable est irr\'eductible. D'o\`u aussi (avec (*) appliqu\'e \`a $\pi$ et \`a $\pi'$) $Jac_{x, \cdots, y} \pi'\neq 0$. Le r\'esultat connu pour $\pi'$ montre l'existe de $(\rho,A'',B'',\zeta'')\in Jord(\psi')$ avec $\zeta'' B'' \in [x,y]$ et $A''\geq \vert y \vert$; n\'ecessairement $(\rho,A'',B'',\zeta'')\in Jord(\psi)$ et on a le r\'esultat pour $\pi$. Cela termine la preuve.

\subsection{Propri\'et\'e d'irr\'eductibilit\'e\label{irreductibilite}}
Soit $\psi$ comme pr\'ec\'edemment et soit $\pi$ dans le paquet associ\'e \`a $\psi$ et soit aussi $\rho$ une repr\'esentation cuspidale unitaire d'un groupe GL et $x\in {\mathbb R}$. Un probl\`eme important est de savoir quand l'induite $\rho\vert\,\vert^{x}\times \pi$ est r\'eductible. Je n'ai pas la r\'eponse en g\'en\'eral; on a besoin ici du cas suivant:
\begin{lem} On suppose que $\rho$ est autoduale que $\psi$ est de bonne parit\'e et que $x\neq 0$ v\'erifie: pour tout $(\rho,A,B,\zeta)\in Jord(\psi)$ soit $A<\vert x\vert -1$ soit $B>\vert x\vert $. Alors pour tout $\pi$ dans le paquet associ\'e \`a $\psi$, l'induite $\rho\vert\,\vert^{x}\times \pi$ est irr\'eductible.
\end{lem}
On fixe $\psi_{>}$ un morphisme dominant $\psi$ et tel que $Jord(\psi_{>})$ contient tous les \'el\'ements $(\rho,A,B,\zeta)$ de $Jord(\psi)$ tel que $A<\vert x\vert -1$ et par contre pour tout $(\rho,A,B,\zeta)\in Jord(\psi)$ tel que $B>\vert x\vert$ l'image correspondant dans $Jord(\psi_{>})$, \'ecrite $(\rho,A+T_{\rho,A,B,\zeta},B+T_{\rho,A,B,\zeta},\zeta)$ (cf. \ref{descriptionpaquets}) v\'erifie $T_{\rho,A,B,\zeta}>>0$. On fixe $\pi$ dans le paquet associ\'e \`a $\psi$ et on note $\pi_{>}$ la repr\'esentation dans le paquet associ\'e \`a $\psi_{>}$ qui correspond \`a $\pi$. On remarque que $\rho\vert\,\vert^{x}\times \pi$ et $\rho\vert\,\vert^{x}\times \pi_{>}$ ont un unique sous-module irr\'eductible not\'e $\sigma$ et $\sigma_{>}$: ceci r\'esulte de ce que $x\neq 0$ et  $Jac_{x}\pi=Jac_{x}\pi_{>}=0$; cette derni\`ere assertion r\'esulte de  \ref{proprietedujac} appliqu\'e \`a $y=x$ et du fait que pour tout $(\rho,A',B',\zeta')\in Jord(\psi)$ $\zeta'B'\neq x$ par hypoth\`ese. Il en est de m\^eme des induites $\rho\vert\,\vert^{-x}\times \pi$ et $\rho\vert\,\vert^{-x}\times \pi_{>}$. Pour d\'emontrer l'irr\'eductibilit\'e, il suffit donc de montrer que $Jac_{-x}\sigma\neq 0$ et $Jac_{-x}\sigma_{>}\neq 0$: en effet $\sigma$ et $\sigma_{>}$ sont isomorphes \`a des quotients de $\rho\vert\,\vert^{-x}\times \pi$ et $\rho\vert\,\vert^{-x}\times \pi_{>}$ et par un simple calcul de module de Jacquet n'intervienent dans ces induites comme sous-quotient qu'avec multiplicit\'e 1. Admettons momentan\'ement que $Jac_{-x}\sigma_{>}\neq 0$ et montrons que $Jac_{x}\sigma\neq 0$: en effet on passe de $\pi_{>}$ \`a $\pi$ en prenant des $Jac_{\zeta' C, \cdots ,\zeta' D}$ pour des demi-entiers $C\leq D$ qui v\'erifie tous qu'il existe $(\rho,A',B',\zeta')\in Jord( \psi)$ avec $B'>\vert x\vert$ et $C\geq B'+1$. Donc en particulier pour tout $C'\in [C,D]$ l'induite $\rho\vert\,\vert^{\pm x}\times \rho\vert\,\vert^{C'}$ est irr\'eductible dans le GL convenable. Ceci montre que $\sigma$ s'obtient \`a partir de $\sigma_{>}$ en prenant les m\^emes modules de Jacquet; de plus s'il existe $\sigma'_{>}$ avec une inclusion
$\sigma_{>}\hookrightarrow \rho\vert\,\vert^{-x}\times \sigma'_{>}$ on a alors aussi par exactitude des modules de Jacquet l'existence d'une repr\'esentation $\sigma'$, tel que  $\sigma\hookrightarrow \rho\vert\,\vert^{-x}\times \sigma$. Ainsi   l'irr\'eductibilit\'e de $\rho\vert\,\vert^{-x}\times \pi$ r\'esulte de l'irr\'eductibilit\'e de $\rho\vert\,\vert^{x}\times \pi_{>}$.

Montrons l'irr\'eductibilit\'e de $\rho\vert\,\vert^{x}\times \pi_{>}$; on peut \'evidemment supposer que $x>0$. On sait alors que l'op\'erateur d'entrelacement normalis\'e $\rho\vert\,\vert^{x}\times \pi \rightarrow \rho\vert\,\vert^{-x}\times \pi$ est holomorphe (\cite{holomorphie}). Le facteur de normalisation n'a pas de p\^ole: avec les notations de \ref{poles}, $a_{0}=1$, $b_{0}=2x+1$ d'o\`u $A_{0}=x=B_{0}$ et $\zeta_{0}=-$. Ainsi l'op\'erateur d'entrelacement normalis\'e co\"{\i}ncide \`a un scalaire non nul pr\`es \`a l'op\'erateur d'entrelacement standard qui est donc lui aussi holomorphe. Il suffit donc de d\'emontrer que cet op\'erateur d'entrelacement standard est injectif. Pour avoir ce r\'esultat, on montre qu'il existe 

un morphisme $\psi'$ tel que 
la restriction de $\psi'$ \`a $W_{F}$ fois la diagonale de $SL(2,{\mathbb C})\times SL(2,{\mathbb C})$ est sans multiplicit\'e
$Jord(\psi')$ contient tous les \'el\'ements $(\rho,\tilde{A},\tilde{B},\zeta)$ avec $\tilde{B}>\vert x\vert$ de $Jord(\psi_{>})$ et d'autres \'el\'ements $(\rho,C',C'',\zeta''')$ v\'erifiant $C'<\vert x\vert-1$

une repr\'esentation $\pi'$ dans le paquet associ\'e \`a $\psi'$ et un ensemble ${\mathcal E}$ comme dans la preuve de \ref{proprietedujac} avec pour tout $z\in {\mathcal E}$, $\vert z\vert<x-1$

et une inclusion $
\pi_{>}\hookrightarrow \times_{z\in {\mathcal E}}\rho\vert\,\vert^{z}\times \pi'.
$

La repr\'esentation $\rho\vert\,\vert^{x}\times \pi'$ est irr\'eductible parce que les hypoth\`eses de \cite{discret} 6.4 sont satisfaites; les induites $\rho\vert\,\vert^z\times \rho\vert\,\vert^{\pm x}$ sont irr\'eductibles pour tout $z\in {\mathcal E}$ dans le GL convenable car $\vert x-z\vert >1$. Donc l'entrelacement d\'ecrit appliqu\'e au membre de droite de l'inclusion ci-dessus est un isomorphisme et sa restriction \`a $\rho\vert\,\vert^{x}\times \pi$ est injective. Ainsi $\sigma_{>}$ est un sous-module de $\rho\vert\,\vert^{-x}\times \pi_{>}$ et c'est ce que l'on cherchait.


\begin{thebibliography}{}
\bibitem{aj} J. Adams, J. Johnson \sl Endoscopic groups and packets of non-tempered representations \rm Compositio Math., 64, 1987, pp. 271-309
\bibitem{arthurnouveau} J. Arthur, \textit{An introduction to the trace formula} Clay Mathematics Proceedings, volume 4, 2005, pp. 1-253
\bibitem{bc} N. Bergeron, L. Clozel, \sl  Spectre automorphe des vari\'et\'es hyperboliques et applications topologiques \rm Ast\'erisque, 303, 2005, 214 pages
\bibitem{franke98} J Franke {\sl Harmonic analysis in weighted $L^2$-spaces} Ann. Sci. de l'ENS (4), 31, 1998, pp. 181-279
\bibitem{franke} J. Franke, {\sl A topological model for some summand of the Eisenstein cohomology of congruence subgroups}, pr\'epublication 1991, paru in Eisenstein series and applications
\'edit\'e par W. T. Gan, S. Kudla, Y. Tschinkel, Progress in Math 258, Birkh\"auser, 2007
\bibitem{g} G. Gotsbacher {\sl Eisenstein cohomology for congruence subgroups of SO(n,2)}, pr\'epublication 2008, archiv: 0803.0762
\bibitem{hg} H. Grobner {\sl The automorphic cohomology and the residual spectrum of Hermintian groups of rank one}, pr\'epublication 2009
\bibitem{gg} G. Gotsbacher, H. Grobner, \sl On the Eisenstein cohomology of odd orthogonal groups, \rm pr\'epublication 2009, arXiv:0904.2562v111
\bibitem{gs} N. Grbac, J. Schwermer, \sl On residual cohomology classes attached to relative rank one Eisenstein series for the symplectic group, \rm pr\'epublication 2009
\bibitem{harris} M. Harris, R. Taylor, \sl The geometry and
cohomology of some simple Shimura varieties,\rm  Annals of Math
Studies, 151, Princeton Univ. Press, 2001
\bibitem{henniart} G. Henniart, \sl Une preuve simple des
conjectures de Langlands pour $GL_n$ sur un corps p-adique,\rm
Invent. Math., 139, 2000, pp. 439-455
\bibitem{jpss} H. Jacquet, I. I. Piatetskii-Shapiro, J. A. Shalika, \sl Rankin-Selberg Convolutions, \rm American Journal of Mathematics, Vol. 105, No. 2, 1983, pp. 367-464 
\bibitem{johnson} J. Johnson \sl Stable base change $\mathbb{C}/\mathbb{R}$  of certain derived functor modules \rm Math. Ann., 287, 1990, pp. 467-493
\bibitem{kasten} H. Kasten {\sl Cohomological Representations and Twisted Rankin-Selberg 
Convolutions},  Dissertation,  Univ Karlsruhe, 2007
\bibitem{elementaire} C. M{\oe}glin,  \textit {Sur certains paquets d'Arthur et involution
d'Aubert-Schneider-Stuhler g\'en\'eralis\'ee
}   ERT, volume 10, 2006, pp. 86-129
\bibitem{unitaire} C. M{\oe}glin \sl Classification et Changement de base pour les  
s\'eries discr\`etes des groupes unitaires 
p-adiques, \rm Pacific Journal of Math, 233-1,  2007, pp 159-204
\bibitem{discret} C. M{\oe}glin, \textit {
Paquets d'Arthur discrets pour un groupe classique p-adique}
 pr\'epublication 2004, \`a para\^{\i}tre dans le volume en l'honneur de S. Gelbart, AMS
\bibitem{pourshahidi} C. M{\oe}glin, \textit{Multiplicit\'e 1 dans les paquets d'Arthur aux places p-adiques} pr\'epublication 2007, \`a para\^itre Volume en l'honneur de F. Shahidi
\bibitem{manuscripta} C. M{\oe}glin, \textit{Formes automorphes de carr\'e int\'egrable non cuspidales} 
Manuscripta Math, 127, 2008, pp. 411-467
\bibitem{holomorphie}C. M{\oe}glin, \sl Holomorphie des op\'erateurs d'entrelacement normalis\'es \`a l'aide des param\`etres d'Arthur, \rm \`a para\^{\i}tre au Canadian Journal of Math
\bibitem{mvw} C. M{\oe}glin, M.-F. Vign\'eras, J.-L. Waldspurger \sl Correspondance de Howe sur un corps p-adique, \rm LN 129, 1987, Springer Verlag
\bibitem{bible} C. M{\oe}glin, J. -L. Waldspurger, \sl D\'ecomposition spectrale et s\'eries d'Eisenstein; une paraphrase de l'Ecriture \rm Birkh\"auser, PM 113, 1994
\bibitem{mw}C. M{\oe}glin, J.-L. Waldspurger, \textit {Le spectre r\'esiduel de GL(n)}  Ann. de l'ENS, 22, 1989, pp. 605-674
\bibitem{shahidi} F. Shahidi, \sl Local coefficients and normalization of intertwining operators for GL(n) \rm Comp. Math., 48, 1983, pp. 271-295
\bibitem{shahididuke} F. Shahidi, \sl Local coefficients as Artin factors for real groups, \rm Duke Math. J., 52, 1985, pp. 973-1007
\bibitem{vogan} D. Vogan, \sl Representations of real reductive Lie groups, \rm PM 15, Birkh\"auser, 1981
\bibitem{vz} D. Vogan, G. Zuckerman, \sl Unitary representations with non zero cohomology \rm Compositio Math 53, 1984, pp. 51-90
\end{thebibliography}
 \end{document}